\documentclass[a4paper,11pt,onecolumn,twoside]{article}
\usepackage{mathrsfs}
\usepackage{mathrsfs}
\usepackage{amsfonts}
\usepackage{fancyhdr}
\usepackage{amsmath,amsfonts,amssymb,graphicx,amsthm}
\begin{document}

\title{\bf Existence and stability of stationary solutions to the full compressible Navier-Stokes-Korteweg system }
\author{{\bf Zhengzheng Chen}\\
School of Mathematics and Statistics\\
Wuhan University, Wuhan 430072, China\\[2mm]
{\bf Huijiang Zhao}\thanks{Corresponding author.E-mail:
hhjjzhao@hotmail.com}\\
School of Mathematics and Statistics, \\Wuhan University, Wuhan
430072, China }
\date{}

\vskip 0.2cm

\maketitle

\vskip 0.2cm \arraycolsep1.5pt
\newtheorem{Lemma}{Lemma}[section]
\newtheorem{Theorem}{Theorem}[section]
\newtheorem{Definition}{Definition}[section]
\newtheorem{Proposition}{Proposition}[section]
\newtheorem{Remark}{Remark}[section]
\newtheorem{Corollary}{Corollary}[section]

\begin{abstract}
This paper is concerned with the existence, uniqueness and nonlinear stability of stationary solutions to the Cauchy problem of the full compressible Navier-Stokes-Korteweg system effected by external force of general form in $\mathbb{R}^3$. Based on the weighted-$L^2$ method and some elaborate $L^\infty$ estimates of solutions to the linearized problem, the existence and uniqueness of stationary solution are obtained by the contraction mapping principle. The proof of the stability result is given by an elementary energy method and relies on some intrinsic properties of the full compressible Navier-Stokes-Korteweg system.

 \bigbreak
\noindent


\end{abstract}

\section{Introduction }
\setcounter{equation}{0}
In this paper, we are interested in the following nonisothermal compressible fluid models of Korteweg type, which can be derived from a Cahn-Hilliard like free energy( see the pioneering work by  Dunn and Serrin \cite{J. E. Dunn-J. Serrin-1985}, and also \cite{D. M. Anderson-G. B. McFadden-G. B.
Wheeler-1998, J. W. Cahn-J. E. Hilliard-1998, M. E. Gurtin-D. Polignone-J. Vinals-1996}).
\begin{eqnarray}\label{1.1}
\left\{\begin{array}{ll}
\rho_t+\nabla\cdot(\rho v)=G(x),\\[2mm]
(\rho v)_t+\nabla\cdot(\rho v\bigotimes v)=\nabla\cdot(S+K)+\rho F(x),\\[2mm]
\displaystyle\left[\rho\left(e+\frac{v^2}{2}\right)\right]_t+\nabla\cdot\left[\rho v\left(e+\frac{v^2}{2}\right)\right]=\nabla\cdot(\tilde{\alpha}\nabla\theta)+\nabla\cdot\left((S+K)\cdot v\right)+\rho v\cdot F(x)+H(x),
 \end{array}\right.
\end{eqnarray}
Here $(x, t)\in\mathbb{R}^3\times\mathbb{R}^+$, $\rho>0, v=(v_1, v_2, v_3), \theta>0$ and $e$ denote the density, the velocity, the internal energy and the temperature of the fluids respectively. $\tilde{\alpha}$ is the heat conduction coefficient. $F(x)=(F_1(x), F_2(x), F_3(x)), G(x), H(x)$ are
the given external force, mass source and energy source, respectively. The viscous stress tensor $S$ and the Korteweg stress tensor $K$ are given by \begin{eqnarray}\label{1.2}
\left\{\begin{array}{ll}
S_{i,j}=(\mu^\prime\nabla\cdot v-P(\rho, e))\delta_{ij}+2\mu d_{ij}(v)\\[2mm]
K_{i,j}=\displaystyle\frac{\kappa}{2}(\Delta\rho^2-|\nabla\rho|^2)\delta_{ij}-\kappa\partial_i\rho\partial_j\rho,
 \end{array}\right.
\end{eqnarray}
where $d_{ij}(v)=(\partial_iv_j+\partial_jv_i)/2$ is the strain tensor, $P$ is the pressure, $\mu$ and $\mu^\prime$ are the viscosity coefficients, and $\kappa$ is the capillary coefficient. Notice that when $\kappa=0$, system (\ref{1.1}) is reduced to the compressible Navier-Stokes system. In this paper, we consider the case of $e=C_\triangledown\theta$, where $C_\triangledown$ is the heat capacity at the constant volume. Our basic assumptions are as follows: $\bar{\rho}, \bar{\theta}, \kappa, \mu, \mu^\prime $ and $\tilde{\alpha}$ are the constants satisfying $\bar{\rho}, \bar{\theta}, \kappa, \mu, \tilde{\alpha}>0$ and $\frac{2}{3}\mu+\mu^\prime\geq0$; $C_\triangledown>0$ is a constant and $P=P(\rho, \theta)>0$ is a smooth function of $\rho, \theta>0$ satisfying $P_\rho(\rho,\theta), P_\theta(\rho,\theta)>0$.

The main purpose of this manuscript is to study the nonlinear stability of stationary solutions to the Cauchy problem of the compressible Navier-Stokes-Korteweg system (1.1). It is convenient to study the Cauchy problem for the following form which is equivalent to (\ref{1.1}) for classical solutions,
\begin{eqnarray}\label{1.3}
\left\{\begin{array}{ll}
\rho_t+\nabla\cdot(\rho v)=G(x),\\[2mm]
\displaystyle \rho v_t+\rho(v\cdot\nabla)v=\displaystyle\mu\Delta v+(\mu+\mu^\prime)\nabla(\nabla\cdot v)-\nabla P(\rho, \theta)+\kappa\rho\nabla\Delta\rho+\rho F(x)- vG(x),\\[2mm]
\displaystyle\rho C_\triangledown\left(\theta_t+(v\cdot\nabla)\theta\right)+\theta P_\theta(\rho, \theta)\nabla\cdot v=\tilde{\alpha}\Delta\theta+\Psi(v)+\Phi(\rho,v)+H(x)+\frac{v^2}{2}G(x)-C_\triangledown G(x)\theta,
 \end{array}\right.
\end{eqnarray}
with the initial date
\begin{equation}\label{1.4}
(\rho, v, \theta)(t,x)|_{t=0}=(\rho_0, v_0, \theta_0)(x)\rightarrow(\bar{\rho}, 0, \bar{\theta}) \quad as\,|x|\rightarrow+\infty.
\end{equation}
Here
\begin{eqnarray}\label{1.5}
\left\{\begin{array}{ll}\Psi(v)=\mu^\prime(\nabla\cdot v)^2+2\mu\mathbb{D}v:\mathbb{D}v,\,\, \mathbb{D}v=(d_{ij}(v))_{i,j=1}^3, \\[2mm]
\Phi(\rho,v)=\kappa\left(\frac{|\nabla\rho|^2}{2}+\rho\Delta\rho\right)\nabla\cdot v-\kappa(\nabla\rho\bigotimes\nabla\rho):\nabla v
\end{array}\right.
\end{eqnarray}
The stationary problem corresponding to the initial value problem (\ref{1.3}), (\ref{1.4}) is
\begin{eqnarray}\label{1.6}
\left\{\begin{array}{ll}
\nabla\cdot(\rho v)=G(x),\\[2mm]
\displaystyle (v\cdot\nabla)v=\displaystyle\frac{1}{\rho}\left\{\mu\Delta v+(\mu+\mu^\prime)\nabla(\nabla\cdot v)-\nabla P(\rho, \theta)\right\}+\kappa\nabla\Delta\rho+ F(x)- \frac{v}{\rho}G(x),\\[3mm]
\displaystyle(v\cdot\nabla)\theta+\frac{\theta P_\theta(\rho, \theta)}{\rho C_\triangledown}\nabla\cdot v=\frac{1}{\rho C_\triangledown}\left\{\tilde{\alpha}\Delta\theta+\Psi(v)+\Phi(\rho,v)+H(x)+\frac{v^2}{2}G(x)-C_\triangledown G(x)\theta\right\},
 \end{array}\right.
\end{eqnarray}

Before stating our main results, we explain some notations as follows, which are borrowed from \cite{Y. Shibata-K. Tanaka-2003} and \cite{J. Z. Qian-H. Yin-2007}.

{\bf Notations:} Throughout this paper, we use the standard notation in vector analysis. For example, we put for scalar $u$, vectors $v=(v_1,v_2,v_3), w=(w_1, w_2, w_3)$ and matrix $f=(f_{ij})_{1\leq i,j\leq3}$.
\[\Delta u=\sum_{i=1}^{3}\frac{\partial^2u}{\partial x_i^{2}},\quad \Delta v=(\Delta v_1,\Delta v_2,\Delta v_3), \quad (v\cdot\nabla)u=\sum_{i=1}^{3}v_i\frac{\partial u}{\partial x_i},\]
\[(v\cdot\nabla)w=\left((v\cdot\nabla)w_1, (v\cdot\nabla)w_2, (v\cdot\nabla)w_3\right),\]
\[\displaystyle\nabla^ku=\left\{\partial_x^\alpha u||\alpha|=k\right\},\quad\nabla^kv=\left\{\partial_x^\alpha v_i||\alpha|=k, i=1,2,3\right\},\]
\[\nabla\cdot v=\sum_{i=1}^{3}\frac{\partial v_i}{\partial x_i},\quad\nabla \cdot f=\left(\sum_{j=1}^{3}\frac{\partial f_{1j}}{\partial x_j}, \sum_{j=1}^{3}\frac{\partial f_{2j}}{\partial x_j}, \sum_{j=1}^{3}\frac{\partial f_{3j}}{\partial x_j}\right)\]
Here $\alpha=(\alpha_1, \alpha_2, \alpha_3)$ is a multi-index, $|\alpha|=\alpha_1+\alpha_2+\alpha_3$ and $\partial_x^\alpha=\partial^{|\alpha|}/\partial x_1^{\alpha_1}\partial x_2^{\alpha_2}\partial x_3^{\alpha_3}$. Moreover, we also use the notations
\[\nabla u=\left(\frac{{\partial u}}{\partial x_1}, \frac{{\partial u}}{\partial x_2}, \frac{{\partial u}}{\partial x_3}\right),\quad\nabla u=\left(\frac{\partial^2 u}{\partial x_{j}\partial x_{i}}\right)_{1\leq i,\,j\leq3}, \]
and denote $\frac{\partial u}{\partial x_{i}}$ by $\partial_iu$ or $u_{x_i}$ without any confusion.

 Next, we introduce some function spaces. Let $L^p$ denote the usual $L^p$ space, put for scalars $u_1, u_2$ and vectors $v=(v_1, v_2, \cdots, v_n)$, $w=(w_1, w_2,  \cdots, w_n)$,
 \[\|u_1\|_{L^p}=\left(\int_{\mathbb{R}^3}|u_1(x)|^pdx\right)^{\frac{1}{p}}, \quad \|v\|_{L^p}=\left(\sum_{i=1}^n\|v_i\|^p_{L^p}\right)^{\frac{1}{p}}, \quad(1\leq p<\infty),\]
 \[\|u_1\|_{L^\infty}=\sup_{\mathbb{R}^3}|u_1(x)|, \quad \|v\|_{L^\infty}=\max_{1\leq i\leq n}\|v_i(x)\|_{L^\infty}, \quad \langle u_1, u_2\rangle=\int_{\mathbb{R}^3}u_1u_2dx,\]
\[\langle v, w\rangle=\sum_{i=1}^{n}\langle v_i, w_i\rangle, \quad \|v\|_k=\left(\sum_{0\leq l\leq k}\|\nabla^lv\|^2\right)^{\frac{1}{2}}\, with \,\|\cdot\|=\|\cdot\|_{L^2}\]
\[H^k=\{u\in L^{1}_{loc}|\,\|u\|_k<\infty\}, \quad \hat{H}^k=\{u\in L^{1}_{loc}|\,\nabla u\in H^{k-1}\}, \]
where $u$ is either a vector or scalar. Further we put
\[\mathcal{H}^{k,l}=\big\{(\sigma,v)|\,\sigma\in H^k,\,\, v\in H^l\big\},\qquad \hat{\mathcal{H}}^{k,l}=\big\{(\sigma,v)|\,\sigma\in \hat{H}^k,\,\, v\in \hat{H}^l\big\}, \]
\[\mathcal{H}^{j,k,l}=\big\{(\sigma, v, \vartheta)|\,\sigma\in H^j,\,\, v\in H^k, \vartheta\in H^l\big\}, \]
\[\hat{\mathcal{H}}^{j,k,l}=\big\{(\sigma, v, \vartheta)|\,\sigma\in \hat{H}^j,\,\, v\in \hat{H}^k, \vartheta\in \hat{H}^l\big\}, \]
and
\[\|(\sigma, v)\|=\|\sigma\|_k+\|v\|_l, \qquad\|(\sigma, v, \vartheta)\|_{j,k,l}=\|\sigma\|_j+\|v\|_k+\|\vartheta\|_l. \]
\begin{Definition}
\[I_\epsilon^k=\displaystyle\big\{\sigma\in H^k|\,\|\sigma\|_{I^k}<\epsilon\big\}, \quad J_\epsilon^k=\displaystyle\big\{v\in H^k|\,\|v\|_{J^k}<\epsilon\big\}, \quad N_\epsilon^k=\displaystyle\big\{\vartheta\in H^k|\,\|\vartheta\|_{N^k}<\epsilon\big\}, \]
where
\[\|\sigma\|_{I^k}=\|\sigma\|_{L^6}+\sum_{\nu=1}^k\left\|(1+|x|)^{\nu}(\nabla^\nu\sigma, \nabla^{\nu+1}\sigma, \nabla^{\nu+2}\sigma)\right\|+\left\|(1+|x|)^{2}\sigma\right\|_{L^\infty}+\left\|(1+|x|)^{2}\nabla\sigma\right\|_{L^\infty}, \]
\[\|v\|_{J^k}=\|v\|_{\hat{J}^k}+\sum_{\nu=1}^k\left\|(1+|x|)^{\nu-1}\nabla^\nu v\right\|,\quad \|\vartheta\|_{N^k}=\|\vartheta\|_{\hat{J}^k}+\sum_{\nu=1}^k\left\|(1+|x|)^{\nu-1}(\nabla^\nu \vartheta, \nabla^{\nu+1}\vartheta)\right\|,\]
and $\|\cdot\|_{\hat{J}^k}$ is defined by
\[\|u\|_{\hat{J}^k}=\|u\|_{L^6}+\sum_{\nu=0}^1\left\|(1+|x|)^{\nu+1}\nabla^\nu u\right\|_{L^\infty}+\displaystyle\left\|(1+|x|)^{2}\nabla^2u\right\|_{L^\infty}.\]

\end{Definition}
Moreover, we put
\[\Lambda_{\epsilon}^{j, k, l}=\big\{(\sigma, v, \vartheta)|\,\sigma\in I_\epsilon^j, v\in J_\epsilon^k, \vartheta\in N_\epsilon^l,\|(\sigma,v,\vartheta)\|_{\Lambda^{j, k, l}}<\epsilon\big\},\]
\[\left\|(\sigma, v, \vartheta)\right\|_{\Lambda^{j, k, l}}=\|\sigma\|_{I^k}+\|v\|_{J^k}+\|\vartheta\|_{N^k},
\]
\[\begin{array}{ll}\dot{\Lambda}_{\epsilon}^{j, k, l}=\big\{(\sigma, v, \vartheta)&\in\Lambda_{\epsilon}^{j,k,l}|\,\nabla\cdot v=\nabla\cdot V_1+V_2\, for\, some \,V_1, V_2\\[2mm]
&\left.\,such\, that \, \left\|(1+|x|)^3V_1\right\|_{L^\infty}+\left\|(1+|x|)^{-1}V_2\right\|_{L^1}\leq\epsilon\right\},\end{array}
\]
\[\mathcal{L}=\{U|U=\nabla\cdot U_1+U_2\, for\, some \, U_1, U_2\, and \,satisfies\, \|U\|_{\mathcal{L}}<\infty\},\]
where
\[\|U\|_\mathcal{L}=\sum_{\nu=1}^3\left\|(1+|x|)^{\nu+1}\nabla^\nu U\right\|+\left\|(1+|x|)^{3}(U, \nabla U)\right\|_{L^\infty}+\left\|(1+|x|)^{2}U_1\right\|_{L^\infty}+\|U_2\|_{L^1}.\]

In this paper, we consider the case where the mass source $G$, the external force $F$ and energy source $H$ are given by the following form
\[
\left(\begin{array}{ll}
        G\\[1mm]
        F\\[1mm]
     H\\[1mm]
 \end{array}\right)=\nabla\cdot\left(\begin{array}{ll}
        G_1\\[1mm]
        F_1\\[1mm]
        H_1\\[1mm]
 \end{array}\right)
 +\left(\begin{array}{ll}
        G_2\\[1mm]
        F_2\\[1mm]
        H_2\\[1mm]
 \end{array}\right)
\]
where $F_1=(F_{1,\, ij}(x))_{1\leq i,\,j\leq3}$, $F_2=(F_{2,\, i}(x))_{1\leq i\leq3}$; $G_1=(G_{1,\,i}(x))_{1\leq i\leq3}$, $G_2=G_2(x)$; $H_1=(H_{1,\,i}(x))_{1\leq i\leq3}$, $H_2=H_2(x)$.

 Now we begin to state our main results. As \cite{J. Z. Qian-H. Yin-2007}, regarding $\rho$ as a smooth function $(P, \theta)$, Our first Theorem is concerning the existence of stationary solution to (\ref{1.6}), and its weighted-$L^2$ and $L^\infty$ estimates.
\begin{Theorem}
Let $\bar{\rho}$, $\bar{\theta}$ be any positive constants, and set $\bar{P}=P(\bar{\rho},\bar{\theta})$. There exists small constants $c_0>0$ and $\epsilon_0>0$ depending on $\bar{\rho}$ and $\bar{\theta}$, such that if $(G, F, H)\in\mathcal{H}^{4,3,4}$ and satisfies the estimate:
 \[\|(G, F, H)\|_{\mathcal{L}}+\left\|(1+|x|)^4\nabla^4(G, H)\right\|+\left\|(1+|x|)^{-1}G\right\|_{L^1}\leq c_0\epsilon\]
 for some positive constant $\epsilon\leq\epsilon_0$, then (\ref{1.6}) admits a solution of the form: $(P, v, \theta)=(\bar{P}+\sigma, v , \bar{\theta}+\vartheta)$ where $(\sigma, v , \theta)\in \dot{\Lambda}_\epsilon^{4, 5, 5}$. Furthermore the solution is unique in the following sense: if there is another solution $(\bar{P}+\sigma_1, v_1, \bar{\theta}+\vartheta_1)$ satisfying (\ref{1.6}) with the same $(G, F, H)$, and $\|(\sigma_1, v_1, \vartheta_1)\|_{{\Lambda}^{4, 5, 5}}\leq\epsilon$, then $(\sigma_1, v_1, \vartheta_1)=(\sigma, v, \vartheta)$.
\end{Theorem}
Next, we consider the stability of the stationary solution of (\ref{1.6}) with respect to the initial disturbance. Let $(\rho^*, v^*, \vartheta^*)$ be the stationary solution obtained in Theorem 1.1, then the stability of $(\rho^*, v^*, \vartheta^*)$ means the solvability of the non-stationary problem (\ref{1.3}), (\ref{1.4}). Let us introduce first the class of functions which solutions of (\ref{1.3}), (\ref{1.4}) belong to.

\begin{Definition}
\[\mathcal{C}(0, T; \mathcal{H}^{j,k,l})=\left\{(\sigma, w, \vartheta)(t,x)\left|\begin{array}{c} \sigma(t,x)\in C^0(0, T; H^j)\bigcap C^1(0, T; H^{j-2}),\\[2mm]
w(t,x)\in C^0(0, T; H^k)\bigcap C^1(0, T; H^{k-2}),\\[2mm]
\vartheta(t,x)\in C^0(0, T; H^l)\bigcap C^1(0, T; H^{l-2})
\end{array}\right.\right\}
\]
\end{Definition}
Then, we have the following Theorem.
\begin{Theorem}
There exist $C>0$ and $\delta>0$ such that if $\|(\rho_0-\rho^*, v_0-v^*, \theta_0-\vartheta^*)\|_{4,3,3}\leq\delta$, then the Cauchy problem (\ref{1.3}), (\ref{1.4}) admits a unique solution $(\rho, v, \theta)=(\rho^*+\sigma, v^*+w, \theta^*+\vartheta)$ globally in time, where $(\sigma, w, \vartheta)\in \mathcal{C}(0, \infty; \mathcal{H}^{4,3,3})$, $\nabla\sigma\in L^2(0, \infty; H^4)$, $\nabla w, \nabla \vartheta\in L^2(0, \infty; H^3)$. Moreover, the solution $(\sigma, w, \vartheta)$ satisfies the estimate:
\begin{equation}\label{1.7}
\|(\sigma, w, \vartheta)(t)\|^2_{4,3,3}+\int_{0}^t\left\|\nabla(\sigma, w, \vartheta)(s)\right\|^2_{4,3,3}\,ds\leq C\|(\sigma, w, \vartheta)(0)\|^2_{4,3,3}.
\end{equation}
for any $t>0$ and
\[
\|(\sigma, v, \vartheta)(t)\|_{L^\infty}\rightarrow0\quad as \,t\rightarrow\infty.
\]
\end{Theorem}

The compressible Navier-Stokes-Kortewg system has been attracted many attentions due to its applications in fluid mechanics
as well as mathematical challenge. A lot of mathematical results on such system have been obtained. More precisely, Hattori and
Li \cite{H. Hattori-D. Li-1994, H. Hattori-D. Li-1996} proved the local existence and the global existence of smooth solutions
for the compressible fluid models of Korteweg type in Sobolev space. Danchin and Desjardins \cite{R. Danchin-B. Desjardins-2001}
proved existence and uniqueness results of suitably smooth solutions for the compressible fluid models of Korteweg type in critical
Besov space. Bresch, Desjardins and Lin \cite{D. Bresch-B. Desjardins-C. K. Lin-2003} showed the global existence of weak solution
to the compressible fluid models of Korteweg type, then Haspot improved their results in \cite{B. Haspot-2011}. The local existence
of strong solutions for the compressible fluid models of Korteweg type was proved by M. Kotschote \cite{M. Kotschote-2008}.  Wang and
Tan \cite{Y. J. Wang-Z. Tan-2011} established the optimal $L^2$ decay rates of global smooth solutions for the compressible fluid models
of Korteweg type without external force. Recently, Li \cite{Y. P. Li-2011} discussed the global existence of smooth solution to the
following Cauchy problem of the isothermal compressible fluid models of Korteweg type with potential external force.
\begin{eqnarray}\label{1.8}
\left\{\begin{array}{ll}
          \rho_t+\nabla\cdot(\rho u)=0,\\[2mm]
          (\rho u)_t+\nabla\cdot(\rho u\bigotimes u)=\nabla\cdot(S+K)+\rho F(x), \\[2mm]
          (\rho,  u)|_{t=0}=(\rho_0, u_0).
 \end{array}\right.
\end{eqnarray}
Here $F(x)=-\nabla\phi$ with $\phi$ being a scalar function and $S, K$ are defined as in (\ref{1.2}). He proved that there exists a unique stationary solution $(\tilde{\rho}(x), 0)$ to problem (\ref{1.8}) if $\phi(x)$ satisfies some smallness condition in the $H^3$ norm. The nonlinear stability of the stationary solution $(\tilde{\rho}(x), 0)$  and the optimal $L^2$-decay rate of smooth solutions to (\ref{1.8}) were also proved in \cite{Y. P. Li-2011}. Motivated by the work Y. Shibata and K. Tanaka \cite{Y. Shibata-K. Tanaka-2003} for the study of compressible Navier-Stokes equations, when the external force is given by the general form $F=\nabla\cdot F_1+F_2$ and also mass source $G$ appears, it is expect that the stationary solution is nontrivial in general. On the other hand, all the above results are concerning about the isothermal compressible fluid models of Korteweg type, for the nonisothermal compressible fluid models of Korteweg type, fewer results have been obtained. To our knowledge, the only available result for the nonisothermal case is \cite{B. Haspot-2009}, where the existence and uniqueness of strong solutions was proved in critical space. Based on these observations, we consider in this paper the nonlinear stability of stationary solutions to the full compressible Navier-Stokes-Korteweg system (\ref{1.1}).

 Now we outline the main ideas used in proving our main results. The proof of Theorem 1.1 is motivated by the method developed by Y. Shibata and K. Tanaka \cite{Y. Shibata-K. Tanaka-2003}. Firstly, as mentioned before, we choose $(P, v, \theta)$ as the independent variables and regarding $\rho$ as a smooth function of $(P,\theta)$. Then in the same sprit as \cite{Y. Shibata-K. Tanaka-2003}, we need to establish the corresponding linear theory in the $L^2$-framework for (\ref{1.6}) by employing the Banach closed range theorem. Compared with the case of compressible Navier-Stokes system, the appearance of the third order terms $\nabla\Delta\sigma$ and $\nabla\Delta\vartheta$  in the velocity equation $(2.7)_2$ result in more difficulties when we estimate the $L^2$ norm of the solutions to the approximate problem. In particular, an additional term $\nabla(\nabla\cdot v)$ appears in the energy estimate. To close the $L^2$ energy type estimate, we frequently use the structures of the approximate system. Then by choosing some suitably space-weights and multipliers, the weighted-$L^2$ estimate of solutions to the linearized problem is also obtained. In order to deal with the nonlinear problem, we have to derive the weighted-$L^\infty$ estimates for solutions $(\sigma, v, \vartheta)$ to the linearized equation ({\ref{2.81}). The weighted-$L^\infty$ estimates for $v$ and $\vartheta$ can be deduced in the same way as that of compressible Navier-Stokes equations. However, for the weighted-$L^\infty$ estimates of $\sigma$, we need to perform some delicate estimates related to the Bessel potential(see (\ref{2.98}) for detail). Moreover, the highly nonlinear terms $\Psi(\tilde{v})$ and $\Phi(\tilde{\rho}, \tilde{v})$ in $(\ref{2.82})$ are overcome by some delicate analysis. Having obtained the weighted-$L^2$ and weighted-$L^\infty$ estimates of solutions to the linearized problem, Theorem 1.1 follows by the contraction mapping principle. As for the nonlinear stability of the stationary solution obtained above, the key step is to deduce some certain a priori estimates for solutions to the initial value problem (\ref{3.1}),(\ref{3.2}) in the $H^3$ framework. Based on the properties  we obtained on the stationary solution and some delicate estimates, we can deduce the desired a priori estimates. It is worth to point out that, for the compressible Navier-Stokes-Korteweg system (\ref{1.2}), the appearance of the Korteweg tensor $\rho\nabla\Delta\rho$ results in more regularity for the density than the velocity and internal energy (see (\ref{1.7})). In fact, we frequently use integration by parts and the equation $(\ref{3.1})_1$ when we estimate the the terms containing $\nabla\Delta\sigma$. As a result, the Korteweg term is split into the energy and the terms small in the $L^2$ norm.

Another interesting problem is to investigate the convergence rate of the non-stationary solutions to the stationary solutions when the time goes to infinity. As mentioned before, this problem has been studied by some authors for the isothermal compressible Navier-Stokes-Korteweg system with $G=0, F=-\nabla\phi$  or without any external force(cf.\cite{Y. P. Li-2011}, \cite{Y. J. Wang-Z. Tan-2011, Z. Tan-H. Q. Wang-X. J. Kai-2012}). But to obtain the convergence rate in our case, it appear to be more delicate since the stationary solution is nontrivial generally. We will consider this problem in a forthcoming paper.

Before concluding this section, we also mention that the nonlinear stability of stationary solution for the compressible Navier-Stokes system has been studied by many authors. For the non-isentropic case, we refer to \cite{A. Matsumura-T. Nishida-1979, A. Matsumura-T. Nishida-1980} for the stability of constant state $(\bar{\rho}, 0, \bar{\theta})$ in $\mathbb{R}^3$, \cite{A. Matsumura-T. Nishida-1983} for the stability of nontrivial stationary solution $(\rho^*(x), 0, \bar{\theta})$ in an exterior domain of $\mathbb{R}^3$ and \cite{J. Z. Qian-H. Yin-2007, A. Matsumura-T. Nishida-1989} for the stability of generally nontrivial stationary solution $(\rho^*(x), v^*(x), \theta^*(x))$ in $\mathbb{R}^3$ and an exterior domain of $\mathbb{R}^3$, respectively. For the isentropic case, the interesting readers are referred to  \cite{Y. Shibata-K. Tanaka-2003, A. Novotny-M. Padula-1997,K. Tanaka-2006} for the stability of generally nontrivial stationary solution $(\rho^*(x), v^*(x))$ in $\mathbb{R}^3$ or an exterior domain of $\mathbb{R}^3$ and \cite{A. Novotny-M. Padula-1996} for the stability of nontrivial stationary solution $(\rho^*(x),0)$ in an exterior domain of $\mathbb{R}^3$.

The rest of this paper is organized as follows. In Section 2, we study the stationary problem.  The non-stationary problem will be studied in Section 3.

\section{Stationary problem}
\setcounter{equation}{0}
This section is devoted to the stationary problem ({\ref{1.6}). Take any two constants $\bar{\rho}, \bar{\theta}>0$. As mentioned in Section 1, by regarding $\rho$ as the function of $(P, \theta)$, changing the variables $(P, v, \theta)\rightarrow(\bar{P}+\sigma, v, \bar{\theta}+\vartheta)$, and rewriting the third equation by using the first one, (\ref{1.6}) can be then reformulated as
\begin{eqnarray}\label{2.1}
\left\{\begin{array}{ll}
\nabla\cdot v\displaystyle+\frac{\rho_P}{\rho}(v\cdot\nabla)\sigma=\displaystyle-\frac{\rho_\theta}{\rho}(v\cdot\nabla) \vartheta+\frac{G(x)}{\rho},\\[2mm]
\displaystyle-\mu\Delta v-(\mu+\mu^\prime)\nabla(\nabla\cdot v)+\nabla \sigma-\kappa\gamma_1\nabla\Delta\sigma-\kappa\gamma_2\nabla\Delta\vartheta=-\rho(v\cdot\nabla)v+\hat{f},\\[2mm]
\displaystyle-\tilde{\alpha}\Delta\vartheta=-\eta_1(v\cdot\nabla)\vartheta+\eta_2(v\cdot\nabla)\sigma+\Psi(v)+\hat{\Phi}-\eta_3G+H+\frac{v^2}{2}G(x)-C_\triangledown (\vartheta+\bar{\theta})G,
 \end{array}\right.
\end{eqnarray}
where
\begin{eqnarray}\label{2.2}
\left\{\begin{array}{ll}
\gamma_1=\bar{\rho}\bar{\rho}_P,\,\,\gamma_2=\bar{\rho}\bar{\rho}_\theta,\,\, \bar{\rho}_P=\rho_P(\bar{P},\bar{\theta}), \,\,\bar{\rho}_\theta=\rho_\theta(\bar{P},\bar{\theta}),\\[2mm]
\displaystyle\eta_1=\displaystyle\rho C_\triangledown-\frac{\theta\rho_\theta^2}{\rho\rho_P},\quad\eta_2=\displaystyle\frac{\theta\rho_\theta}{\rho},\quad\eta_3=\displaystyle\frac{\theta\rho_\theta}{\rho\rho_P},\\[2mm]
\hat{f}=\kappa\rho\left(\nabla\sigma\cdot\nabla^2\rho_P+\nabla\rho_P\cdot\nabla^2\sigma+\nabla\rho_P\Delta\sigma\right)+\kappa\left(\rho\rho_P-\bar{\rho}\bar{\rho}_P
\right)\nabla\Delta\sigma\\[2mm]
\qquad+\kappa\rho\left(\nabla\vartheta\cdot\nabla^2\rho_\theta+\nabla\rho_\theta\cdot\nabla^2\vartheta+\nabla\rho_\theta\Delta\vartheta\right)+\kappa\left(\rho\rho_\theta-\bar{\rho}\bar{\rho}_\theta
\right)\nabla\Delta\vartheta,\\[2mm]
\hat{\Phi}=\kappa\left[\frac{1}{2}|\rho_P\nabla\sigma+\rho_\theta\nabla\vartheta|^2+\rho(\nabla\rho_P\cdot\nabla\sigma
+\rho_P\Delta\sigma+\nabla\rho_\theta\cdot\nabla\vartheta+\rho_\theta\Delta\vartheta)\right]\nabla\cdot v\\[2mm]
\qquad-\kappa\left[(\rho_P\nabla\sigma+\rho_\theta\nabla\vartheta)\bigotimes(\rho_P\nabla\sigma+\rho_\theta\nabla\vartheta)\right]:\nabla v.
\end{array}\right.
\end{eqnarray}

Our goal of this section is to prove Theorem 1.1 by application of weighted $L^2$-method to the linearized problem for (\ref{2.1}).

\subsection{Weighted $L^2$ theory for linearized problem}
We shall consider the linearized equation of (\ref{2.1}):
\begin{eqnarray}\label{2.3}
\left\{\begin{array}{ll}
\nabla\cdot v\displaystyle+(a\cdot\nabla)\sigma=g,\\[2mm]
\displaystyle-\mu\Delta v-(\mu+\mu^\prime)\nabla(\nabla\cdot v)+\nabla \sigma-\kappa\gamma_1\nabla\Delta\sigma-\kappa\gamma_2\nabla\Delta\vartheta=f,\\[2mm]
\displaystyle-\tilde{\alpha}\Delta\vartheta=h,
 \end{array}\right.
\end{eqnarray}
where $a=(a_1(x), a_2(x), a_3(x))$, $(g, f, h)\in\mathcal{H}^{4, 3, 3}$ are given. Throughout this subsection, we put
\[f=-(b_1\cdot\nabla) c_1+\tilde{f}, \qquad h=-(b_2\cdot\nabla)c_2+\tilde{h}.\]
and assume that
\begin{equation}\label{2.4}
a\in\hat{H}^4, \quad \left\|(1+|x|)a\right\|_{L^\infty}+\sum_{\nu=1}^4\left\|(1+|x|)^{\nu-1}\nabla^\nu a\right\|\leq\delta
\end{equation}
\begin{equation}\label{2.5}\left(g,\tilde{f},\tilde{h}\right)\in\mathcal{H}^{4,3,4},  \quad b_1, b_2, c_1\in J^{5}, c_2\in N^{5},
\end{equation}
\begin{equation}\label{2.6}
\|(1+|x|)(g,\tilde{h})\|+\displaystyle\sum_{\nu=1}^4(1+|x|)^\nu\nabla^\nu(g,\tilde{h})\|+\displaystyle\sum_{\nu=0}^{3}(1+|x|)^{\nu+1}\nabla^\nu(\tilde{f}, \tilde{h})\|\leq\infty
\end{equation}

\subsubsection{Solution to approximate problem}
First, we solve the approximate problem:
\begin{eqnarray}\label{2.7}
\left\{\begin{array}{ll}
\nabla\cdot v\displaystyle+(a\cdot\nabla)\sigma-\epsilon\Delta\sigma+\epsilon\sigma=g,\\[2mm]
\displaystyle-\mu\Delta v-(\mu+\mu^\prime)\nabla(\nabla\cdot v)+\nabla \sigma-\kappa\gamma_1\nabla\Delta\sigma-\kappa\gamma_2\nabla\Delta\vartheta+\epsilon v=f,\\[2mm]
\displaystyle-\tilde{\alpha}\Delta\vartheta+\epsilon\vartheta=h,
 \end{array}\right.
\end{eqnarray}
in $\mathcal{H}^{3, 2, 3}$. In the following lemma, we prove some fundamental a priori estimate needed later.
\begin{Lemma}
Suppose that $(\sigma, v, \vartheta)\in\mathcal{H}^{3, 2, 3}$ is a solution to (\ref{2.7}). Then there exists two positive constants $\delta_0=\delta_0(\gamma_1, \gamma_2, \kappa, \mu, \mu^\prime, \tilde{\alpha})$ and $\epsilon_0=\epsilon_0(\gamma_1, \gamma_2, \kappa, \mu, \mu^\prime, \tilde{\alpha})<1$ such that if $\delta$ in (\ref{2.4}) satisfies $\delta\leq\delta_0$ and $0<\epsilon<\epsilon_0$,  we have the following estimate:
\begin{equation}\label{2.8}
\left\|\nabla(\sigma, v, \vartheta)\right\|^2_{2,1,2}+\epsilon\left\|(\sigma, v, \vartheta)\right\|^2_{2,1,1}\leq C\epsilon^{-1}\|(g,f,h)\|^2+C\|\nabla(g, h)\|^2.
\end{equation}
Here, $C>0$ is a constant depending only on $\gamma_1, \gamma_2, \kappa, \mu, \mu^\prime$ and $\tilde{\alpha}$.
\end{Lemma}
\noindent{\bf Proof.} The proof consists of four steps.

{\bf Step 1}. Taking the $L^2$ inner product with $\sigma$ and $v$ on $(\ref{2.7})_1$, $(\ref{2.7})_2$ , respectively, using integration by parts and canceling the term $\langle \nabla\sigma, v\rangle$ by adding the two resultant equations together, we have
\begin{equation}\label{2.9}
\begin{array}{rl}
&\mu\|\nabla v\|^2+(\mu+\mu^\prime)\|\nabla\cdot v\|^2+\epsilon\|(\sigma, v)\|^2_{1,0}\\[2mm]
&=\langle g, \sigma\rangle+\langle f, v\rangle+\kappa\gamma_1\langle \nabla\Delta\sigma, v\rangle+\kappa\gamma_2\langle \nabla\Delta\vartheta, v\rangle-\langle(a\cdot\nabla)\sigma, \sigma\rangle.
\end{array}
\end{equation}
Differentiating $(\ref{2.7})_1$ and $(\ref{2.7})_2$, and employing the same argument, we have
\begin{equation}\label{2.10}
\begin{array}{rl}
&\mu\|\nabla^2 v\|^2+(\mu+\mu^\prime)\|\nabla(\nabla\cdot v)\|^2+\epsilon\|\nabla(\sigma, v)\|^2_{1,0}\\[2mm]
&=\langle \nabla g, \nabla\sigma\rangle+\langle\nabla f, \nabla v\rangle+\kappa\gamma_1\langle \nabla(\nabla\Delta\sigma), \nabla v\rangle+\kappa\gamma_2\langle \nabla(\nabla\Delta\vartheta), \nabla v\rangle-\langle\nabla((a\cdot\nabla)\sigma), \nabla\sigma\rangle.
\end{array}
\end{equation}
Adding (\ref{2.10}) to (\ref{2.9}) yields
\begin{equation}\label{2.11}
\begin{array}{rl}
&\mu\|\nabla v\|^2_1+(\mu+\mu^\prime)\|(\nabla\cdot v)\|^2_1+\epsilon\|(\sigma, v)\|^2_{2,1}\\[2mm]
=&\displaystyle\sum_{\nu=0}^1\left\{\langle \nabla^\nu g, \nabla^\nu\sigma\rangle+\langle\nabla^\nu f, \nabla^\nu v\rangle+\kappa\gamma_1\langle \nabla^\nu(\nabla\Delta\sigma), \nabla^\nu v\rangle\right.\\[2mm]&\left.
\quad\quad\quad+\kappa\gamma_2\langle \nabla^\nu(\nabla\Delta\vartheta), \nabla^\nu v\rangle-\langle\nabla^\nu((a\cdot\nabla)\sigma), \nabla^\nu\sigma\rangle\right\}
=I_1+I_2+I_3+I_4+I_5.
\end{array}
\end{equation}
It follows from the Cauchy inequality that
\begin{eqnarray}\label{2.12}
\begin{array}{ll}
&I_1\leq\displaystyle\frac{\epsilon}{4}\|\sigma\|^2+C_\epsilon\|g\|^2+\eta\|\Delta\sigma\|^2+C_\eta\|g\|^2,\\[3mm]
&I_2\leq\displaystyle\frac{\epsilon}{4}\|v\|^2+C_\epsilon\|f\|^2+\eta\|\nabla(\nabla\cdot v)\|^2+C_\eta\|f\|^2,\\[3mm]
&I_3\leq\displaystyle\eta(\|\nabla\sigma\|^2+\|\nabla\Delta\sigma\|^2)+C_\eta\|\nabla(\nabla\cdot v)\|^2,
\end{array}
\end{eqnarray}
and
\begin{equation}\label{2.13}
I_4\leq\displaystyle\eta(\|\nabla\vartheta\|^2+\|\nabla\Delta\vartheta\|^2)+C_\eta\|\nabla(\nabla\cdot v)\|^2.
\end{equation}
Here and hereafter, $\eta>0$ denotes a sufficiently small constant and $C_\epsilon, C_\eta$ denote some positive constants depending only on $\epsilon$ and $\eta$, respectively.  Moreover, the Cauchy-Schwartz inequality and the Hardy inequality imply that
\begin{equation}\label{2.14}
\begin{array}{ll}
I_5&\leq\displaystyle C\left(\left|\langle(a\cdot\nabla)\sigma, \sigma\rangle\right|+\left|\langle(a\cdot\nabla)\sigma, \Delta\sigma\rangle\right|\right)\\[2mm]
&\leq \displaystyle C\|(1+|x|)a\|_{L^\infty}\displaystyle\left(\|\nabla\sigma\|\left\|\frac{\sigma}{|x|}\right\|+\displaystyle\left\|\frac{\nabla\sigma}{|x|}\right\|\|\Delta\sigma\|\right)\\[2mm]
&\leq C\delta\|\nabla\sigma\|_1^2.
\end{array}
\end{equation}
Combining (\ref{2.11})-(\ref{2.14}), we obtain
\begin{equation}\label{2.15}
\begin{array}{ll}
&\mu\|\nabla v\|^2_1+(\mu+\mu^\prime)\|\nabla\cdot v\|^2_1+\epsilon\|(\sigma, v)\|^2_{2,1}\\[2mm]
&\leq C\eta\|\nabla(\sigma, \vartheta)\|^2_2+C\delta\|\nabla\sigma\|^2_1+C_\eta\|\nabla(\nabla\cdot v)\|^2+(C_\epsilon+C_\eta)\|(g, f)\|^2.
\end{array}
\end{equation}

{\bf Step 2}. Differentiating $(\ref{2.7})_1$, we get
\[\nabla(\nabla\cdot v)=-\nabla((a\cdot\nabla)\sigma)+\epsilon\nabla\Delta\sigma-\epsilon\nabla\sigma+\nabla g.\]
which together with the sobolev inequality imply that
\begin{equation}\label{2.16}
\begin{array}{ll}
\|\nabla(\nabla\cdot v)\|^2&\leq\displaystyle C\left(\|(a, \nabla a)\|_{L^\infty}\|\nabla\sigma\|_1^2+\epsilon^2\|\nabla\Delta\sigma\|^2+\epsilon^2\|\nabla\sigma\|^2+\|\nabla g\|^2\right)\\[2mm]
&\leq C\left(\delta^2\|\nabla\sigma\|_1^2+\epsilon^2\|\nabla\Delta\sigma\|^2+\epsilon^2\|\nabla\sigma\|^2+\|\nabla g\|^2\right)
\end{array}
\end{equation}

{\bf Step 3}. Taking the $L^2$ inner product with $\nabla\sigma$ on $(\ref{2.7})_2$, we have from the Cauchy inequality that
\begin{equation}\label{2.17}
\begin{array}{ll}
 &\|\nabla\sigma\|^2+\kappa\gamma_1\|\Delta\sigma\|^2\\[2mm]
 &=\mu\langle\Delta v, \nabla\sigma\rangle+(\mu+\mu^\prime)\langle\nabla(\nabla\cdot v), \nabla\sigma\rangle+\kappa\gamma_2\langle \nabla\Delta\vartheta, \nabla \sigma\rangle-\epsilon\langle v, \nabla\sigma\rangle+\langle f, \nabla\sigma\rangle\\[2mm]
 &\leq \displaystyle \frac{1}{2}\|\nabla\sigma\|^2+C\left(\|\Delta v\|^2+\|\nabla(\nabla\cdot v)\|^2+\|\nabla\Delta\vartheta\|^2+\epsilon^2\|v\|^2+\|f\|^2\right)
\end{array}
\end{equation}
Consequently,
\begin{equation}\label{2.18}
\|\nabla\sigma\|^2+\|\Delta\sigma\|^2\leq C\left(\|\Delta v\|^2+\|\nabla\Delta\vartheta\|^2+\epsilon^2\|v\|^2+\|f\|^2\right).
\end{equation}
On the other hand, it follows from $(\ref{2.7})_2$ that
\begin{equation}\label{2.19}
\|\nabla\Delta\sigma\|^2\leq C\left(\|\Delta v\|^2+\|\nabla\Delta\vartheta\|^2+\|\nabla\sigma\|^2+\epsilon^2\|v\|^2+\|f\|^2\right).
\end{equation}
Therefore, we have from a linear combination of (\ref{2.18}) and (\ref{2.19}) that
 \begin{equation}\label{2.20}
\|\nabla\sigma\|^2_2\leq C\left(\|\Delta v\|^2+\|\nabla\Delta\vartheta\|^2+\epsilon^2\|v\|^2+\|f\|^2\right).
\end{equation}

{\bf Step 4}. By using the same argument as (\ref{2.9}) and (\ref{2.10}), one can get from $(\ref{2.7})_3$ that
 \begin{equation}\label{2.21}
 \begin{array}{ll}
\tilde{\alpha}\|\nabla\vartheta\|^2_1+\epsilon\|\vartheta\|_1^2&=\langle h, \vartheta\rangle+\langle \nabla h, \nabla\vartheta\rangle\\[2mm]
 &\leq\displaystyle \frac{\epsilon}{2}\|\vartheta\|^2+\frac{\tilde{\alpha}}{2}\|\nabla^2\vartheta\|^2+(C_\epsilon+C)\|h\|^2
 \end{array}
\end{equation}
which implies
\begin{equation}\label{2.22}
\tilde{\alpha}\|\nabla\vartheta\|^2_1+\epsilon\|\vartheta\|_1^2\leq\displaystyle(C_\epsilon+C)\|h\|^2
\end{equation}
On the other hand,
 \begin{equation}\label{2.23}
\tilde{\alpha}\|\nabla\Delta\vartheta\|^2=\|-\epsilon\nabla\vartheta+\nabla h\|^2\leq\epsilon^2\|\nabla\vartheta\|^2+\|\nabla h\|^2
\end{equation}
Combining (\ref{2.22}) and (\ref{2.23}), we obtain
\begin{equation}\label{2.24}
\|\nabla\vartheta\|^2_2+\epsilon\|\vartheta\|_1^2\leq\displaystyle C_\epsilon\|h\|^2+C\|\nabla h\|^2.
\end{equation}
Thus, by some suitably linear combinations of (\ref{2.15}), (\ref{2.16}), (\ref{2.20}) and (\ref{2.24}) and using the smallness of $\epsilon, \eta $ and $\delta$, we can get (\ref{2.8}). This completes the proof of Lemma 2.1.

Now, we employ the closed range theorem to prove the existence of solution to (\ref{2.7}). We introduce the operator $A$ defined on $D(A)\subset L^2$ into $H^1\times L^2\times H^1$ by
\[A(\sigma, v, \vartheta)=\left(A_1(\sigma, v, \vartheta), A_2(\sigma, v, \vartheta), A_3(\sigma, v, \vartheta)\right)\]
where $D(A)=\mathcal{H}^{3,2,3}$ and
\[\left\{\begin{array}{ll}
 &A_1(\sigma, v, \vartheta)=\nabla\cdot v\displaystyle+(a\cdot\nabla)\sigma-\epsilon\Delta\sigma+\epsilon\sigma,\\[2mm]
 &A_2(\sigma, v, \vartheta)=-\mu\Delta v-(\mu+\mu^\prime)\nabla(\nabla\cdot v)+\nabla \sigma-\kappa\gamma_1\nabla\Delta\sigma-\kappa\gamma_2\nabla\Delta\vartheta+\epsilon v,\\[2mm]
 &A_3(\sigma, v, \vartheta)= -\tilde{\alpha}\Delta\vartheta+\epsilon\vartheta.
\end{array}\right.
\]
Clearly, $A$ is closed operator. Furthermore, Lemma 2.1 implies that for each $0<\epsilon<\epsilon_0$, the range of $A$ is closed.
\begin{Proposition}
 There exists two positive constants $\delta_0=\delta_0(\gamma_1, \gamma_2, \kappa, \mu, \mu^\prime, \tilde{\alpha})$ and $\epsilon_0=\epsilon_0(\gamma_1, \gamma_2, \kappa, \\ \mu, \mu^\prime, \tilde{\alpha})<1$ such that if $\delta$ in (\ref{2.4}) satisfies $\delta\leq\delta_0$ and $0<\epsilon<\epsilon_0$,  then (\ref{2.7}) has a solution $(\sigma, v, \vartheta)\in \mathcal{H}^{3,2,3}$, which satisfies
\begin{equation}\label{2.25}
\left\|(\sigma, v, \vartheta)\right\|_{3,2,3}\leq C(\epsilon)(\|(g,f,h)\|+\|\nabla(g, h)\|).
\end{equation}
where $C(\epsilon)>0$ is a constant depending only on $\gamma_1, \gamma_2, \kappa, \mu, \mu^\prime, \tilde{\alpha}$ and $\epsilon$, and $C(\epsilon)\rightarrow\infty$ as $\epsilon\rightarrow0$.
\end{Proposition}
\noindent{\bf Proof.} Firstly, for any $(\sigma, v, \vartheta)\in \mathcal{H}^{3,2,3}$ and $(\sigma_*, v_*, \vartheta_*)\in \mathcal{H}^{\infty,\infty,\infty}$, we have from integration by parts that
\begin{eqnarray}\label{2.26}
\begin{array}{ll}
 \langle A(\sigma, v, \vartheta), (\sigma_*, v_*, \vartheta_*)\rangle&=\langle \nabla\cdot v\displaystyle+(a\cdot\nabla)\sigma-\epsilon\Delta\sigma+\epsilon\sigma, \sigma_*\rangle+\langle -\tilde{\alpha}\Delta\vartheta+\epsilon\vartheta, \vartheta_*\rangle\\[2mm]
 &\quad+\langle -\mu\Delta v-(\mu+\mu^\prime)\nabla(\nabla\cdot v)+\nabla \sigma-\kappa\gamma_1\nabla\Delta\sigma-\kappa\gamma_2\nabla\Delta\vartheta+\epsilon v, v_*\rangle\\[2mm]
 &=\langle \sigma, -\nabla\cdot v_*-\nabla\cdot(a\sigma_*)-\epsilon\Delta\sigma_*+\epsilon\sigma_*+\kappa\gamma_1\Delta(\nabla\cdot v_*)\rangle\\[2mm]
 &\quad+\langle v, -\mu\Delta v_*-(\mu+\mu^\prime)\nabla(\nabla\cdot v_*)-\nabla\sigma_*+\epsilon v_*\rangle\\[2mm]
 &\quad+\langle \vartheta, -\tilde{\alpha}\Delta \vartheta_*+\epsilon\vartheta_*+\kappa\gamma_2\Delta(\nabla\cdot v_*)\rangle
 \end{array}
\end{eqnarray}
Therefore, $D(A^*)=\mathcal{H}^{2,3,2}$ and for any $(\sigma_*, v_*, \vartheta_*)\in\mathcal{H}^{2,3,2}$,
\[A^*(\sigma_*, v_*, \vartheta_*)=\left(A^*_1(\sigma_*, v_*, \vartheta_*), A^*_2(\sigma_*, v_*, \vartheta_*), A^*_3(\sigma_*, v_*, \vartheta_*)\right),\]
where
\begin{equation}\label{2.27}
\left\{\begin{array}{ll}
 &A_1^*(\sigma_*, v_*, \vartheta_*)=-\nabla\cdot v_*-\nabla\cdot(a\sigma_*)-\epsilon\Delta\sigma_*+\epsilon\sigma_*+\kappa\gamma_1\Delta(\nabla\cdot v_*),\\[2mm]
 &A_2^*(\sigma_*, v_*, \vartheta_*)=-\mu\Delta v_*-(\mu+\mu^\prime)\nabla(\nabla\cdot v_*)-\nabla\sigma_*+\epsilon v_*,\\[2mm]
 &A_3^*(\sigma_*, v_*, \vartheta_*)=-\tilde{\alpha}\Delta \vartheta_*+\epsilon\vartheta_*+\kappa\gamma_2\Delta(\nabla\cdot v_*).
\end{array}\right.
\end{equation}
Employing the same argument as in the proof of Lemma 2.1 and using the equation:
\[\Delta(\nabla\cdot v_*)=-\frac{1}{2\mu+\mu^\prime}\left(\Delta\sigma_*-\epsilon\nabla\cdot v_*+\nabla\cdot A_2^*\right)\]
which follows by taking the divergence "$\nabla\cdot$ " on both side of $(\ref{2.27})_2$, one can get
\begin{equation}\label{2.28}
\|\Delta(\nabla\cdot v_*)\|+\|\nabla(\sigma_*, v_*, \vartheta_*)\|_1+\epsilon\|(\sigma_*, v_*, \vartheta_*)\|_{2,1,1}\leq C_\epsilon\|(A_1^*, A_2^*, A_3^*, \nabla\cdot A_2^*)\|
\end{equation}
Hence the closed range theorem implies the existence of solution to (\ref{2.7}). (\ref{2.25}) follows directly from (\ref{2.8}). This completes the proof of Proposition 2.1.

\subsubsection{Solution to linearized problem (\ref{2.3}) and its $L^2$ estimate}
In the following Lemma, we discuss the estimate for solution to (\ref{2.7}) independent of $\epsilon$.
\begin{Lemma}
Assume that $(\sigma, v, \vartheta)\in\mathcal{H}^{3,2,3}$ is a solution of the approximate problem (\ref{2.7}). Then there exists a constant $\delta_0=\delta_0(\gamma_1, \gamma_2, \kappa, \mu, \mu^\prime, \tilde{\alpha})>0$ such that such that if $\delta$ in (\ref{2.4}) satisfies $\delta\leq\delta_0$,  we have the estimate
\begin{equation}\label{2.29}
\|\nabla(\sigma, v, \vartheta)\|_{5, 4, 5}\leq C\left\{\|(1+|x|)(g,f,h)\|+\|\nabla(g,f,h)\|_{3, 2, 3}\right\},
\end{equation}
where the constant $C>0$ depends only on $\gamma_1, \gamma_2, \kappa, \mu, \mu^\prime$ and $\tilde{\alpha}$.
\end{Lemma}
\noindent{\bf Proof.~~}Using the Friedrichs mollifier, we may assume that $(\sigma, v, \vartheta)\in\mathcal{H}^{\infty, \infty, \infty}$. By the same argument as in the proof of Lemma 2.1, we have
\begin{equation}\label{2.30}
\|\nabla(\sigma, v, \vartheta)\|^2_{2, 1,2} \leq C_1\Big\{\|f\|^2+\|\nabla (g,f)\|^2+\sum_{\nu=0}^1\left(\langle \nabla^\nu g, \nabla^\nu\sigma\rangle+\langle\nabla^\nu f, \nabla^\nu v\rangle+\langle\nabla^\nu h, \nabla^\nu\vartheta\rangle\right)\Big\}
\end{equation}
For the third term on the right hand of (\ref{2.30}), the Cauchy inequality and the Hardy inequality imply that
\begin{equation}\label{2.31}
\begin{array}{ll}
&\displaystyle\sum_{\nu=0}^1\left(\langle \nabla^\nu g, \nabla^\nu\sigma\rangle+\langle\nabla^\nu f, \nabla^\nu v\rangle+\langle\nabla^\nu h, \nabla^\nu\vartheta\rangle\right)\\[2mm]
&\quad\quad\leq\displaystyle\frac{1}{2C_1}\|\nabla(\sigma, v, \vartheta)\|_1^2+C\|(1+|x|)(g,f,h)\|^2
\end{array}
\end{equation}
Consequently,
\begin{equation}\label{2.32}
\|\nabla(\sigma, v, \vartheta)\|^2_{2, 1,2} \leq C\left\{\|(1+|x|)(g,f,h)\|^2+\|\nabla(g,h)\|^2\right\}
\end{equation}
where the constant $C$ depends only on $\gamma_1, \gamma_2, \kappa, \mu, \mu^\prime$ and $\tilde{\alpha}$.

Moreover, for any multi-index $\alpha$ with $1\leq|\alpha|\leq k-1$, applying $\partial_x^\alpha$ to (\ref{2.7}), we have
\begin{eqnarray}\label{2.33}
\left\{\begin{array}{ll}
\nabla\cdot \partial_x^\alpha v\displaystyle+(a\cdot\nabla)\partial_x^\alpha\sigma-\epsilon\Delta \partial_x^\alpha\sigma+\epsilon \partial_x^\alpha\sigma= \partial_x^\alpha g-I_\alpha,\\[2mm]
\displaystyle-\mu\Delta\partial_x^\alpha v-(\mu+\mu^\prime)\nabla(\nabla\cdot\partial_x^\alpha v)+\nabla \partial_x^\alpha\sigma-\kappa\gamma_1\nabla\Delta \partial_x^\alpha\sigma-\kappa\gamma_2\nabla\Delta \partial_x^\alpha\vartheta+\epsilon \partial_x^\alpha v= \partial_x^\alpha f,\\[2mm]
\displaystyle-\tilde{\alpha}\Delta \partial_x^\alpha\vartheta+\epsilon \partial_x^\alpha\vartheta= \partial_x^\alpha h,
 \end{array}\right.
\end{eqnarray}
where \[I_\alpha=\displaystyle\sum_{\beta<\alpha}C_\alpha^\beta\left(\partial_x^{\alpha-\beta}a\cdot\nabla\right)\partial^\beta_x\sigma\]
with $C^\alpha_\beta$ being the binomial coefficients
corresponding to multi-indices. Notice that the third term on the right hand of (\ref{2.30}) can also be estimated as follows:
\begin{equation}\label{2.34}
\begin{array}{ll}
&\displaystyle\sum_{\nu=0}^1\left(\langle \nabla^\nu g, \nabla^\nu\sigma\rangle+\langle\nabla^\nu f, \nabla^\nu v\rangle+\langle\nabla^\nu h, \nabla^\nu\vartheta\rangle\right)\\[2mm]
&\quad\leq\displaystyle\frac{1}{2C_1}\|\nabla^2(\sigma, v, \vartheta)\|^2+C\left\{\|(g,f,h)\|^2+\|(\sigma, v, \vartheta)\|^2\right\}
\end{array}
\end{equation}
Thus, it follows from (\ref{2.30}) and (\ref{2.34}) that
\begin{equation}\label{2.35}
\|\nabla(\sigma, v, \vartheta)\|^2_{2, 1,2} \leq C\left\{\|(g,f,h)\|^2+\|\nabla(g,h)\|^2+\|(\sigma, v, \vartheta)\|^2\right\}
\end{equation}
Applying (\ref{2.35}) to (\ref{2.33}), we have
\begin{equation}\label{2.36}
\|\nabla\partial_x^\alpha(\sigma, v, \vartheta)\|^2_{2, 1,2} \leq C\left\{\|\partial_x^\alpha(g,f,h)\|^2+\|\nabla\partial_x^\alpha(g,h)\|^2+\|\partial_x^\alpha(\sigma, v, \vartheta)\|^2+\|I_\alpha\|_1^2\right\}
\end{equation}
Since
\begin{equation}\label{2.37}
\|I_\alpha\|^2_{1}\leq C\delta^2\|\nabla\sigma\|^2_{|\alpha|}
\end{equation}
as follows from the Sobolev inequality and the assumption (\ref{2.4}). We get from (\ref{2.36}) and (\ref{2.37}) that
\begin{equation}\label{2.38}
\begin{array}{ll}
&\|\nabla^{|\alpha|+3}\sigma, \nabla^{|\alpha|+2}v, \nabla^{|\alpha|+3}\vartheta\|^2\\[2mm]
&\leq C\left\{\|\partial_x^\alpha(g,f,h)\|^2+\|\nabla\partial_x^\alpha(g,h)\|^2+\|\partial_x^\alpha(\sigma, v, \vartheta)\|^2+\delta^2\|\nabla\sigma\|^2_{|\alpha|}\right\}
\end{array}
\end{equation}
Combining (\ref{2.32}) and (\ref{2.38}), we obtain (\ref{2.29}) if $\delta>0$ is small enough. This completes the proof of Lemma 2.2.

Now, we are ready to show the existence of solution to the linearized problem (\ref{2.3}) by using (\ref{2.29}).
\begin{Proposition}
There exists $\delta_0=\delta_0(\gamma_1, \gamma_2, \kappa, \mu, \mu^\prime, \tilde{\alpha})>0$ such that such that if $\delta$ in (\ref{2.4}) satisfies $\delta\leq\delta_0$, then the linearized problem (\ref{2.3}) admits a solution $(\sigma, v, \vartheta)\in\hat{\mathcal{H}}^{6, 5, 6}$ which satisfies the estimate:
\begin{equation}\label{2.39}
\|(\sigma, v, \vartheta)\|_{L^6}+\|\nabla(\sigma, v, \vartheta)\|_{5, 4, 5}\leq C\left\{\|(1+|x|)(g,f,h)\|+\|\nabla(g, f, h)\|_{3, 2, 3}\right\}
\end{equation}
where the constant $C>0$ depends only on $\gamma_1, \gamma_2, \kappa, \mu, \mu^\prime$ and $\tilde{\alpha}$.
\end{Proposition}
\noindent{\bf Proof.~~} Set \[K=\|(1+|x|)(g,f,h)\|+\|\nabla(g, f, h)\|_{3, 2, 3}\]
From Proposition 2.1 and Lemma 2.2, it follows that for each $0<\epsilon<\epsilon_0$($\epsilon_0$ is given in Proposition 2.1), (\ref{2.7}) admits a solution $(\sigma^\epsilon, v^\epsilon, \vartheta^\epsilon)\in\mathcal{H}^{6, 5, 6}$ which satisfies
\[\|\nabla(\sigma^\epsilon, v^\epsilon, \vartheta^\epsilon)\|_{5, 4, 5}\leq CK.\]
The Sobolev inequality imply that
\[\|(\sigma^\epsilon, v^\epsilon, \vartheta^\epsilon)\|_{L^6}\leq C\|\nabla(\sigma^\epsilon, v^\epsilon, \vartheta^\epsilon)\|\leq CK.\]
Choosing an appropriate subsequence, there exist $(\sigma, v, \vartheta)\in L^6$, $(\hat{\sigma}_i, \hat{v}_i, \hat{\vartheta}_i)\in\mathcal{H}^{5, 4, 5}$ such that
\[
(\sigma^\epsilon, v^\epsilon, \vartheta^\epsilon)\rightharpoonup(\sigma, v, \vartheta) \quad weakly \,\,in \,\,L^6
\]
\[
\left(\frac{\partial\sigma^\epsilon}{\partial x_i}, \frac{\partial v^\epsilon}{\partial x_i}, \frac{\partial\vartheta^\epsilon}{\partial x_i}\right)\rightharpoonup\left(\hat{\sigma}_i, \hat{v}_i, \hat{\vartheta}_i\right)\quad weakly\,\, in\,\, \mathcal{H}^{5, 4, 5}.
\]
as $\epsilon\rightarrow0$, then one can check easily that
\[
\left(\frac{\partial\sigma}{\partial x_i}, \frac{\partial v}{\partial x_i}, \frac{\partial\vartheta}{\partial x_i}\right)=\left(\hat{\sigma}_i, \hat{v}_i, \hat{\vartheta}_i\right),\]
and
\[\|(\sigma, v, \vartheta)\|+\|\nabla(\sigma, v, \vartheta)\|_{5, 4, 5}\leq CK.
\]
On the other hand, we have
\[\begin{array}{ll}
&\nabla\cdot v^\epsilon\displaystyle+(a\cdot\nabla)\sigma^\epsilon-\epsilon\Delta\sigma^\epsilon+\epsilon\sigma^\epsilon\longrightarrow\nabla\cdot v\displaystyle+(a\cdot\nabla)\sigma,\\[2mm]
&-\mu\Delta v^\epsilon-(\mu+\mu^\prime)\nabla(\nabla\cdot v^\epsilon)+\nabla \sigma^\epsilon-\kappa\gamma_1\nabla\Delta\sigma^\epsilon-\kappa\gamma_2\nabla\Delta\vartheta^\epsilon+\epsilon v^\epsilon\\[2mm]&\longrightarrow
-\mu\Delta v-(\mu+\mu^\prime)\nabla(\nabla\cdot v)+\nabla \sigma-\kappa\gamma_1\nabla\Delta\sigma-\kappa\gamma_2\nabla\Delta\vartheta,\\[2mm]
&-\tilde{\alpha}\Delta\vartheta^\epsilon+\epsilon\vartheta^\epsilon\longrightarrow-\tilde{\alpha}\Delta\vartheta.
\end{array}\]
in distribution sense. This completes the proof of Proposition 2.2.

\subsubsection{Weighted $L^2$ estimate for solution to the linearized equation (\ref{2.3})}
In this subsection, we give the weighted $L^2$ estimate for the solution to (\ref{2.3}).
\begin{Lemma}
Let $(\sigma, v, \vartheta)\in\hat{\mathcal{H}}^{6, 5, 6}$ be a solution to (\ref{2.3}) which satisfies (\ref{2.39}). Then there exists a constant $\delta_0=\delta_0(\gamma_1, \gamma_2, \kappa, \mu, \mu^\prime, \tilde{\alpha})>0$ such that such that if $\delta$ in (\ref{2.4}) satisfies $\delta\leq\delta_0$, we have for any $1\leq l\leq 4$ that
\begin{equation}\label{2.40}
\begin{array}{ll}
&\displaystyle\sum_{\nu=1}^l\Big\{\left\|(1+|x|)^\nu\left(\nabla^\nu\sigma, \nabla^{\nu+1}\sigma, \nabla^{\nu+2}\sigma\right)\right\|+\left\|(1+|x|)^\nu\left(\nabla^{\nu+1}v, \nabla^{\nu+1}\vartheta, \nabla^{\nu+2}\vartheta\right)\right\|\Big\}\\[2mm]
&\leq C\left\{\|\nabla(\sigma, v, \vartheta)\|+\|b_1\|_{J^{5}}\|c_1\|_{J^{5}}+\|b_2\|_{J^{5}}\|c_2\|_{N^{5}}+\displaystyle\sum_{\nu=1}^l\left\|(1+|x|)^\nu\nabla^{\nu-1}(\tilde{f}, \tilde{h})\right\|\right.\\[2mm]
&\left.\qquad+\displaystyle\sum_{\nu=1}^l\left\|(1+|x|)^\nu\nabla^{\nu}(g, \tilde{h})\right\|\right\}
\end{array}
\end{equation}
where $C>0$ is a constant depending only on $\gamma_1, \gamma_2, \kappa, \mu, \mu^\prime$ and $\tilde{\alpha}$.
\end{Lemma}
\noindent{\bf Proof.~~}The proof is divided into four steps.

{\bf Step 1.} Using the Friedrichs mollifier, we may assume that $(\sigma, v, \vartheta)\in\mathcal{H}^{\infty, \infty, \infty}$. For any multi-index $\alpha$ with $|\alpha|=l-1$, applying $\partial_x^\alpha$ to $(\ref{2.3})_2$, then taking the $L^2$ inner product with $(1+|x|)^{2l}\nabla\partial_x^\alpha\sigma$ on the resultant equation and summing up $\alpha$, we have
\begin{equation}\label{2.41}
\begin{array}{ll}
&\displaystyle\left\|(1+|x|)^l(\nabla^l\sigma, \nabla^{l+1}\sigma)\right\|^2\\[2mm]
&\leq C\left\{\left|\langle|\nabla^{l+1}v|, (1+|x|)^{2l}\nabla^l\sigma\rangle\right|+\left|\langle\nabla^{l}\sigma, (1+|x|)^{2l-1}\nabla^{l+1}\sigma\rangle\right|\right.\\[2mm]
&\left.\,\,\,\,\,+\left|\langle \nabla^l\Delta\vartheta,  (1+|x|)^{2l}\nabla^{l}\sigma\rangle\right|+\left|\langle \nabla^{l-1}\{-(b_1\cdot\nabla)c_1+\tilde{f}\},  (1+|x|)^{2l}\nabla^{l}\sigma\rangle\right|\right\}
\end{array}
\end{equation}
From $(\ref{2.3})_2$, we also obtain
\begin{equation}\label{2.42}
\begin{array}{ll}
\displaystyle\left\|(1+|x|)^l \nabla^{l+2}\sigma\right\|^2
&\leq C\left\{\left\|(1+|x|)^l\nabla^{l}\sigma\right\|^2+\left\|(1+|x|)^l\nabla^{l+1}v\right\|^2+\left\|(1+|x|)^l\nabla^{l-1}\tilde{f}\right\|^2\right.\\[2mm]
&\quad+\left|\langle \nabla^l\Delta\vartheta+\nabla^{l-1}\{(b_1\cdot\nabla)c_1\}, (1+|x|)^{2l}\nabla^{l+2}\sigma\rangle\right|\bigg\}
\end{array}
\end{equation}
Thus, it follows from a linear combination of (\ref{2.41}) and (\ref{2.42}) that
\begin{equation}\label{2.43}
\begin{array}{ll}
&\displaystyle\left\|(1+|x|)^l(\nabla^l\sigma, \nabla^{l+1}\sigma, \nabla^{l+2}\sigma)\right\|^2\\[2mm]
&\leq C\left\{\left\|(1+|x|)^l\nabla^{l+1}v\right\|^2+\left\|(1+|x|)^l\nabla^{l-1}\tilde{f}\right\|^2+\Big(\left|\langle|\nabla^{l+1}v|, (1+|x|)^{2l}\nabla^l\sigma\rangle\right|
\right.\\[3mm]
&\quad+\left|\langle\nabla^{l}\sigma, (1+|x|)^{2l-1}\nabla^{l+1}\sigma\rangle\right|+\left|\langle \nabla^l\Delta\vartheta,  (1+|x|)^{2l}\nabla^{l}\sigma\rangle\right|\\[3mm]
&\quad+\left|\langle \nabla^{l-1}\tilde{f},  (1+|x|)^{2l}\nabla^{l}\sigma\rangle\right|+\left|\langle \nabla^l\Delta\vartheta, (1+|x|)^{2l}\nabla^{l+2}\sigma\rangle\right|\Big)\\[3mm]
&\quad+\left|\langle \nabla^{l-1}\{(b_1\cdot\nabla)c_1\}, (1+|x|)^{2l}\nabla^{l+2}\sigma\rangle\right|+\left|\langle \nabla^{l-1}\{(b_1\cdot\nabla)c_1\}, (1+|x|)^{2l}\nabla^{l}\sigma\rangle\right|\bigg\}\\[3mm]&
=C\left\{\left\|(1+|x|)^l\nabla^{l+1}v\right\|^2+\left\|(1+|x|)^l\nabla^{l-1}\tilde{f}\right\|^2+I_1+I_2+I_3\right\}.
\end{array}
\end{equation}
The cauchy inequality implies that
\begin{equation}\label{2.44}
\begin{array}{ll}
I_1&\leq\eta\left\|(1+|x|)^l(\nabla^l\sigma, \nabla^{l+1}\sigma, \nabla^{l+2}\sigma)\right\|^2+C_\eta\Big\{\|(1+|x|)^l\nabla^{l+1}v\|^2\\[2mm]
&\quad+\|(1+|x|)^{l-1}\nabla^{l}\sigma\|^2+\|(1+|x|)^l\nabla^{l+2}\vartheta\|^2+\|(1+|x|)^l\nabla^{l-1}\tilde{f}\|^2\Big\}
\end{array}
\end{equation}
For $I_2$, notice that
\begin{equation}\label{2.45}
\left\|(1+|x|)^{|\alpha|+|\beta|}|\partial_x^\alpha b_1||\partial_x^\beta c_1|\right\|\leq C\|b\|_{J^{5}}\|c\|_{J^{5}}
\end{equation}
for any multi-index $\alpha, \beta$ with $|\alpha|\leq1$ or $|\beta|\leq1$ and $|\alpha|, |\beta|\leq 5$. If $1\leq l\leq3$, since
\[(1+|x|)^l\nabla^{l-1}\{(b_1\cdot\nabla)c_1\}\in L^2\]
as follows from (\ref{2.45}), we have
\begin{equation}\label{2.46}
I_{2}\leq\eta\left\|(1+|x|)^l\nabla^{l+2}\sigma\right\|^2+C_\eta\|b_1\|^2_{J^{5}}\|c_1\|^2_{J^{5}}.
\end{equation}
If $l=4$, we get from integration by parts that
\begin{equation}\label{2.47}
\begin{array}{ll}
I_2&\leq\left|\langle \nabla^{l}\{(b_1\cdot\nabla)c_1\}, (1+|x|)^{2l}\nabla^{l+1}\sigma\rangle\right|+\left|\langle \nabla^{l-1}\{(b_1\cdot\nabla)c_1\}, (1+|x|)^{2l-1}\nabla^{l+1}\sigma\rangle\right|\\[2mm]
&=I_{2,1}+I_{2,2}
\end{array}
\end{equation}
Using the Leibniz formula, we have
\begin{equation}\label{2.48}
\begin{array}{ll}
I_{2,1}&\leq\displaystyle\sum_{|\alpha|=4}\displaystyle\sum_{\beta\leq\alpha,\,\, |\beta|=1,2}C_\alpha^\beta\left|\langle (\partial_x^{\alpha-\beta}b_1\cdot\nabla)\partial_x^{\beta}c_1, (1+|x|)^{2l}\nabla^{l+1}\sigma\rangle\right|\\[3mm]
&\quad+\displaystyle\sum_{|\alpha|=4}\displaystyle\sum_{\beta\leq\alpha,\,\, |\beta|=0,3,4}C_\alpha^\beta\left|\langle (\partial_x^{\alpha-\beta}b_1\cdot\nabla)\partial_x^{\beta}c_1, (1+|x|)^{2l-1}\nabla^{l+1}\sigma\rangle\right|\\[3mm]
&=I^1_{2,1}+I^2_{2,1}
\end{array}
\end{equation}
By integration by parts, $I^1_{2,1}$ can be estimated as follows
\begin{equation}\label{2.49}
\begin{array}{ll}
I^1_{2,1}&\leq C\left\|(1+|x|)^2\nabla b_1\right\|_{L^\infty}\Big\{\left\|(1+|x|)^2\nabla^3 c_1\right\|\left\|(1+|x|)^4\nabla^6\sigma\right\|\\[3mm]
&\quad+\displaystyle\sum_{\nu=1}^2\left\|(1+|x|)^\nu\nabla^{\nu+2}c_1\right\|\left\|(1+|x|)^4\nabla^{5}\sigma\right\|\Big\}\\[3mm]
&\quad+(\mbox{the same term except for the exchang of}\,\, b_1 \,\,and\,\, c_1)\\[3mm]
&\leq\eta\left\|(1+|x|)^4(\nabla^{5}\sigma, \nabla^{6}\sigma)\right\|^2+C_\eta\|b_1\|^2_{J^{5}}\|c_1\|^2_{J^{5}}
\end{array}
\end{equation}
For $I^2_{2,1}$, we deduce from (\ref{2.45}) that
\begin{equation}\label{2.50}
I^2_{2,1}\leq\eta\left\|(1+|x|)^4\nabla^5\sigma\right\|^2+C_\eta\|b_1\|^2_{J^{5}}\|c_1\|^2_{J^{5}}.
\end{equation}
Combining (\ref{2.47})-(\ref{2.49}), we obtain
\begin{equation}\label{2.51}
I_{2,1}\leq\eta\left\|(1+|x|)^4(\nabla^{5}\sigma, \nabla^{6}\sigma)\right\|^2+C_\eta\|b_1\|^2_{J^{5}}\|c_1\|^2_{J^{5}}.
\end{equation}
Similarly, we can also get
\begin{equation}\label{2.52}
I_{2,2}\leq\eta\left\|(1+|x|)^4(\nabla^{5}\sigma, \nabla^{6}\sigma)\right\|^2+C_\eta\|b_1\|^2_{J^{5}}\|c_1\|^2_{J^{5}}.
\end{equation}
Thus, it follows from (\ref{2.46}), (\ref{2.47}), (\ref{2.51}) and (\ref{2.52}) that
\begin{equation}\label{2.53}
I_{2}\leq\eta\left\|(1+|x|)^l(\nabla^{l+1}\sigma, \nabla^{l+2}\sigma)\right\|^2+C_\eta\|b_1\|^2_{J^{5}}\|c_1\|^2_{J^{5}}.
\end{equation}
Similar to the estimate of $I_2$, we have if $1\leq l\leq3$,
\begin{eqnarray}\label{2.54}
I_{3}\leq\eta\left\|(1+|x|)^l\nabla^{l}\sigma\right\|^2+C_\eta\|b_1\|^2_{J^{5}}\|c_1\|^2_{J^{5}},
\end{eqnarray}
and if $l=4$,
\begin{eqnarray}\label{2.55}
I_{3}\leq\eta\left\|(1+|x|)^4\nabla^{4}\sigma\right\|^2+C_\eta\left(\|b_1\|^2_{J^{5}}\|c_1\|^2_{J^{5}}+
\left\|(1+|x|)^3\nabla^3\sigma\right\|^2\right).
\end{eqnarray}
Substituting (\ref{2.44}), (\ref{2.53})-(\ref{2.55}) into (\ref{2.43}), we arrive at
\begin{equation}\label{2.56}
\begin{array}{ll}
&\displaystyle\left\|(1+|x|)^l(\nabla^l\sigma, \nabla^{l+1}\sigma, \nabla^{l+2}\sigma)\right\|^2\\[2mm]
&\leq C\left\{\left\|(1+|x|)^l\nabla^{l+1}v\right\|^2+\left\|(1+|x|)^l\nabla^{l-1}\tilde{f}\right\|^2+\left\|(1+|x|)^l\nabla^{l+1}\vartheta\right\|^2\right.\\[2mm]
&\left.\quad+\left\|(1+|x|)^{l-1}(\nabla^{l-1}\sigma, \nabla^l\sigma)\right\|^2+\|b_1\|^2_{J^{5}}\|c_1\|^2_{J^{5}}\right\}.
\end{array}
\end{equation}

{\bf Step 2.} For any multi-index $\alpha$ with $|\alpha|=l$, applying $\partial_x^\alpha$ to $(\ref{2.3})_2$, then taking the $L^2$ inner product with $(1+|x|)^{2l}\partial_x^\alpha v$ on the resultant equation, we have from integration by parts that
\begin{equation}\label{2.57}
\begin{array}{ll}
&\mu\|(1+|x|)^l\nabla\partial_x^\alpha v\|^2+\mu\langle\nabla\partial_x^\alpha v, 2l(1+|x|)^{2l-1}\frac{x}{|x|}\partial_x^\alpha v\rangle+(\mu+\mu^\prime)\|(1+|x|)^l\nabla\cdot\partial_x^\alpha v\|^2
\\[2mm]
&+(\mu+\mu^\prime)\langle\nabla\cdot\partial_x^\alpha v, 2l(1+|x|)^{2l-1}\frac{x}{|x|}\partial_x^\alpha v\rangle+\langle\nabla\partial_x^\alpha \sigma, (1+|x|)^{2l}\partial_x^\alpha v\rangle\\[2mm]
&-\kappa\gamma_1\langle\nabla\Delta\partial_x^\alpha \sigma, (1+|x|)^{2l}\partial_x^\alpha v\rangle-\kappa\gamma_2\langle\nabla\Delta\partial_x^\alpha \vartheta, (1+|x|)^{2l}\partial_x^\alpha v\rangle\\[2mm]
&=\langle\nabla\Delta\partial_x^\alpha \{-(b_1\cdot\nabla)c_1+\tilde{f}\}, (1+|x|)^{2l}\partial_x^\alpha v\rangle
\end{array}
\end{equation}
Applying $\partial_x^\alpha$ to $(\ref{2.3})_1$, then taking the $L^2$ inner product with $(1+|x|)^{2l}\partial_x^\alpha v$ on the resultant equation, we have from integration by parts that
\begin{equation}\label{2.58}
\begin{array}{ll}
&-\langle\partial_x^\alpha v, (1+|x|)^{2l}\nabla\partial_x^\alpha \sigma\rangle
-\langle\partial_x^\alpha v, 2l(1+|x|)^{2l-1}\frac{x}{|x|}\partial_x^\alpha \sigma\rangle\\[2mm]
&\hspace{30mm}+\langle\nabla\partial_x^\alpha ((a\cdot\nabla)\sigma), (1+|x|)^{2l}\partial_x^\alpha \sigma\rangle=\langle\partial_x^\alpha g,  (1+|x|)^{2l}\partial_x^\alpha \sigma\rangle
\end{array}
\end{equation}
Canceling the term $-\langle\partial_x^\alpha v, (1+|x|)^{2l}\nabla\partial_x^\alpha \sigma\rangle$ by adding (\ref{2.58}) to (\ref{2.57}), and taking summation with respect to $\alpha$, we obtain
\begin{equation}\label{2.59}
\begin{array}{ll}
\displaystyle\left\|(1+|x|)^l\nabla^{l+1}v\right\|^2&\leq C\left\{\Big(\left|\langle\nabla^{l+1}v, (1+|x|)^{2l}\nabla^lv\rangle\right|
+\left|\langle\nabla^{l}v, (1+|x|)^{2l-1}\nabla^{l}\sigma\rangle\right|\right.\\[2mm]
&\quad+\left|\langle \nabla^lg,  (1+|x|)^{2l}\nabla^{l}\sigma\rangle\right|+\left|\langle \nabla^{l}\tilde{f},  (1+|x|)^{2l}\nabla^{l}v\rangle\right|\Big)\\[2mm]
&\quad+\left|\langle \nabla^{l}\{(a\cdot\nabla)\sigma\},  (1+|x|)^{2l}\nabla^{l}\sigma\rangle\right|
+\left|\langle \nabla^{l}\{(b_1\cdot\nabla) c_1\}, (1+|x|)^{2l}\nabla^{l}v\rangle\right|\\[2mm]
&\quad+\left|\langle \nabla\Delta\nabla^{l}\sigma, (1+|x|)^{2l}\nabla^{l}v\rangle\right|+\left|\langle \nabla\Delta\nabla^{l}\vartheta, (1+|x|)^{2l}\nabla^{l}v\rangle\right|\Big\}
\\[2mm]&=C\{I_4+I_5+I_6+I_7+I_8\}
\end{array}
\end{equation}
Integration by parts and the Cauchy inequality imply that
\begin{equation}\label{2.60}
\begin{array}{ll}
I_4&\leq\eta\left\|(1+|x|)^l(\nabla^l\sigma, \nabla^{l+1}v)\right\|^2\\[2mm]&
\quad+C_\eta\left\{\|(1+|x|)^{l-1}\nabla^lv\|^2+\|(1+|x|)^{l}\nabla^lg\|^2+\|(1+|x|)^{l}\nabla^{l-1}\tilde{f}\|^2\right\}.
\end{array}
\end{equation}
Similar to the estimate of $I_{2,1}$, we have
\begin{equation}\label{2.61}
I_6\leq\eta\left\|(1+|x|)^l\nabla^{l+1}v\right\|^2+C_\eta\left\{\|(1+|x|)^{l-1}\nabla^lv\|^2+\|b_1\|^2_{J^{k+1}}\|c_1\|^2_{J^{k+1}}\right\}.
\end{equation}
For $I_7$, we deduce from integration by parts and $(\ref{2.3})_1$ that
\begin{equation}\label{2.62}
\begin{array}{ll}
I_7&\leq\left|\langle \Delta\nabla^{l}\sigma, 2l(1+|x|)^{2l-1}\frac{x}{|x|}\cdot\nabla^{l}v\rangle\right|+\left|\langle \Delta\nabla^{l}\sigma, (1+|x|)^{2l}\nabla^{l}(\nabla\cdot v)\rangle\right|\\[2mm]
&\leq2l\left|\langle \nabla^{l+2}\sigma, (1+|x|)^{2l-1}|\nabla^{l}v|\rangle\right|+\left|\langle \nabla^{l+2}\sigma, (1+|x|)^{2l}\nabla^{l}\{(a\cdot\nabla)\sigma\}\rangle\right|\\[2mm]
&\quad+\left|\langle \nabla^{l+2}\sigma, (1+|x|)^{2l-1}\nabla^{l}g|\rangle\right|\\[2mm]
&\leq\eta\left\|(1+|x|)^l\nabla^{l+2}\sigma\right\|^2+C_\eta\left(\|(1+|x|)^{l-1}\nabla^lv\|^2+\|(1+|x|)^{l}\nabla^lg\|^2\right)+A
\end{array}
\end{equation}
where
\begin{equation}\label{2.63}
\begin{array}{ll}
A&=\left|\langle \nabla^{l+2}\sigma, (1+|x|)^{2l}\nabla^{l}\{(a\cdot\nabla)\sigma\}\rangle\right|\\[2mm]
&=\displaystyle\sum_{|\alpha|=l}\left|\langle \nabla^{l+2}\sigma, (1+|x|)^{2l}\partial_x^\alpha\{(a\cdot\nabla)\sigma\}\rangle\right|\\[2mm]
&\leq\displaystyle\sum_{|\alpha|=l}\left|\langle \nabla^{l+2}\sigma, (1+|x|)^{2l}(a\cdot\nabla)\partial_x^\alpha\sigma\rangle\right|
\\[2mm]&\quad+\displaystyle\sum_{|\alpha|=l}\sum_{\beta<\alpha}C_\alpha^\beta\left|\langle \nabla^{l+2}\sigma, (1+|x|)^{2l}\partial_x^{\alpha-\beta}a\cdot\partial_x^\beta\nabla\sigma)\rangle\right|=A_1+A_2.
\end{array}
\end{equation}
The Sobolev inequality and the Cauchy inequality imply that
\begin{equation}\label{2.64}
\begin{array}{ll}
A_1&\leq C\|a\|_{L^\infty}\left\|(1+|x|)^{l}\nabla^{l+2}\sigma\right\|\left\|(1+|x|)^{l}\nabla^{l+1}\sigma\right\|\\[2mm]
&\leq C\delta\left\|(1+|x|)^{l}(\nabla^{l+1}\sigma, \nabla^{l+2}\sigma)\right\|^2
\end{array}
\end{equation}
\begin{equation}\label{2.65}
\begin{array}{ll}
A_2&\leq \displaystyle\sum_{|\alpha|=l}\bigg\{\displaystyle\sum_{\beta<\alpha,\,|\alpha-\beta|\leq\frac{|\alpha|}{2}}+\displaystyle\sum_{\beta<\alpha,\,|\alpha-\beta|>\frac{|\alpha|}{2}}\bigg\}
C_\alpha^\beta\left|\langle \nabla^{l+2}\sigma, (1+|x|)^{2l}\partial_x^{\alpha-\beta}a\cdot\partial_x^\beta\nabla\sigma)\rangle\right|\\[3mm]
&\leq\displaystyle\sum_{|\alpha|=l}\displaystyle\sum_{\beta<\alpha,\,|\alpha-\beta|\leq\frac{|\alpha|}{2}}C_\alpha^\beta\left\|(1+|x|)^{|\alpha-\beta|-1}\partial_x^\alpha a\right\|_{L^\infty}\left\|(1+|x|)^{|\beta|+1}\nabla^{|\beta|+1}\sigma\right\|\left\|(1+|x|)^{l}\nabla^{l+2}\sigma\right\|\\[3mm]
&\quad+\displaystyle\sum_{|\alpha|=l}\displaystyle\sum_{\beta<\alpha,\,|\alpha-\beta|>\frac{|\alpha|}{2}}C_\alpha^\beta\left\|(1+|x|)^{|\alpha-\beta|-1}\partial_x^\alpha a\right\|\left\|(1+|x|)^{|\beta|+1}\nabla^{|\beta|+1}\sigma\right\|_{L^\infty}\left\|(1+|x|)^{l}\nabla^{l+2}\sigma\right\|\\[3mm]
&\leq C\delta\displaystyle\sum_{\nu=l}^l\left\|(1+|x|)^l(\nabla^l\sigma, \nabla^{l+1}\sigma, \nabla^{l+2}\sigma)\right\|^2.
\end{array}
\end{equation}
Combining (\ref{2.62})-(\ref{2.65}), we obtain
\begin{equation}\label{2.66}
\begin{array}{ll}
I_7&\leq\eta\left\|(1+|x|)^l\nabla^{l+2}\sigma\right\|^2+C_\eta\left(\left\|(1+|x|)^{l-1}\nabla^lv\right\|^2+\left\|(1+|x|)^{l}\nabla^lg\right\|^2\right)\\[2mm]
&\quad+C\delta\displaystyle\sum_{\nu=1}^l\left\|(1+|x|)^l(\nabla^l\sigma, \nabla^{l+1}\sigma, \nabla^{l+2}\sigma)\right\|^2.
\end{array}
\end{equation}
Finally, similar to the estimate of $A$ and $I_7$, respectively,  we have
\begin{equation}\label{2.67}
I_5\leq C\delta\displaystyle\sum_{\nu=1}^l\left\|(1+|x|)^l(\nabla^l\sigma, \nabla^{l+1}\sigma, \nabla^{l+2}\sigma)\right\|^2.
\end{equation}
and
\begin{equation}\label{2.68}
\begin{array}{ll}
I_8&\leq\eta\left\|(1+|x|)^l\nabla^{l+2}\vartheta\right\|^2+C_\eta\left(\left\|(1+|x|)^{l-1}\nabla^lv\right\|^2+\left\|(1+|x|)^{l}\nabla^lg\right\|^2\right)\\[2mm]
&\quad+C\delta\displaystyle\sum_{\nu=1}^l\left\|(1+|x|)^l(\nabla^l\sigma, \nabla^{l+1}\sigma, \nabla^{l+2}\vartheta)\right\|^2.
\end{array}
\end{equation}
Substituting (\ref{2.60}), (\ref{2.61}), (\ref{2.66})-(\ref{2.68}) into (\ref{2.59}), by the smallness of $\eta$, we arrive at
\begin{equation}\label{2.69}
\begin{array}{ll}
\displaystyle\left\|(1+|x|)^l\nabla^{l+1}v\right\|^2&\leq C\left\{\eta\left\|(1+|x|)^l\nabla^l\sigma\right\|^2+(\eta+\delta)\left\|(1+|x|)^{l+2}\nabla^{l+2}(\sigma,\vartheta)\right\|^2\right.\\[2mm]
&\quad+C_\eta\Big(\left\|(1+|x|)^l\nabla^lg\right\|^2+\left\|(1+|x|)^{l-1}\nabla^lv\right\|^2+\left\|(1+|x|)^l\nabla^{l-1}\tilde{f}\right\|^2\Big)
\\[2mm]&\quad+\displaystyle\delta\sum_{\nu=1}^l\left\|(1+|x|)^\nu(\nabla^\nu\sigma, \nabla^{\nu+1}\sigma)\right\|^2\Big\}
\end{array}
\end{equation}

{\bf Step 3.~~}For any multi-index $\alpha$ with $|\alpha|=l$, applying $\partial_x^\alpha$ to $(\ref{2.3})_3$, then taking the $L^2$ inner product with $(1+|x|)^{2l}\partial_x^\alpha \vartheta$ on the resultant equation,  integrating by parts and summing up $\alpha$, we can get
\begin{equation}\label{2.70}
\begin{array}{ll}
\tilde{\alpha}\left\|(1+|x|)^l\nabla^{l+1}\vartheta\right\|^2&\leq C\left|\langle \nabla^{l+1}\vartheta, (1+|x|)^{2l-1}\nabla^{l}\vartheta|\rangle\right|+C\left|\langle \nabla^{l}\{-(b_2\cdot\nabla)c_2+\tilde{h}\}, (1+|x|)^{2l}\nabla^{l}\vartheta\rangle\right|\\[2mm]
&=I_9+I_{10}.
\end{array}
\end{equation}
For $I_9$, the Cauchy inequality imply that
\begin{equation}\label{2.71}
I_9\leq\frac{\tilde{\alpha}}{4}\left\|(1+|x|)^l\nabla^{l+1}\vartheta\right\|^2+C\left\|(1+|x|)^{l-1}\nabla^{l}\vartheta\right\|^2.
\end{equation}
Similar to the estimate of $I_7$, $I_{10}$ can be estimated as follows
\begin{equation}\label{2.72}
\begin{array}{ll}
I_{10}&\leq\displaystyle\frac{\tilde{\alpha}}{4}\left\|(1+|x|)^l\nabla^{l+1}\vartheta\right\|^2
\\[2mm]&\qquad\qquad\qquad\quad+C\left\{\left\|(1+|x|)^{l-1}\nabla^{l}\vartheta\right\|^2+\|b_2\|^2_{J^{5}}\|c_2\|^2_{N^{5}}+\left\|(1+|x|)^l\nabla^{l-1}\tilde{h}\right\|^2\right\}.
\end{array}
\end{equation}
Putting (\ref{2.71}) and (\ref{2.72}) into (\ref{2.70}) gives
\begin{equation}\label{2.73}
\left\|(1+|x|)^l\nabla^{l+1}\vartheta\right\|^2\leq C\left\{\left\|(1+|x|)^{l-1}\nabla^{l}\vartheta\right\|^2+\|b_2\|^2_{J^{5}}\|c_2\|^2_{N^{5}}
+\left\|(1+|x|)^l\nabla^{l-1}\tilde{h}\right\|^2\right\}.
\end{equation}
On the other hand, we also get from $(\ref{2.3})_3$ that
\begin{equation}\label{2.74}
\left\|(1+|x|)^l\nabla^{l+2}\vartheta\right\|^2\leq C\left\|(1+|x|)^{l}\nabla^{l}h\right\|^2\leq C\left(\|b_2\|^2_{J^{5}}\|c_2\|^2_{N^{5}}
+\left\|(1+|x|)^l\nabla^{l}\tilde{h}\right\|^2\right).
\end{equation}
Consequently, we deduce from (\ref{2.73}) and (\ref{2.74}) that
\begin{equation}\label{2.75}
\begin{array}{ll}
&\left\|(1+|x|)^l(\nabla^{l+1}\vartheta, \nabla^{l+2}\vartheta)\right\|^2\\[2mm]&\leq C\left\{\left\|(1+|x|)^{l-1}\nabla^{l}\vartheta\right\|^2+\|b_2\|^2_{J^{5}}\|c_2\|^2_{N^{5}}
+\left\|(1+|x|)^l(\nabla^{l-1}\tilde{h}, \nabla^{l}\tilde{h})\right\|^2\right\}.
\end{array}
\end{equation}

{\bf Step 4.~} Now, we begin to prove (\ref{2.40}) by using the estimates in the above three steps. We use the method of induction. First, for the case of $l=1$, we derive from
(\ref{2.43}), (\ref{2.44}), (\ref{2.46}) and (\ref{2.54}) that
\begin{equation}\label{2.76}
\begin{array}{ll}
\displaystyle\left\|(1+|x|)(\nabla\sigma, \nabla^{2}\sigma, \nabla^{3}\sigma)\right\|^2
&\leq C\left\{\left\|(1+|x|)\nabla^{2}v\right\|^2+\left\|(1+|x|)\tilde{f}\right\|^2+\left\|(1+|x|)\nabla^3\vartheta\right\|^2\right.
\\[2mm]&\left.\quad\quad+\left\|\nabla\sigma\right\|^2+\|b_1\|^2_{J^{5}}\|c_1\|^2_{J^{5}}\right\}.
\end{array}
\end{equation}
which, together with (\ref{2.69}) and (\ref{2.75}) with $l=1$ gives
\begin{equation}\label{2.77}
\begin{array}{ll}
&\displaystyle\left\|(1+|x|)(\nabla\sigma, \nabla^{2}\sigma, \nabla^{3}\sigma)\right\|^2+\left\|(1+|x|)\nabla^2v\right\|^2+\left\|(1+|x|)(\nabla^{2}\vartheta, \nabla^{3}\vartheta)\right\|^2\\[2mm]
&\leq C\left\{\left\|\nabla(\sigma, v, \vartheta)\right\|^2+\|b_1\|^2_{J^{5}}\|c_1\|^2_{J^{5}}+\|b_2\|^2_{J^{5}}\|c_2\|^2_{N^{5}}
+\left\|(1+|x|)(\tilde{f}, \tilde{h}, \nabla g, \nabla\tilde{h})\right\|^2\right\}.
\end{array}
\end{equation}
by the smallness of $\eta$ and $\delta$. Thus, we assume for $l\geq2$ that
\begin{equation}\label{2.78}
\begin{array}{ll}
&\displaystyle\left\|(1+|x|)^{l-1}(\nabla^{l-1}\sigma, \nabla^{l}\sigma, \nabla^{l+1}\sigma)\right\|^2+\left\|(1+|x|)^{l-1}\nabla^lv\right\|^2+\left\|(1+|x|)^{l-1}(\nabla^{l}\vartheta, \nabla^{l+1}\vartheta)\right\|^2\\[2mm]
&\hspace{10mm}\leq C\left\{\left\|\nabla(\sigma, v, \vartheta)\right\|^2+\|b_1\|^2_{J^{5}}\|c_1\|^2_{J^{5}}+\|b_2\|^2_{J^{5}}\|c_2\|^2_{N^{5}}\right.\\[2mm]
&\displaystyle\hspace{20mm}+\sum_{\nu=1}^{l-1}\left\|(1+|x|)^\nu(\nabla^{\nu-1}\tilde{f}, \nabla^{\nu}\tilde{h}, \nabla^{\nu}g, \nabla^{\nu-1}\tilde{h})\right\|^2\Big\}.
\end{array}
\end{equation}
Furthermore, the linear combination $[M_1\times(\ref{2.69})+(\ref{2.56})]+M_2\times(\ref{2.75})$ for $M_1>0$ and $M_2>0$ large enough in turn gives
\begin{equation}\label{2.79}
\begin{array}{ll}
&\displaystyle\left\|(1+|x|)^l\left(\nabla^l\sigma, \nabla^{l+1}\sigma, \nabla^{l+2}\sigma\right)\right\|^2+\left\|(1+|x|)^l\left(\nabla^{l+1}v, \nabla^{l+1}\vartheta, \nabla^{l+2}\vartheta\right)\right\|^2\\[2mm]
&\leq C\left\{\left\|(1+|x|)^{l-1}(\nabla^{l-1}\sigma, \nabla^{l}\sigma, \nabla^{l}v, \nabla^{l}\vartheta)\right\|^2+\|b_1\|^2_{J^{5}}\|c_1\|^2_{J^{5}}+\|b_2\|^2_{J^{5}}\|c_2\|^2_{N^{5}}\right.\\[2mm]
&\quad+\displaystyle\left\|(1+|x|)^l\nabla^{l-1}(\tilde{f}, \tilde{h})\right\|+\left\|(1+|x|)^l\nabla^{l}(g, \tilde{h})\right\|+\delta\sum_{\nu=1}^{l-1}\left\|(1+|x|)^\nu\left(\nabla^\nu\sigma, \nabla^{\nu+1}\sigma\right)\right\|^2\Big\}
\end{array}
\end{equation}
provided that $\eta$ and $\delta$ are small enough. Combining (\ref{2.78}) with (\ref{2.79}), if $\delta>0$ is small enough, we can get  (\ref{2.40}). This completes the proof of Lemma 2.3.

Combining Proposition 2.2 and Lemma 2.3, we have the following theorem.
\begin{Theorem}
There exists $\delta_0=\delta_0(\gamma_1, \gamma_2, \kappa, \mu, \mu^\prime, \tilde{\alpha})>0$ such that such that if $\delta$ in (\ref{2.4}) satisfies $\delta\leq\delta_0$, then (\ref{2.3}) admits a solution $(\sigma, v, \vartheta)\in\hat{\mathcal{H}}^{6, 5, 6}$ which satisfies the estimate:
\begin{equation}\label{2.80}
\begin{array}{ll}
&\|(\sigma, v, \vartheta)\|_{L^6}+\displaystyle\sum_{\nu=1}^4\left\|(1+|x|)^\nu\left(\nabla^\nu\sigma, \nabla^{\nu+1}\sigma, \nabla^{\nu+2}\sigma\right)\right\|+\displaystyle\sum_{\nu=1}^{5}\left\|(1+|x|)^{\nu-1}\nabla^{\nu}v\right\|\\[2mm]
&+\displaystyle\sum_{\nu=1}^{5}\left\|(1+|x|)^{\nu-1}\left(\nabla^{\nu}\vartheta, \nabla^{\nu+1}\vartheta\right)\right\|\\[2mm]
&\quad \leq C\left\{\|b_1\|_{J^{5}}\|c_1\|_{J^{5}}+\|b_2\|_{J^{5}}\|c_2\|_{N^{5}}+\|(1+|x|)(g,\tilde{h})\|\right.\\[2mm]
&\qquad \quad\displaystyle+\sum_{\nu=1}^{4}\left\|(1+|x|)^\nu\nabla^{\nu}(g, \tilde{h})\right\|+\displaystyle\sum_{\nu=0}^{3}\displaystyle\left\|(1+|x|)^{\nu+1}\nabla^{\nu}(\tilde{f}, \tilde{h})\right\|\Big\}
\end{array}
\end{equation}
where $C>0$ is a constant depending only on $\gamma_1, \gamma_2, \kappa, \mu, \mu^\prime$ and $\tilde{\alpha}$.
\end{Theorem}

\subsection{Proof of Theorem 1.1}
In this subsection, we shall construct a solution to (\ref{1.6}) by the contraction mapping principle in $\dot{\Lambda}_{\epsilon}^{4, 5, 5}$. To this end, we consider the following iteration system
\begin{eqnarray}\label{2.81}
\left\{\begin{array}{ll}
\nabla\cdot v\displaystyle+\frac{\tilde{\rho}_P}{\tilde{\rho}}(\tilde{v}\cdot\nabla)\sigma=g,\\[2mm]
\displaystyle-\mu\Delta v-(\mu+\mu^\prime)\nabla(\nabla\cdot v)+\nabla \sigma-\kappa\gamma_1\nabla\Delta\sigma-\kappa\gamma_2\nabla\Delta\vartheta=-\bar{\rho}(\tilde{v}\cdot\nabla)\tilde{v}+\tilde{f},\\[2mm]
\displaystyle-\tilde{\alpha}\Delta\vartheta=-\bar{\eta}_1(\tilde{v}\cdot\nabla)\tilde{\vartheta}+\tilde{h},
 \end{array}\right.
\end{eqnarray}
where
\begin{eqnarray}\label{2.82}
\left\{\begin{array}{ll}
g=\displaystyle-\frac{\tilde{\rho}_\theta}{\tilde{\rho}}(\tilde{v}\cdot\nabla) \tilde{\vartheta}+\frac{G(x)}{\tilde{\rho}},\\[2mm]
\tilde{f}=-(\tilde{\rho}-\bar{\rho})(\tilde{v}\cdot\nabla)\tilde{v}+\kappa\tilde{\rho}\left(\nabla\tilde{\sigma}\cdot\nabla^2\tilde{\rho}_P+\nabla\tilde{\rho}_P\cdot\nabla^2\tilde{\sigma}+\nabla\tilde{\rho}_P\Delta\tilde{\sigma}\right)
+\kappa\left(\tilde{\rho}\tilde{\rho}_P-\bar{\rho}\bar{\rho}_P
\right)\nabla\Delta\tilde{\sigma}\\[2mm]
\qquad+\kappa\tilde{\rho}\left(\nabla\tilde{\vartheta}\cdot\nabla^2\tilde{\rho}_\theta+\nabla\tilde{\rho}_\theta\cdot\nabla^2\tilde{\vartheta}+\nabla\tilde{\rho}_\theta\Delta\tilde{\vartheta}\right)
+\kappa\left(\tilde{\rho}\tilde{\rho}_\theta-\bar{\rho}\bar{\rho}_\theta
\right)\nabla\Delta\tilde{\vartheta}+\tilde{\rho}F-\tilde{v}G,\\[2mm]
\tilde{h}=-(\tilde{\eta}_1-\bar{\eta}_1)(\tilde{v}\cdot\nabla)\tilde{\vartheta}+\tilde{\eta}_2(\tilde{v}\cdot\nabla)\tilde{\sigma}\displaystyle+\Psi(\tilde{v})-\tilde{\eta}_3G+H
+\Phi(\tilde{\rho}, \tilde{v})+\frac{\tilde{v}^2}{2}G-C_\triangledown \tilde{\theta}G,\\[2mm]
\displaystyle\tilde{\eta}_1=\displaystyle\tilde{\rho} C_\triangledown-\frac{\tilde{\theta}\tilde{\rho}_\theta^2}{\tilde{\rho}\tilde{\rho}_P},\quad\tilde{\eta}_2=\displaystyle\frac{\tilde{\theta}\tilde{\rho}_\theta}{\tilde{\rho}},\quad
\tilde{\eta}_3=\displaystyle\frac{\tilde{\theta}\tilde{\rho}_\theta}{\tilde{\rho}\tilde{\rho}_P},\quad \tilde{\theta}=\bar{\theta}+\tilde{\vartheta}.\\[2mm]
\end{array}\right.
\end{eqnarray}
Here, $(\tilde{\sigma}, \tilde{v}, \tilde{\vartheta})\in\dot{\Lambda}_{\epsilon}^{4, 5, 5}$ is given,  and $\tilde{\rho}_P=\rho_P(\bar{P}+\tilde{\sigma}, \bar{\theta}+\tilde{\vartheta}), \bar{\eta}_1=\eta_1(\bar{P}, \bar{\theta})$, etc.

\subsubsection{Introduction of solution map $T$ for (\ref{2.1})}
Firstly, we apply Theorem 2.1 to (\ref{2.81}) to get the weighted $L^2$ estimate. Let
\begin{eqnarray}\label{2.83}
a=\displaystyle-\frac{\tilde{\rho}_P}{\tilde{\rho}}\tilde{v},\quad b_1=c_1={\bar{\rho}}^{\frac{1}{2}}\tilde{v},\quad b_2=\bar{\eta}_1\tilde{v},\quad c_2=\tilde{\vartheta},
\end{eqnarray}
and $g, \tilde{f}, \tilde{h}$ in Theorem 2.1 be defined as in (\ref{2.82}). We choose $\epsilon>0$ sufficiently small such that $\frac{\bar{\rho}}{2}<\tilde{\rho}<2\bar{\rho}$, as follows from the sobolev inequality. Assume that the assumptions of Theorem 2.1 hold and denote
\begin{equation}\label{2.84}
K_0=\displaystyle\sum_{\nu=0}^{3}\left\|(1+|x|)^{\nu+1}\nabla^\nu(G, F, H)\right\|+\left\|(1+|x|)^{4}\nabla^4(G, H)\right\|<\infty
\end{equation}
then we can check (\ref{2.5}) and (\ref{2.6}) easily and additionally we have
\[\|(1+|x|)(g, \tilde{h})\|+\displaystyle\sum_{\nu=1}^4(1+|x|)^\nu\nabla^\nu(g,\tilde{h})\|+\displaystyle\sum_{\nu=0}^{3}(1+|x|)^{\nu+1}\nabla^\nu(\tilde{f}, \tilde{h})\|\leq C\left(\epsilon^2+K_0\right).
\]
for some constant $C=C(\bar{\rho}, \bar{\theta}, \mu, \mu^\prime, \kappa)>0$. Applying Theorem 2.1 to (\ref{2.81}), we have the following lemma.
\begin{Lemma}
Let $(G, F, H)\in H^{4,3,4}$ satisfy (\ref{2.84}). Then there exists a constant $\epsilon_0>0$ such that if $\epsilon\leq\epsilon_0$, (\ref{2.81}) with $(\tilde{\sigma}, \tilde{v}, \tilde{\vartheta})\in\dot{\Lambda}_{\epsilon}^{6, 5, 6}$ admits a solution $(\sigma, v, \vartheta)\in\hat{\mathcal{H}}^{6,5,6}$ which satisfies the estimate:
\begin{equation}\label{2.85}
\begin{array}{ll}
&\|(\sigma, v, \vartheta)\|_{L^6}+\displaystyle\sum_{\nu=1}^4\left\|(1+|x|)^\nu\left(\nabla^\nu\sigma, \nabla^{\nu+1}\sigma, \nabla^{\nu+2}\sigma\right)\right\|+\displaystyle\sum_{\nu=1}^{5}\left\|(1+|x|)^{\nu-1}\nabla^{\nu}v\right\|\\[3mm]
&\quad+\displaystyle\sum_{\nu=1}^{5}\left\|(1+|x|)^{\nu-1}\left(\nabla^{\nu}\vartheta, \nabla^{\nu+1}\vartheta\right)\right\|\leq C\left(\epsilon^2+K_0\right)
\end{array}
\end{equation}
where the constant $C$ depends only on $\bar{\rho}, \bar{\theta}, \mu, \mu^\prime, \kappa$ and $\tilde{\alpha}$.
\end{Lemma}
Based on Lemma 2.4, we can define the solution map $T: \dot{\Lambda}_{\epsilon}^{4, 5, 5}\rightarrow\hat{\mathcal{H}}^{6, 5, 6}$ by $(\sigma, v, \vartheta)=T(\tilde{\sigma}, \tilde{v}, \tilde{\vartheta})$. Since the contraction mapping principle will be applied to prove Theorem 1.1, we have to show that $T(\tilde{\sigma}, \tilde{v}, \tilde{\vartheta})=(\sigma, v, \vartheta)\in\dot{\Lambda}_{\epsilon}^{4, 5, 5}$. To this end, we first cite the following lemma which will play an important role when we estimate the solution by the $L^\infty$ norm.
\begin{Lemma}(\cite{Y. Shibata-K. Tanaka-2003})
Let $E(x)$ be a scalar function satisfying
\[\left|\partial_x^\alpha E(x)\right|\leq\frac{C_\alpha}{|x|^{|\alpha|+1}},\quad|\alpha|=0, 1, 2.
\]
(i) If $\phi(x)$ is a smooth scalar function of the form $\phi=\nabla\cdot\phi_1+\phi_2$ satisfying
\[L_1(\phi)\equiv\left\|(1+|x|)^3\phi\right\|_{L^\infty}+\left\|(1+|x|)^2\phi_1\right\|_{L^\infty}+\|\phi_2\|_{L^1}<\infty,\]
then we have for any multi-index $\alpha$ with $|\alpha|=0,1$
\[\left|\partial_x^\alpha (E\ast\phi)(x)\right|\leq\frac{C_\alpha}{|x|^{|\alpha|+1}}L_1(\phi).
\]
(ii) If $\phi(x)$ is a smooth scalar function of the form $\phi=\phi_1\phi_2$ satisfying
\[L_2(\phi)\equiv\left\|(1+|x|)^2\phi\right\|_{L^\infty}+\left\|(1+|x|)^3(\nabla\phi_1)\phi_2\right\|_{L^\infty}
+\left\|(1+|x|)^3\phi_1(\nabla\phi_2)\right\|_{1}<\infty,\]
then we have for any multi-index $\alpha$ with $|\alpha|=1, 2$
\[\left|\partial_x^\alpha (E\ast\phi)(x)\right|\leq\frac{C_\alpha}{|x|^{|\alpha|}}L_2(\phi).
\]
Here $C_\alpha$ denotes a constant depending only on $\alpha$.
\end{Lemma}

With the aid of the Helmholtz decomposition and the Fourier transform, the solution of (\ref{2.81}) can be formulated as follows, cf. \cite{Y. Shibata-K. Tanaka-2003}.
\begin{equation}\label{2.86}
v=w+\nabla p,\quad \displaystyle\sigma-\kappa\gamma_1\Delta\sigma=\Phi+\kappa\gamma_2\Delta\vartheta+(2\mu+\mu^\prime)\Delta p,\quad\vartheta=E\ast\Theta,
\end{equation}
where
\begin{eqnarray}\label{2.87}
\left\{\begin{array}{ll}
w_j(x)=\displaystyle\sum_{i=1}^3E_{ij}\ast f_i(x),\\[2mm]
p(x)=E_0\ast R(x),\\[2mm]
\Phi=E_0\ast(\nabla\cdot f).
 \end{array}\right.
\end{eqnarray}
and
\begin{eqnarray}\label{2.88}
\left\{\begin{array}{ll}
E_{ij}(x)=\displaystyle\frac{1}{8\pi\mu}\left(\frac{\delta_{ij}}{|x|}-\frac{x_ix_j}{|x|^3}\right),\quad E_0=-\frac{1}{4\pi|x|}\\[3mm]
f_i=-\bar{\rho}(\tilde{v}\cdot\nabla)\tilde{v}_i+\tilde{f}_i,\\[2mm]
R(x)=\displaystyle-\frac{\tilde{\rho}_\theta}{\tilde{\rho}}(\tilde{v}\cdot\nabla)\sigma+g,\\[2mm]
\Theta=\displaystyle\frac{1}{\tilde{\alpha}}\{\bar{\eta}_1(\tilde{v}\cdot\nabla)\tilde{\vartheta}-\tilde{h}\}.
 \end{array}\right.
\end{eqnarray}
Now, we shall estimate the $L^\infty$ norm of the solution to (\ref{2.81}) by using Lemma 2.5.
\begin{Lemma}
Let $(G, F, H)\in H^{4,3,4}$ satisfy the following estimate:
\begin{equation}\label{2.89}
K\equiv\|(G, F, H)\|_{\mathcal{L}}+\left\|(1+|x|)^4\nabla^4(G, H)\right\|<\infty
\end{equation}
If $(\sigma, v, \vartheta)\in\hat{\mathcal{H}}^{6,5,6}$ is a solution to (\ref{2.81}) with
$(\tilde{\sigma}, \tilde{v}, \tilde{\vartheta})\in\dot{\Lambda}_{\epsilon}^{4, 5, 5}$ and satisfies (\ref{2.85}), then $(\sigma, v, \vartheta)$ satisfies the estimate:
\begin{equation}\label{2.90}
\begin{array}{ll}
&\displaystyle\sum_{\nu=0}^1\left\|(1+|x|)^2\nabla^\nu\sigma\right\|_{L^\infty}+\displaystyle\sum_{\nu=0}^1\left\|(1+|x|)^{\nu+1}\nabla^\nu(v, \vartheta)\right\|_{L^\infty}\\[4mm]
&\qquad+\left\|(1+|x|)^2\nabla^2(v, \vartheta)\right\|_{L^\infty}\leq C\left(\epsilon^2+K\right)
\end{array}
\end{equation}
where the constant $C>0$ depends only on $\bar{\rho}, \bar{\theta}, \mu, \mu^\prime, \kappa$ and $\tilde{\alpha}$.
\end{Lemma}
\noindent{\bf Proof.} First, we deduce an estimate on $f$. Since $(\tilde{\sigma}, \tilde{v}, \tilde{\vartheta})\in\dot{\Lambda}_{\epsilon}^{4, 5, 5}$, there exits $\tilde{V}_1=(\tilde{V}^i_1)_{1\leq i\leq3}$ and $\tilde{V}_2$ such that $\nabla\cdot\tilde{v}=\nabla\cdot\tilde{V}_1+\tilde{V}_2$, and
\[
\left\|(1+|x|)^3\tilde{V}_1\right\|_{L^\infty}+\left\|(1+|x|)^{-1}\tilde{V}_2\right\|_{L^1}\leq\epsilon
\]
Thus we have
\[
\begin{array}{ll}
f_i&=-\tilde{\rho}(\tilde{v}\cdot\nabla)\tilde{v}_i+\kappa\tilde{\rho}[\nabla\tilde{\sigma}\cdot\nabla^2\tilde{\rho}_P+\nabla\tilde{\rho}_P\cdot\nabla^2\tilde{\sigma}+\nabla\tilde{\rho}_P\Delta\tilde{\sigma}]_i+\kappa\left(\tilde{\rho}\tilde{\rho}_\theta-\bar{\rho}\bar{\rho}_\theta
\right)\Delta\tilde{\theta}_{x_i}\\[2mm]
&\quad+\kappa\tilde{\rho}[\nabla\tilde{\vartheta}\cdot\nabla^2\tilde{\rho}_\theta+\nabla\tilde{\rho}_\theta\cdot\nabla^2\tilde{\vartheta}+\nabla\tilde{\rho}_\theta\Delta\tilde{\vartheta}]_i
+\kappa\left(\tilde{\rho}\tilde{\rho}_P-\bar{\rho}\bar{\rho}_P
\right)\Delta\tilde{\sigma}_{x_i}+\tilde{\rho}F_i-\tilde{v}_iG\\[2mm]
&=\nabla\cdot\left(-\tilde{\rho}\tilde{v}_i\tilde{v}+\tilde{\rho}\tilde{v}_i\tilde{V}_1+\rho F_{1,i}\right)
+\left\{-\tilde{\rho}(\tilde{V}_1\cdot\nabla)\tilde{v}_i-\tilde{v}_i(\tilde{V}_1\cdot\nabla)\tilde{\rho}+\tilde{\rho}\tilde{v}_i\tilde{V}_2
+\tilde{\rho} F_{2,i}\right.\\[2mm]
&\quad+\kappa\tilde{\rho}[\nabla\tilde{\vartheta}\cdot\nabla^2\tilde{\rho}_\theta+\nabla\tilde{\rho}_\theta\cdot\nabla^2\tilde{\vartheta}+\nabla\tilde{\rho}_\theta\Delta\tilde{\vartheta}]_i
+\kappa\tilde{\rho}[\nabla\tilde{\sigma}\cdot\nabla^2\tilde{\rho}_P+\nabla\tilde{\rho}_P\cdot\nabla^2\tilde{\sigma}+\nabla\tilde{\rho}_P\Delta\tilde{\sigma}]_i
\\[2mm]&\left.\quad
+\tilde{v}_i(\tilde{v}\cdot\nabla)\tilde{\rho}-\nabla\tilde{\rho}\cdot F_{1,i}+\kappa\left(\tilde{\rho}\tilde{\rho}_P-\bar{\rho}\bar{\rho}_P
\right)\Delta\tilde{\sigma}_{x_i}-\tilde{v}_iG+\kappa\left(\tilde{\rho}\tilde{\rho}_\theta-\bar{\rho}\bar{\rho}_\theta
\right)\Delta\tilde{\vartheta}_{x_i}\right\}\\[2mm]
&=\nabla\cdot f_{i,1}+f_{i,2}.
\end{array}
\]
Here $[\cdots ]_i$ denotes the $i-th$ component of the vector $[\cdots]$.

By $(\tilde{\sigma}, \tilde{v}, \tilde{\vartheta})\in\dot{\Lambda}_{\epsilon}^{4, 5, 5}$ and (\ref{2.85}), using the Sobolev inequality and mean value theorem, we have
\[
\left\|(1+|x|)^3f_i\right\|_{L^\infty}+\left\|(1+|x|)^2f_{1, i}\right\|_{L^\infty}+\left\|f_{2,i}\right\|_{L^1}\leq C(\epsilon^2+K_1)
\]
and
\[
\left\|(1+|x|)^3\nabla f_i\right\|_{L^\infty}+\left\|(1+|x|)^2f_{i}\right\|_{L^\infty}\leq C(\epsilon^2+K_1)
\]
where \[K_1=\left\|(1+|x|)^3(F, G, \nabla F, \nabla G)\right\|_{L^\infty}+\left\|(1+|x|)^2F_1\right\|_{L^\infty}+\|F_2\|_{L^1}.\]
Hence, by (i) and (ii) of Lemma 2.5, we obtain
\begin{equation}\label{2.91}
|x||w_j|,\,\, |x|^2\left(|\Phi|, |\nabla\Phi|, |\nabla w_j|, |\nabla^2 w_j|\right)\leq C(\epsilon^2+K_1).
\end{equation}
As for $\nabla p, \nabla^2p, \nabla^3p$, due to \cite{J. Z. Qian-H. Yin-2007},
\begin{equation}\label{2.92}
|x||\nabla p|,\,\, |x|^2\left(|\nabla^2p|, |\nabla^3p|\right)\leq C(\epsilon^2+K_0+K_2).
\end{equation}
where \[K_2=\left\|(1+|x|)^2G\right\|_{L^\infty}+\left\|(1+|x|)^3(\nabla G, \nabla^2 G)\right\|_{L^\infty}\]
Combining $(\ref{2.86})_1$, (\ref{2.91}) and (\ref{2.92}) yields
\begin{equation}\label{2.93}
|x||v|,\,\, |x|^2\left(|\nabla v|, |\nabla^2v|\right)\leq C(\epsilon^2+K_0+K_1+K_2).
\end{equation}
Next, we turn to estimate $\vartheta$. To this end, we rewrite $\Theta$ as
\[
\begin{array}{ll}
\Theta&=\displaystyle\frac{1}{\tilde{\alpha}}\left\{\tilde{\eta}_1(\tilde{v}\cdot\nabla)\tilde{\vartheta}-\tilde{\eta}_2(\tilde{v}\cdot\nabla)\tilde{\sigma}\displaystyle-\Psi(\tilde{v})+\tilde{\eta}_3G-H
-\Phi(\tilde{\rho}, \tilde{v})-\frac{\tilde{v}^2}{2}G+C_\triangledown (\tilde{\vartheta}+\bar{\theta})G\right\}\\[2mm]
&=\displaystyle\frac{1}{\tilde{\alpha}}\nabla\cdot\left\{(\tilde{\eta}_1\tilde{\vartheta}-\tilde{\eta}_2\tilde{\sigma})(\tilde{v}-\tilde{V}_1)
+\tilde{\eta}_3G_1-H_1-\frac{\tilde{v}^2}{2}G_1+C_\triangledown \tilde{\vartheta}G_1\right\}\\[2mm]
&\quad+\displaystyle\frac{1}{\tilde{\alpha}}\nabla\cdot\left\{(\tilde{V}_1\cdot\nabla)(\eta_1\tilde{\vartheta}-\eta_2\tilde{\sigma})
-(\tilde{\eta}_1\tilde{\vartheta}-\tilde{\eta}_2\tilde{\sigma})\tilde{V}_2-\nabla\tilde{\eta}_1\tilde{v}\tilde{\vartheta}
+\nabla\tilde{\eta}_2\tilde{v}\tilde{\sigma}-\Psi(\tilde{v})-\Phi(\tilde{\rho}, \tilde{v})\right.\\[2mm]&
\quad\left.-\nabla\tilde{\eta}_3\cdot G_1+\tilde{\eta}_3G_2+H_2+
\nabla\cdot\tilde{v}\tilde{v}\cdot G_1-\displaystyle\frac{\tilde{v}^2}{2}G_2+C_\triangledown\nabla\tilde{\vartheta}\cdot G_1-C_\triangledown G_2(\tilde{\vartheta}+\bar{\theta})\right\}
\\[2mm]
&=\nabla\cdot \Theta_1+\Theta_2.
\end{array}
\]
and
\[
\begin{array}{ll}
\Theta&=\displaystyle\sum_{i=1}^3(-\frac{1}{\tilde{\alpha}}\tilde{\eta}_2\tilde{v}_i)\tilde{\sigma}_{x_i}
+\frac{1}{\tilde{\alpha}}\left\{\tilde{\eta}_1(\tilde{v}\cdot\nabla)\tilde{\vartheta}\displaystyle-\Psi(\tilde{v})+\tilde{\eta}_3G-H
-\Phi(\tilde{\rho}, \tilde{v})-\frac{\tilde{v}^2}{2}G+C_\triangledown (\tilde{\vartheta}+\bar{\theta})G\right\}
\\[2mm]
&=\displaystyle\sum_{i=1}^3\Theta^i_1\Theta^i_2+\Theta_3
\end{array}
\]
Since $(\tilde{\sigma}, \tilde{v}, \tilde{\vartheta})\in\dot{\Lambda}_{\epsilon}^{4, 5, 5}$, it follows from (\ref{2.85}) and the Sobolev inequality that
\[
\left\|(1+|x|)^3\Theta\right\|_{L^\infty}+\left\|(1+|x|)^2\Theta_1\right\|_{L^\infty}+\left\|\Theta_2\right\|_{L^1}\leq C(\epsilon^2+K_3),
\]
\[
\left\|(1+|x|)^3\nabla\Theta_3\right\|_{L^\infty}+\left\|(1+|x|)^2\Theta_3\right\|_{L^\infty}\leq C(\epsilon^2+K_3),
\]
\[
\left\|(1+|x|)^3\Theta^i_1\Theta^i_2\right\|_{L^\infty}+\left\|(1+|x|)^3(\nabla\Theta^i_1)\Theta^i_2\right\|_{L^\infty}
+\left\|(1+|x|)^3\Theta^i_1(\nabla\Theta^i_2)\right\|_{L^\infty}\leq C\epsilon^2,
\]
where \[K_3=\left\|(1+|x|)^3(G, H, \nabla G, \nabla H)\right\|_{L^\infty}+\left\|(1+|x|)^2(G_1, H_1)\right\|_{L^\infty}+\|(G_2, H_2)\|_{L^1}.\]
Thus, it follows from (i) and (ii) of Lemma 2.5 that
\begin{equation}\label{2.94}
|x||\vartheta|,\,\, |x|^2\left(|\nabla \vartheta|, |\nabla^2\vartheta|\right)\leq C(\epsilon^2+K_3).
\end{equation}
Finally, for the estimate of $\vartheta$, taking the Fourier transform on both side of $(\ref{2.86})_2$, we have
\begin{equation}\label{2.95}
\begin{array}{ll}
\sigma(x)&=\displaystyle\frac{1}{(2\pi)^{3/2}}\mathcal{F}^{-1}\left(\frac{1}{1+\kappa\gamma_1|\xi|^2}\right)\ast L(\Phi, p, \vartheta)\\[3mm]
&=\displaystyle\frac{1}{(4\pi\kappa\gamma_1)^{3/2}}\int_{\mathbb{R}^3}\left(\int_0^\infty e^{-t-\frac{|y|^2}{4\gamma_1t}}t^{-\frac{3}{2}}\,dt\right)
(\Phi+\kappa\gamma_2\Delta\vartheta+(2\mu+\mu^\prime)\Delta p)(x-y)\,dy
\end{array}
\end{equation}
where $L(\Phi, p, \vartheta)=\Phi+\kappa\gamma_2\Delta\vartheta+(2\mu+\mu^\prime)\Delta p$ and we have used the fact that (cf. \cite{L. C. Evans-1998})
\begin{equation}\label{2.96}
\displaystyle\mathcal{F}^{-1}\left(\frac{1}{1+\kappa\gamma_1|\xi|^2}\right)=\displaystyle\frac{1}{(2\kappa\gamma_1)^{3/2}}\int_0^\infty e^{-t-\frac{|y|^2}{4\kappa\gamma_1t}}t^{-\frac{3}{2}}\,dt
\end{equation}
Note that the right hand of (\ref{2.96}) is the so-called Bessel potential. We deduce from (\ref{2.91}), (\ref{2.92}) and (\ref{2.94}) that
\begin{equation}\label{2.97}
\begin{array}{ll}
|x|^2|\sigma(x)|&\leq\displaystyle\frac{|x|^2}{(4\pi\kappa\gamma_1)^{3/2}}\int_{\mathbb{R}^3}\left(\int_0^\infty e^{-t-\frac{|y|^2}{4\kappa\gamma_1t}}t^{-\frac{3}{2}}\,dt\right)\frac{1}{|x-y|^2}\,dy\\[3mm]
&\quad\times\displaystyle\sup_{(x-y)\in \mathbb{R}^3}\left\{|x-y|^2\left|\Phi+\kappa\gamma_2\Delta\vartheta+(2\mu+\mu^\prime)\Delta p\right|(x-y)\right\}\\[3mm]
&\leq C(\epsilon^2+K)\cdot A
\end{array}
\end{equation}
where $K$ is defined by $(\ref{2.89})$ and
\begin{equation}\label{2.98}
\begin{array}{ll}
A&=\displaystyle\int_{\mathbb{R}^3}\left(\int_0^\infty e^{-t-\frac{|y|^2}{4\kappa\gamma_1t}}t^{-\frac{3}{2}}\,dt\right)\frac{|x|^2}{|x-y|^2}\,dy\\[3mm]
&=\displaystyle\int_0^\infty e^{-t}t^{-\frac{3}{2}}\left(\int_{\mathbb{R}^3}\frac{|x|^2}{|x-y|^2}e^{-t-\frac{|y|^2}{4\kappa\gamma_1t}}\,dy\right)dt\\[3mm]
&\leq\displaystyle2\int_0^\infty e^{-t}t^{-\frac{3}{2}}\left(\int_{\mathbb{R}^3}e^{-\frac{|y|^2}{4\kappa\gamma_1t}}\,dy\right)dt+2\int_0^\infty e^{-t}t^{-\frac{3}{2}}\left(\int_{\mathbb{R}^3}\frac{|y|^2}{|x-y|^2}e^{-t-\frac{|y|^2}{4\kappa\gamma_1t}}\,dy\right)dt\\[3mm]
&\leq\displaystyle16(\kappa\gamma_1)^{3/2}\int_0^\infty e^{-t}\,dt\int_{\mathbb{R}^3}e^{-|z|^2}\,dz+2\int_0^\infty e^{-t}t^{-\frac{3}{2}}\left(\int_{B(x,\,1)}\frac{|y|^2}{|x-y|^2}e^{-\frac{|y|^2}{4\kappa\gamma_1t}}\,dy\right)dt\\[3mm]
&\quad\displaystyle+2\int_0^\infty e^{-t}t^{-\frac{3}{2}}\left(\int_{\mathbb{R}^3\setminus B(x,\,1)}\frac{|y|^2}{|x-y|^2}e^{-\frac{|y|^2}{4\kappa\gamma_1t}}\,dy\right)dt\\[3mm]
&\leq \displaystyle C+8\kappa\gamma_1\int_0^\infty e^{-t}t^{-\frac{1}{2}}\left(\int_{B(x,\,1)}\frac{1}{|x-y|^2}\,dy\right)dt\\[3mm]
&\displaystyle\quad+2\int_0^\infty e^{-t}t^{-\frac{3}{2}}\left(\int_{\mathbb{R}^3\setminus B(x,\,1)}|y|^2e^{-\frac{|y|^2}{4\kappa\gamma_1t}}\,dy\right)dt\\[3mm]
&\leq C
\end{array}
\end{equation}
Here $B(x,\,1)$ denotes the unit ball in $\mathbb{R}^3$. Consequently, it follows from (\ref{2.97}) and (\ref{2.98}) that
\begin{equation}\label{2.99}
|x|^2|\sigma(x)|\leq C(\epsilon^2+K)
\end{equation}
Differentiating the equation (\ref{2.95}) and notice that
\[
\left\||x|^2|\nabla^3\vartheta|\right\|_{L^\infty}\leq C\left\|\nabla(|x|^2|\nabla^3\vartheta|)\right\|_1\leq C(\epsilon^2+K_0),
\]
by using the same argument as above, we can also obtain
\begin{equation}\label{2.100}
|x|^2|\nabla\sigma(x)|\leq C(\epsilon^2+K)
\end{equation}
Next, we consider the case of $|x|<1$. The Sobolev inequality and (\ref{2.85}) imply that
\begin{equation}\label{2.101}
\begin{array}{ll}
&\|(\sigma, v, \vartheta)\|_{L^\infty}\leq C\|\nabla(\sigma, v, \vartheta)\|_1\leq C(\epsilon^2+K_0)\\[2mm]
&\|\nabla^\nu(\sigma, v, \vartheta)\|_{L^\infty}\leq C\|\nabla^{\nu+1}(\sigma, v, \vartheta)\|_1\leq C(\epsilon^2+K_0),\quad \nu=1,\,2.
\end{array}
\end{equation}
(\ref{2.90}) thus follows from (\ref{2.93}), (\ref{2.94}), (\ref{2.99}), (\ref{2.100}) and (\ref{2.101}). This completes the proof of Lemma 2.6.

In the following Proposition, we show  that the solution $(\sigma, v, \vartheta)\in\dot{\Lambda}_{\epsilon}^{4, 5, 5}$.
\begin{Proposition}
There exits $c_0>0$ such that for any sufficiently small constant $\epsilon>0$, if $(G, F, H)\in\mathcal{H}^{4,3,4}$ satisfies
\[K+\left\|(1+|x|)^{-1}G\right\|_{L^1}\leq c_0\epsilon\quad(\mbox{$K$ is defined in Lemma 2.6})\]
then (\ref{2.81}) with $(\tilde{\sigma}, \tilde{v}, \tilde{\vartheta})\in\dot{\Lambda}_{\epsilon}^{4, 5, 5}$ admits a solution $(\sigma, v, \vartheta)=T(\tilde{\sigma}, \tilde{v}, \tilde{\vartheta})\in\dot{\Lambda}_{\epsilon}^{4, 5, 5}$.
\end{Proposition}
\noindent{\bf Proof.} By Lemmas 2.4 and 2.6, it follows that (\ref{2.81}) has a solution $(\sigma, v, \vartheta)\in\hat{\mathcal{H}}^{4,5,5}$ , which satisfies
\[\left\|(\sigma, v, \vartheta)\right\|_{{\Lambda}^{4, 5, 5}}\leq C(\epsilon^2+K)\leq C(\epsilon^2+c_0\epsilon),\]
where the constant $C>0$ depends only on $\bar{\rho}, \bar{\theta}, \mu, \mu^\prime, \kappa$ and $\tilde{\alpha}$. Thus if we take $c_0\leq\frac{1}{2C}$ and $\epsilon>0$ is small enough, it follows that $(\sigma, v, \vartheta)\in\Lambda_{\epsilon}^{4, 5, 5}$. Finally, we define $V_1$ and $V_2$ by
\[V_1=-\frac{\tilde{\rho}_P}{\tilde{\rho}}\tilde{v}\sigma,\quad V_2=\nabla\cdot\left(\frac{\tilde{\rho}_P}{\tilde{\rho}}\tilde{v}\right)\sigma-\frac{\tilde{\rho}_\theta}{\tilde{\rho}}\tilde{v}\cdot\nabla\vartheta+\frac{G}{\tilde{\rho}}
\]
Then it follows from $(\ref{2.81})_1$ that
\[\nabla\cdot v=\nabla\cdot V_1+V_2.
\]
Moreover, by $(\tilde{\sigma}, \tilde{v}, \tilde{\vartheta})\in\dot{\Lambda}_{\epsilon}^{4, 5, 5}$ , (\ref{2.85}) and (\ref{2.90}), we have form the Sobolev inequality that
\[
\begin{array}{ll}
\left\|(1+|x|)^3V_1\right\|_{L^\infty}+\left\|(1+|x|)^{-1}V_2\right\|_{L^1}&\leq C\left\{\epsilon^2+K+\left\|(1+|x|)^{-1}G\right\|_{L^1}\right\}\\[2mm]
&\leq C(\epsilon^2+c_0\epsilon)\leq\epsilon
\end{array}
\]
if $c_0\leq\frac{1}{2C}$ and $\epsilon>0$ is sufficiently small. This completes the proof of Proposition 2.3.

\subsubsection{Contraction of the solution map $T$ }
In this subsection, we shall show that the solution map $T$ for (\ref{2.81}) is contractive. Suppose that  $(\tilde{\sigma}^j, \tilde{v}^j, \tilde{\vartheta}^j)\in\dot{\Lambda}_{\epsilon}^{4, 5, 5}$ and $(\sigma^j, v^j, \vartheta^j)=T(\tilde{\sigma}^j, \tilde{v}^j, \tilde{\vartheta}^j)$ for $j=1,2$, then we deduce from (\ref{2.81}) that
\begin{eqnarray}\label{2.102}
\left\{\begin{array}{ll}
\nabla\cdot (v^1-v^2)\displaystyle+\frac{\tilde{\rho}^1_P}{\tilde{\rho}^1}(\tilde{v}^1\cdot\nabla)(\sigma^1-\sigma^2)=g,\\[2mm]
\displaystyle-\mu\Delta (v^1-v^2)-(\mu+\mu^\prime)\nabla(\nabla\cdot (v^1-v^2))+\nabla (\sigma^1-\sigma^2)-\kappa\gamma_1\nabla\Delta(\sigma^1-\sigma^2)\\
-\kappa\gamma_2\nabla\Delta(\vartheta^1-\vartheta^2)=-\tilde{\rho}^2(\tilde{v}^1-\tilde{v}^2)\cdot\nabla\tilde{v}^1-\tilde{\rho}^2(\tilde{v}^2\cdot\nabla)(\tilde{v}^1-\tilde{v}^2)+\tilde{f},\\[2mm]
\displaystyle-\tilde{\alpha}\Delta(\vartheta^1-\vartheta^2)=-\tilde{\eta}_1^2[(\tilde{v}^1-\tilde{v}^2)\cdot\nabla\tilde{\vartheta}^1+(\tilde{v}^2\cdot\nabla)(\tilde{\vartheta}^1-\tilde{\vartheta}^2)]+\tilde{h},
 \end{array}\right.
\end{eqnarray}
where
\begin{eqnarray}\label{2.103}
\left\{\begin{array}{ll}
g=\displaystyle-\left(\frac{\tilde{\rho}_P^1}{\tilde{\rho}^1}-\frac{\tilde{\rho}_P^2}{\tilde{\rho}^2}\tilde{v}^2\right)\cdot\nabla\sigma^2
-\left(\frac{\tilde{\rho}_\theta^1}{\tilde{\rho}^1}\tilde{v}^1\cdot\nabla\tilde{\vartheta}^1-\frac{\tilde{\rho}_\theta^2}{\tilde{\rho}^2}\tilde{v}^2\cdot\nabla\tilde{\vartheta}^2\right) +\left(\frac{1}{\tilde{\rho}^1}-\frac{1}{\tilde{\rho}^2}\right)G,\\[3mm]
\tilde{f}=-(\tilde{\rho}^1-\tilde{\rho}^2)(\tilde{v}^1\cdot\nabla)\tilde{v}^1-(\tilde{\rho}^1-\tilde{\rho}^2)F-(\tilde{v}^1-\tilde{v}^2)G
+\big\{\kappa\tilde{\rho}^1\left(\nabla\tilde{\sigma}^1\cdot\nabla^2\tilde{\rho}^1_P+\right.
\\[2mm]\qquad\left.+\nabla\tilde{\rho}^1_P\cdot\nabla^2\tilde{\sigma}^1+\nabla\tilde{\rho}^1_P\Delta\tilde{\sigma}^1\right)
-\kappa\tilde{\rho}^2\left(\nabla\tilde{\sigma}^2\cdot\nabla^2\tilde{\rho}^2_P+\nabla\tilde{\rho}^2_P\cdot\nabla^2\tilde{\sigma}^2+\nabla\tilde{\rho}^2_P\Delta\tilde{\sigma}^2\right)\big\}
\\[2mm]\qquad+\kappa\Big\{\tilde{\rho}^1\big(\nabla\tilde{\vartheta}^1\cdot\nabla^2\tilde{\rho}^1_\theta+\nabla\tilde{\rho}^1_\theta\cdot\nabla^2\tilde{\vartheta}^1+\nabla\tilde{\rho}^1_\theta\Delta\tilde{\vartheta}^1\big)
-\tilde{\rho}^2\big(\nabla\tilde{\vartheta}^2\cdot\nabla^2\tilde{\rho}^2_\theta+\nabla\tilde{\rho}^2_\theta\cdot\nabla^2\tilde{\vartheta}^2
\\[2mm]\qquad+\nabla\tilde{\rho}^2_\theta\Delta\tilde{\vartheta}^2\big)
\Big\}
+\kappa\left\{\left(\tilde{\rho}^1\tilde{\rho}^1_\theta-\bar{\rho}\bar{\rho}_\theta\right)\nabla\Delta\tilde{\vartheta}^1-\left(\tilde{\rho}^2\tilde{\rho}^2_\theta-\bar{\rho}\bar{\rho}_\theta\right)\nabla\Delta\tilde{\vartheta}^2\right\}
\\[2mm]\qquad+\kappa\left\{\left(\tilde{\rho}^1\tilde{\rho}^1_P-\bar{\rho}\bar{\rho}_P
\right)\nabla\Delta\tilde{\sigma}^1-\left(\tilde{\rho}^2\tilde{\rho}^2_P-\bar{\rho}\bar{\rho}_P
\right)\nabla\Delta\tilde{\sigma}^2\right\},\\[2mm]
\tilde{h}=-(\tilde{\eta}_1^1-\tilde{\eta}^2_1)(\tilde{v}^1\cdot\nabla)\tilde{\vartheta}^1+(\tilde{\eta}_1^1-\tilde{\eta}^2_2)(\tilde{v}^1\cdot\nabla)\tilde{\sigma}^1+\tilde{\eta}^2_2\big((\tilde{v}^1\cdot\nabla)\tilde{\sigma}^1
-(\tilde{v}^2\cdot\nabla)\tilde{\sigma}^2\big)\\[2mm]
\qquad+\Psi(\tilde{v}^1)-\Psi(\tilde{v}^2)+\Phi(\tilde{\rho}^1, \tilde{v}^1)-\Phi(\tilde{\rho}^2, \tilde{v}^2)
+\frac{1}{2}(\tilde{v}^1+\tilde{v}^2)(\tilde{v}^1-\tilde{v}^2)G\\
[2mm]\qquad-\left(\tilde{\eta}^1_3-\tilde{\eta}^2_3\right)G-C_\triangledown(\tilde{\vartheta}^1-\tilde{\vartheta}^2)G,\\[2mm]
\displaystyle\tilde{\eta}_1=\displaystyle\tilde{\rho}^j C_\triangledown-\frac{\tilde{\theta}^j(\tilde{\rho}^j_\theta)^2}{\tilde{\rho}^j\tilde{\rho}^j_P},\quad\tilde{\eta}^j_2=\displaystyle\frac{\tilde{\theta}^j\tilde{\rho}^j_\theta}{\tilde{\rho}^j},\quad
\tilde{\eta}^j_3=\displaystyle\frac{\tilde{\theta}^j(\tilde{\rho}^j_\theta)^2}{\tilde{\rho}^j\tilde{\rho}^j_P},\quad \tilde{\theta}^j=\bar{\theta}+\tilde{\vartheta}^j,\,\, j=1,2.\\[2mm]
\end{array}\right.
\end{eqnarray}
Since
\begin{equation}\label{2.104}
\begin{array}{ll}
&\left\|(1+|x|)(g,\tilde{h})\right\|\displaystyle+\sum_{\nu=1}^{4}\left\|(1+|x|)^\nu\nabla^{\nu}(g, \tilde{h})\right\|+\displaystyle\sum_{\nu=0}^{3}\displaystyle\left\|(1+|x|)^{\nu+1}\nabla^{\nu}(\tilde{f}, \tilde{h})\right\|\\[2mm]
&\quad \leq C(\epsilon+K)\left\|(\tilde{\sigma}^1-\tilde{\sigma}^2, \tilde{v}^1-\tilde{v}^2, \tilde{\vartheta}^1-\tilde{\vartheta}^2)\right\|_{\Lambda^{4,5,5}}
\end{array}
\end{equation}
as follows from the Sobolev inequality for $K$ defined in (\ref{2.89}). Applying Theorem 2.1 to (\ref{2.102}), we obtain
\begin{equation}\label{2.105}
\begin{array}{ll}
&\|(\sigma^1-\sigma^2, v^1-v^2, \vartheta^1-\vartheta^2)\|_{L^6}+\displaystyle\sum_{\nu=1}^{5}\left\|(1+|x|)^{\nu-1}\nabla^{\nu}(v^1-v^2)\right\|
\\[2mm]&\quad\quad\quad+\displaystyle\sum_{\nu=1}^4\left\|(1+|x|)^\nu\left(\nabla^\nu(\sigma^1-\sigma^2), \nabla^{\nu+1}(\sigma^1-\sigma^2), \nabla^{\nu+2}(\sigma^1-\sigma^2)\right)\right\|\\[2mm]&
\qquad\quad\quad\quad+\displaystyle\sum_{\nu=1}^{5}\left\|(1+|x|)^{\nu-1}\left(\nabla^{\nu}(\vartheta^1-\vartheta^2), \nabla^{\nu+1}(\vartheta^1-\vartheta^2)\right)\right\|\\[2mm]
&\quad\quad\quad\qquad\qquad\leq C(\epsilon+K)\left\|(\tilde{\sigma}^1-\tilde{\sigma}^2, \tilde{v}^1-\tilde{v}^2, \tilde{\vartheta}^1-\tilde{\vartheta}^2)\right\|_{\Lambda^{4,5,5}}
\end{array}
\end{equation}
Similarly, by the same argument as in the proof Lemma 2.6, we can get
\begin{equation}\label{2.106}
\begin{array}{ll}
&\displaystyle\sum_{\nu=0}^{1}\left\|(1+|x|)^{2}\nabla^{\nu}(\sigma^1-\sigma^2)\right\|_{L^\infty}
+\displaystyle\sum_{\nu=0}^{1}\left\|(1+|x|)^{\nu+1}\nabla^{\nu}(v^1-v^2,\,\vartheta^1-\vartheta^2)\right\|_{L^\infty}\\[4mm]
&\displaystyle\quad+\left\|(1+|x|)^2\nabla^2(v^1-v^2,\,\vartheta^1-\vartheta^2)\right\|_{L^\infty}\\[2mm]&
\quad\quad \leq C(\epsilon+K)\left\|(\tilde{\sigma}^1-\tilde{\sigma}^2, \tilde{v}^1-\tilde{v}^2, \tilde{\vartheta}^1-\tilde{\vartheta}^2)\right\|_{\Lambda^{4,5,5}}\\[2mm]&
\quad\quad\quad\quad +C\epsilon\left\{\left\|(1+|x|)^{3}(\tilde{V}_1^1-\tilde{V}_1^2)\right\|_{L^\infty}+\left\|(1+|x|)^{-1}(\tilde{V}_2^1-\tilde{V}_2^2)\right\|_{L^1}\right\}
\end{array}
\end{equation}
where $\tilde{V}^j_1, \tilde{V}^j_2,\, j=1,2$ are functions satisfying
\begin{equation}\label{2.107}
\nabla\cdot\tilde{v}^j=\nabla\cdot\tilde{V}^j_1+\tilde{V}^j_2,\quad \left\|(1+|x|)^{3}\tilde{V}_1^j\right\|_{L^\infty}+\left\|(1+|x|)^{-1}\tilde{V}_2^j\right\|_{L^1}\leq\epsilon
\end{equation}
Moreover, if we define $\tilde{V}^j_1, \tilde{V}^j_2,\, j=1,2$ as
\begin{equation}\label{2.108}
V_1^j=-\frac{\tilde{\rho}^j_P}{\tilde{\rho}^j}\tilde{v}^j\sigma^j,\quad V_2^j=\nabla\cdot\left(\frac{\tilde{\rho}^j_P}{\tilde{\rho}^j}\tilde{v}^j\right)\sigma^j-\frac{\tilde{\rho}^j_\theta}{\tilde{\rho}^j}\tilde{v}^j\cdot\nabla\vartheta^j+\frac{G}{\tilde{\rho}^j}
\end{equation}
then we get from $(\ref{2.102})_1$ that
\[\nabla\cdot(v^1-v^2)=\nabla\cdot(V_1^1-V_1^2)+V_2^1-V_2^2\]
and
\begin{equation}\label{2.109}
\begin{array}{ll} &\left\|(1+|x|)^{3}(V_1^1-V_1^2)\right\|_{L^\infty}+\left\|(1+|x|)^{-1}(V_2^1-V_2^2)\right\|_{L^1}\\[2mm]
&\quad\leq C\left(\epsilon+\left\|(1+|x|)^{-1}G\right\|_{L^1}\right)\left\|(\tilde{\sigma}^1-\tilde{\sigma}^2, \tilde{v}^1-\tilde{v}^2,\tilde{\vartheta}^1-\tilde{\vartheta}^2)\right\|_{\Lambda^{4,5,5}}
\end{array}
\end{equation}
Combining (\ref{2.105})-(\ref{2.109}), we obtain
\begin{equation}\label{2.110}
\begin{array}{ll}
&\left\|(\sigma^1-\sigma^2, v^1-v^2, \vartheta^1-\vartheta^2)\right\|_{\Lambda^{4,5,5}}\\[2mm]&\qquad+\left\|(1+|x|)^{3}(V_1^1-V_1^2)\right\|_{L^\infty}+\left\|(1+|x|)^{-1}(V_2^1-V_2^2)\right\|_{L^1}\\[2mm]
&\quad\leq C\left(\epsilon+K\right)\left\|(\tilde{\sigma}^1-\tilde{\sigma}^2, \tilde{v}^1-\tilde{v}^2,\tilde{\vartheta}^1-\tilde{\vartheta}^2)\right\|_{\Lambda^{4,5,5}}\\[2mm]
&\qquad\quad+C\epsilon\left\{\left\|(1+|x|)^{3}(\tilde{V}_1^1-\tilde{V}_1^2)\right\|_{L^\infty}+\left\|(1+|x|)^{-1}(\tilde{V}_2^1-\tilde{V}_2^2)\right\|_{L^1}\right\}
\end{array}
\end{equation}
Therefore, we have the following Proposition.
\begin{Proposition}
There exits a constant $c_0>0$ such that for any sufficiently small constant $\epsilon>0$, if $(G, F, H)\in\mathcal{H}^{4,3,4}$ satisfies
\[K+\left\|(1+|x|)^{-1}G\right\|_{L^1}\leq c_0\epsilon\quad(\mbox{$K$ is defined in Lemma 2.6}),\]
the for $(\tilde{\sigma}^j, \tilde{v}^j, \tilde{\vartheta}^j)\in\dot{\Lambda}_{\epsilon}^{4, 5, 5}$ and $(\sigma^j, v^j, \vartheta^j)=T(\tilde{\sigma}^j, \tilde{v}^j, \tilde{\vartheta}^j)$,\,j=1,2, we have the following estimates
\begin{equation}\label{2.111}
\begin{array}{ll}
&\left\|(\sigma^1-\sigma^2, v^1-v^2, \vartheta^1-\vartheta^2)\right\|_{\Lambda^{4,5,5}}\\[2mm]&\qquad+\left\|(1+|x|)^{3}(V_1^1-V_1^2)\right\|_{L^\infty}+\left\|(1+|x|)^{-1}(V_2^1-V_2^2)\right\|_{L^1}\\[2mm]
&\quad\leq \displaystyle\frac{1}{2}\left\{\left\|(\tilde{\sigma}^1-\tilde{\sigma}^2, \tilde{v}^1-\tilde{v}^2,\tilde{\vartheta}^1-\tilde{\vartheta}^2)\right\|_{\Lambda^{4,5,5}}\right.\\[2mm]
&\left.\qquad\quad+\left\|(1+|x|)^{3}(\tilde{V}_1^1-\tilde{V}_1^2)\right\|_{L^\infty}+\left\|(1+|x|)^{-1}(\tilde{V}_2^1-\tilde{V}_2^2)\right\|_{L^1}\right\}
\end{array}
\end{equation}
where $(\tilde{V}_1^j, \tilde{V}_2^j),\,j=1,2$ satisfy (\ref{2.107}) and $(V_1^j, V_2^j),\,j=1,2$ are defined by (\ref{2.109}).
\end{Proposition}
 Hence, by Propositions 2.3 and 2.4, the contraction mapping principle implies the existence and uniqueness of solution to (\ref{1.6}). This completes the proof of Theorem 1.1.

\section{Non-stationary problem}
\setcounter{equation}{0}
In this section, we consider the stability of the Stationary solution with respect to the initial disturbance $(\rho_0, v_0, \vartheta_0)$ . Fix $\bar{\rho}, \bar{\theta}$ to be positive constants and let $F, G, H$ be small in the sense of Theorem 1.1. We denote the corresponding stationary solution obtained in Theorem 1.1 by $(P^*, v^*, \theta^*)$ , and set $\rho^*\equiv\bar{\rho}+\sigma^*=\rho(P^*, \theta^*)$. Then by direct calculations, we have the following estimate for $\sigma^*$:
\[
\begin{array}{ll}
\|\sigma^*\|_{N^5}&\equiv\displaystyle\sum_{\nu=1}^5\left\|(1+|x|)^{\nu-1}(\nabla^\nu\sigma^*,\nabla^{\nu+1}\sigma^*)\right\|+\sum_{\nu=0}^1\displaystyle
\left\|(1+|x|)^{\nu+1}\nabla^\nu\sigma^*\right\|_{L^\infty}+\left\|(1+|x|)^2\nabla^2\sigma^*\right\|_{L^\infty}\\[2mm]
&\leq C\epsilon,
\end{array}
\]
where the constant $C>0$ is depending only on $\bar{\rho}$ and $\bar{\theta}$. Thus, we have
\[\left\|(\sigma^*, v^*, \theta^*)\right\|_{\mathcal{F}^{5,5,5}}\equiv\|\sigma^*\|_{N^5}+\|v^*\|_{J^5}+\|\vartheta^*\|_{N^5}\leq (C+1)\epsilon.\]
For simplicity, we assume in this section that $\left\|(\sigma^*, v^*, \theta^*)\right\|_{\mathcal{F}^{5,5,5}}\leq\epsilon$ for $\epsilon$ sufficiently small. Define the new variables
\[\sigma(t,x)=\rho(t,x)-\rho^*,\quad w(t,x)=v(t,x)-v^*,\quad \vartheta(t,x)=\theta(t,x)-\theta^*,\]
then the initial value problem (\ref{1.3}), (\ref{1.4}) is reformulated as
\begin{eqnarray}\label{3.1}
\left\{\begin{array}{ll}
         \sigma_t(t)+\nabla\cdot\{(\rho^*+\sigma(t))w(t)\}=-\nabla\cdot\left(v^*\sigma(t)\right),\\[2mm]
          w(t)_t-\frac{1}{\rho^*}\left[\mu\Delta w(t)+(\mu+\mu^\prime)\nabla(\nabla\cdot w(t))\right]+A(t)\nabla\sigma(t)-\kappa\nabla\Delta\sigma(t)+B(t)\nabla\vartheta(t)=f(t),\\[2mm]
          \vartheta_t(t)-\tilde{\alpha}D^*\Delta\vartheta(t)+E(t)\nabla\cdot w(t)=h(t),
 \end{array}\right.
\end{eqnarray}
with initial date
\begin{equation}\label{3.2}
(\sigma, w, \vartheta)(t, x)|_{t=0}=(\sigma_0, w_0, \vartheta_0)(x)\equiv(\sigma-\sigma^*, v-v^*, \theta-\theta^*)(0,x).
\end{equation}
where
\[
\begin{array}{ll}
f(t)&=-(v^*\cdot\nabla)w(t)-(w(t)\cdot\nabla)(v^*+w(t))-\left(A(t)-A^*\right)\nabla\rho^*-\left(B(t)-B^*\right)\nabla\vartheta^*\\[2mm]
&\displaystyle\quad-\left(\left(\frac{v}{\rho}\right)(t)-\frac{v^*}{\rho^*}\right)G-\frac{\sigma(t)}{\rho^*(\rho^*+\sigma(t))}\left[\mu\Delta(v^*+w(t))+(\mu+\mu^\prime)\nabla
\left(\nabla\cdot(v^*+w(t))\right)\right],\\[3mm]
h(t)&=-(v^*\cdot\nabla)\vartheta(t)-\left(w(t)\cdot\nabla\right)(\theta^*+\vartheta(t))+\tilde{\alpha}\left(D(t)-D^*\right)\Delta(\theta^*+\vartheta(t))+\left(D(t)-D^*\right)H
\\[2mm]
&\quad+\left(D(t)-D^*\right)\left(\Psi(v^*)+\Phi(\rho^*,v^*)\right)+D(t)\left[\Psi(v)(t)+\Phi(\rho, v)(t)-\Psi(v^*)-\Phi(\rho^*, v^*)\right]
\\[2mm]
&\quad\displaystyle+\frac{1}{2}\left[D(t)v^2(t)-D^*v^{*2}\right]G-C_\triangledown\left[D(t)\theta(t)-D^*\theta^*\right]G-(E(t)-E^*)\nabla\cdot v^*,
\end{array}
\]
and \[A(t)=\frac{P_\rho(\rho, \theta)}{\rho},\quad B(t)=\frac{P_\theta(\rho, \theta)}{\rho},\quad D(t)=\frac{1}{C_\triangledown\rho},\quad E(t)=\frac{\theta P_\theta(\rho, \theta)}{C_\triangledown\rho}\]
with $A^*=A(\rho^*, \theta^*), A(t)=A\left(\rho^*+\sigma(t), \theta^*+\vartheta(t)\right),\,etc.$ Moreover, we set $A_i(t), B_i(t), E_i(t), i=1,2$ and  $D_1(t)$ to be functions satisfying: $A(t)-A^*=A_1(t)\sigma(t)+A_2(t)\vartheta(t), B(t)-B^*=B_1(t)\sigma(t)+B_2(t)\vartheta(t)$, $E(t)-E^*=E_1(t)\sigma(t)+E_2(t)\vartheta(t)$ and $D(t)-D^*=D_1(t)\sigma(t)$, respectively.

The aim of this section is to prove Theorem 1.2. The proof consists of the following two steps. The first one is the local existence result:
\begin{Proposition}
Suppose that $(\sigma, w, \vartheta)(0)\in\mathcal{H}^{4,3,3}$. Then there exists a constant $t_0>0$ such that the initial value problem (\ref{3.1})-(\ref{3.2}) admits a unique solution $(\sigma, w, \vartheta)(0)\in\mathcal{C}(0,t_0; \mathcal{H}^{4,3,3})$. Moreover, $(\sigma, w, \vartheta)(t)$ satisfies
\[\left\|(\sigma, w, \vartheta)(t)\right\|^2_{4,3,3}\leq2\left\|(\sigma, w, \vartheta)(0)\right\|^2_{4,3,3}\]
for any $t\in[0,t_0]$.
\end{Proposition}
And the other is an a priori estimate:
\begin{Proposition}
Let $(\sigma, w, \vartheta)(0)\in\mathcal{C}(0, t_1; \mathcal{H}^{4,3,3})$ be a solution to  the initial value problem (\ref{3.1})-(\ref{3.2}) for some positive constant $t_1$. Then there exists a constant $\epsilon_0>0$  such that if $\epsilon\leq\epsilon_0$ and $\sup_{0\leq t\leq t_1}\|(\sigma, w, \vartheta)(t)\|_{4,3,3}$,
$\|(\sigma, w, \vartheta)\|_{\mathcal{F}^{5,5,5}}\leq\epsilon$, then it holds
\begin{equation}\label{3.3}
\left\|(\sigma, w, \vartheta)(t)\right\|^2_{4,3,3}+\int_0^t\left\|(\nabla\sigma, \nabla w, \nabla\vartheta)(s)\right\|^2_{4,3,3}\,ds\leq C\left\|(\sigma, w, \vartheta)(0)\right\|^2_{4,3,3}
\end{equation}
for any $t\in[0,t_1]$, where the constant $C>0$ is depending only on $\bar{\rho}, \bar{\theta}, \mu, \mu^\prime, \kappa$ and $\tilde{\alpha}$.
\end{Proposition}
For the proof of the local existence, we can apply the H. Hattori-D. Li \cite{H. Hattori-D. Li-1994} method directly. So we shall devote the following sections to the proof of Proposition 3.2.

Before proving the a priori estimate (\ref{3.3}), let us introduce the absolute constant $\bar{\epsilon}>0$ such that $C_0\bar{\epsilon}=1/4\min\{\bar{\rho}, \bar{\theta}\}$, where $C_0$ is the constant which appears in the inequality $\|\cdot\|_{L^\infty}\leq C_0\|\cdot\|_{2}$. In the following lemmas and their proofs, the small constant $\epsilon$ is at least taken in such a way that
\[\|(\sigma, w, \vartheta)(t)\|_{4,3,3},\quad\|(\sigma, w, \vartheta)\|_{\mathcal{F}^{5,5,5}}\leq\epsilon\leq\bar{\epsilon}\]
so that \[(\rho, \theta)\in\mathcal{S}(\bar{\rho}, \bar{\theta})\equiv\left\{(\rho, \theta)\Big|\frac{\bar{\rho}}
{2}\leq\rho\leq\frac{3\bar{\rho}}{2}, \frac{\bar{\theta}}
{2}\leq\rho\leq\frac{3\bar{\theta}}{2}\right\}.\]

\subsection{Some estimates for $f(t), h(t)$ and their derivatives}
\begin{Lemma}
Let $(\sigma, w, \vartheta)(t)$ and $(\sigma^*, w^*, \vartheta^*)$ be satisfying
\[\|(\sigma, w, \vartheta)(t)\|_{4,3,3},\quad\|(\sigma, w, \vartheta)\|_{\mathcal{F}^{5,5,5}}\leq\epsilon.\]
Then for a multi-index $\alpha$ with $0\leq|\alpha|\leq3$, we have

(i) If we write $\partial_x^\alpha f(t)$ of the form
\[\partial_x^\alpha f(t)=-\frac{\sigma(t)}{\rho^*(\rho^*+\sigma(t))}\left[
\mu\Delta\partial_x^\alpha w(t)+(\mu+\mu^\prime)\nabla
\left(\nabla\cdot\partial_x^\alpha w(t)\right)\right]+F_\alpha(t),\]
then $F_\alpha(t)$ satisfies the estimate:
\begin{eqnarray}\label{3.4}
F_\alpha(t)\leq C\left\{\begin{array}{ll}
|\nabla v^*||w(t)|+(|v^*|+|w(t)|)|\nabla w(t)|+\left(|\nabla\sigma^*|+|\nabla\theta^*|+|\nabla^2v^*|\right)|\sigma(t)|\\[2mm]
+\left(|\nabla\sigma^*|+|\nabla\theta^*|\right)|\vartheta(t)|+\left(|w(t)|+|\sigma(t)|\right)|G|,\quad if\,\alpha=0,\\[2mm]
\left|\nabla^{|\alpha|+1}v^*\right||w(t)|\displaystyle+\sum_{\nu=1}^{|\alpha|+1}\left|\nabla^\nu w(t)\right|+\displaystyle\sum_{\nu=1}^{|\alpha|+1}\left(\left|\nabla^\nu\sigma^*\right|+\left|\nabla^\nu\theta^*\right|+\left|\nabla^{\nu+1}v^*\right|\right)|\sigma(t)|\\[2mm]
+\displaystyle\sum_{\nu=1}^{|\alpha|+1}\left(\left|\nabla^\nu\sigma^*\right|+\left|\nabla^\nu\theta^*\right|\right)|\vartheta(t)|+\sum_{\nu=1}^{|\alpha|}\left(\left|\nabla^\nu \displaystyle\sigma(t)\right|+\left|\nabla^\nu \vartheta(t)\right|\right)\\[2mm]+(|w(t)|+|\sigma(t)|)\displaystyle\sum_{\nu=0}^{|\alpha|}|\nabla^\nu G|+R^{|\alpha|}_F(t), \quad if\,|\alpha|=1,2,3.
 \end{array}\right.
\end{eqnarray}
Here, $R^{k}_F(t)=0, k=1,2$ and $R^{3}_F(t)$ satisfies
\begin{equation}\label{3.5}
\|R^{3}_F(t)\|_{L^{\frac{3}{2}}}\leq C\epsilon\left\|(\nabla^2\sigma, \nabla^3w)(t)\right\|_{1,0}
\end{equation}

(ii) If we write $\partial_x^\alpha h(t)$ of the form
\[\partial_x^\alpha h(t)=\tilde{\alpha}D_1(t)\sigma(t)\Delta\partial_x^\alpha\vartheta(t)+H_\alpha(t),\]
then $H_\alpha(t)$ satisfies the estimate:
\begin{eqnarray}\label{3.6}
H_\alpha(t)\leq C\left\{\begin{array}{ll}
\left(|\nabla^3\theta^*|+|\nabla v^*|+|\nabla v^*||\nabla^2\sigma(t)|+|\nabla w(t)||\nabla^2\sigma^*|\right)|\sigma(t)|+|v^*||\nabla\vartheta(t)|\\[2mm]
+\left(|\nabla\vartheta(t)|+|\nabla\theta^*|\right)|w(t)|
+\left(|\nabla\sigma^*|+|\nabla\sigma(t)|\right)|\nabla\sigma(t)|+|\vartheta(t)||\nabla v^*|\\[2mm]
+\left(|\nabla\sigma(t)|+|\nabla^2\sigma(t)|+|\nabla\sigma^*|+|\nabla^2\sigma^*|+|\nabla w(t)|+|\nabla v^*|\right)|\nabla w(t)|\\[2mm]+|\sigma(t)||H|+|(\sigma, w, \vartheta)(t)||G| ,\quad if\,\alpha=0,\\[2mm]
\displaystyle\sum_{\nu=1}^{|\alpha|+1}\left(\left|\nabla^{\nu+1}\sigma^*\right|+\left|\nabla^{\nu+1}\theta^*\right|+\left|\nabla^{\nu}v^*\right|+\left|\nabla^{\nu-1}H\right|\right)|\sigma(t)|
+\displaystyle\sum_{\nu=1}^{|\alpha|+2}\left|\nabla^\nu\sigma(t)\right|\\[2mm]
+\displaystyle\sum_{\nu=1}^{|\alpha|+1}\left(\left|\nabla^\nu \displaystyle w(t)\right|+\left|\nabla^\nu \vartheta(t)\right|\right)+|\vartheta(t)|\displaystyle\sum_{\nu=1}^{|\alpha|+1}|\nabla^\nu v^*|
+|w(t)|\big|\nabla^{|\alpha|+1}\theta^*\big|\\[2mm]
+|(\sigma, w, \vartheta)(t)|\displaystyle\sum_{\nu=0}^{|\alpha|}\left|\nabla^\nu G\right|+R^{|\alpha|}_H(t), \quad if\,|\alpha|=1,2,3.
 \end{array}\right.
\end{eqnarray}
Here, $R^{1}_H(t)=0$ and $R^{2}_H(t)$, $R^{3}_H(t)$ satisfies
\begin{equation}\label{3.7}
\|R^{2}_H(t)\|_{L^{\frac{3}{2}}}\leq C\epsilon\left\|\nabla^2w(t)\right\|,\quad\|R^{3}_H(t)\|_{L^{\frac{3}{2}}}\leq C\epsilon\left\|(\nabla\sigma, \nabla w)(t)\right\|_{3,2}
\end{equation}
\end{Lemma}
\noindent{\bf Proof.} By the Leibniz formula and the Sobolev embedding: $H^2\hookrightarrow L^\infty$, we can check (\ref{3.4}), (\ref{3.5}) with
\begin{eqnarray}\label{3.8}
\left\{\begin{array}{ll}
R^k_F(t)=0, \quad if\,k=1,2,\\[2mm]
R^3_F(t)=\left|\nabla^2 w(t)\right|\left|\nabla^3 \sigma(t)\right|+\left|\nabla^2 w(t)\right|^2+\left|\nabla^4\sigma^*\right||\nabla\sigma(t)|,
 \end{array}\right.
\end{eqnarray}
and
\begin{eqnarray}\label{3.9}
\left\{\begin{array}{ll}
R^1_H(t)=0, \quad R^2_H(t)=\left|\nabla^2w(t)\right|^2,\\[2mm]
R^3_H(t)=\left(\left|\nabla^2\vartheta(t)\right|+\left|\nabla^2 w(t)\right|\right)\left|\nabla^2 w(t)\right|+\left(\left|\nabla^3 \sigma(t)\right|+
\left|\nabla^4 \sigma(t)\right|+\left|\nabla^3 w(t)\right|\right)\left|\nabla^2w(t)\right|\\[2mm]
+\left|\nabla^3\sigma(t)\right|\left|\nabla^3w(t)\right|+\left(\left|\nabla^4 v^*(t)\right|+
\left|\nabla^5 \sigma^*(t)\right|\right)|\nabla w(t)|+\left(\left|\nabla\sigma(t)\right|+
\left|\nabla^2\sigma(t)\right|\right)|\nabla^4 v^*|
 \end{array}\right.
\end{eqnarray}
For the proof of (\ref{3.5}) and (\ref{3.7}),  we only give here the estimate of the most difficult term $R_H^3$. The others can be dealt with similarly. Using the Gagliard-Nirenberg inequality, we have
\[
\begin{array}{ll}
\|R^{2}_H(t)\|_{L^{\frac{3}{2}}}&\leq C\Big\{\left\|\nabla^2\vartheta(t)\right\|_{L^6}\left\|\nabla^2w(t)\right\|+ \left\|\nabla^2w(t)\right\|_{L^6}\left\|\nabla^2w(t)\right\|+\left\|\nabla^3\sigma(t)\right\|_{L^6}\left\|\nabla^2w(t)\right\|
\\[2mm]&\qquad+\left\|\nabla^2w(t)\right\|_{L^6}\left\|\nabla^3w(t)\right\|+\left\|\nabla^2w(t)\right\|_{L^6}\left\|\nabla^4\sigma(t)\right\|+
\left\|\nabla^3\sigma(t)\right\|_{L^6}\left\|\nabla^3w(t)\right\|\\[2mm]
&\qquad+\left\|(\nabla^4v^*, \nabla^5\sigma^*)\right\|_{L^6}\left\|\nabla w(t)\right\|+\left\|(\nabla\sigma, \nabla^2\sigma)(t)\right\|\left\|\nabla ^4v^*\right\|_{L^6}\Big\}\\[2mm]
&\leq C\epsilon\left\|(\nabla\sigma, \nabla w)(t)\right\|_{3,2}
\end{array}
\]
which is the desired estimate $(\ref{3.7})_2$. This completes the proof of Lemma 3.1.

\subsection{Estimates for $\nabla w(t), \nabla\vartheta(t)$ and their derivatives up to $\nabla^4w(t), \nabla^4\vartheta(t)$}
\begin{Lemma}
Let $(\sigma, w, \vartheta)(t)\in\mathcal{C}(0, t_1; \mathcal{H}^{4,3,3})$ be a solution to the Cauchy problem (\ref{3.1})-(\ref{3.2}) for some positive constant $t_1$. Then there exists four constants $\epsilon_0, \lambda_0>0$ and $d_1, d_2>0$ such that if $0<\epsilon\leq\epsilon_0$ and
$\|(\sigma, w, \vartheta)(t)\|_{4,3,3}, \|(\sigma, w, \vartheta)\|_{\mathcal{F}^{5,5,5}} \leq\epsilon$, then it holds
\begin{equation}\label{3.10}
\begin{array}{ll}
\displaystyle\frac{d}{dt}\left(\left\|\sigma(t)\right\|^2+\langle\hat{A}(t)\nabla\sigma(t), \nabla\sigma(t)\rangle
+\langle \tilde{A}(t)w(t), w(t)\rangle+\langle \tilde{B}(t)\vartheta(t), \vartheta(t)\rangle\right)\\[2mm]
\qquad+d_1\|\nabla w(t)\|^2+d_2\|\nabla \vartheta(t)\|^2\leq C\epsilon\|\nabla \sigma(t)\|_1^2
\end{array}
\end{equation}
and for $1\leq k\leq3$ and any $\lambda$ with $0<\lambda<\lambda_0$,
\begin{equation}\label{3.11}
\begin{array}{ll}
&\displaystyle\frac{d}{dt}\left(\left\|\nabla^k\sigma(t)\right\|^2+\langle\hat{A}(t)\nabla^{k+1}\sigma(t), \nabla^{k+1}\sigma(t)\rangle
+\langle \tilde{A}(t)\nabla^{k}w(t), \nabla^{k}w(t)\rangle\right.\\[2mm]
&\left.\quad+\langle \tilde{B}(t)\nabla^{k}\vartheta(t), \nabla^{k}\vartheta(t)\rangle\right)+d_1\|\nabla^{k+1} w(t)\|^2+d_2\|\nabla^{k+1} \vartheta(t)\|^2\\[2mm]
&\leq C(\epsilon+\lambda)\|\nabla(\sigma, w, \vartheta)(t)\|^2_{k+1, k-1, k-1}+C\lambda^{-1}\|\nabla^k( w, \vartheta)(t)\|^2
\end{array}
\end{equation}
where the constant $C>0$ is depending only on $\bar{\rho}, \bar{\theta}, \mu, \mu^\prime, \kappa$ and $\tilde{\alpha}$. Setting
\[\hat{A}(t)=\frac{\rho}{P_\rho(\rho, \theta)},\quad \tilde{A}(t)=\frac{\rho^2}{P_\rho(\rho, \theta)},\quad\tilde{B}(t)=\frac{C_\triangledown\rho^2}{\theta P_\rho(\rho, \theta)}\]
then $\hat{A}(t)=\hat{A}(\rho^*+\sigma(t), \vartheta^*+\vartheta(t)), \tilde{A}(t)=\tilde{A}(\rho^*+\sigma(t), \vartheta^*+\vartheta(t))$ and
$\tilde{B}(t)=\tilde{B}(\rho^*+\sigma(t), \vartheta^*+\vartheta(t))$.
\end{Lemma}
\noindent{\bf Proof.} Using the Friedrichs mollifier, we may assume that $(\sigma, w, \vartheta)\in\mathcal{C}(0, t_1; \mathcal{H}^{\infty, \infty, \infty})$. For any multi-index $\alpha$ with $0\leq|\alpha|\leq3$, applying $\partial_x^\alpha$ to $(\ref{3.1})_1$, $(\ref{3.1})_2$, $(\ref{3.1})_3$, then taking the $L^2$ inner product of the resultant equations with $\partial_x^\alpha\sigma(t)$, $\tilde{A}(t)\partial_x^\alpha w(t)$, and $\tilde{B}(t)\partial_x^\alpha \vartheta(t)$, respectively, we have
\[
\displaystyle\frac{1}{2}\frac{d}{dt}\left\|\partial_x^\alpha\sigma(t)\right\|^2-\langle(\rho^*+\sigma(t))\partial_x^\alpha w(t), \nabla\partial_x^\alpha \sigma(t)\rangle=\langle\partial_x^\alpha(v^*\sigma(t))+I_\alpha(t), \nabla\partial_x^\alpha \sigma(t)\rangle,
\]
\[\begin{array}{ll}
&\displaystyle\langle\tilde{A}(t)\partial_x^\alpha w_t(t), \partial_x^\alpha w(t)\rangle-\langle\frac{\tilde{A}(t)}{\rho^*}\partial_x^\alpha\{\mu\Delta w(t)+(\mu+\mu^\prime)\nabla
\left(\nabla\cdot w(t)\right)\}, \partial_x^\alpha w(t)\rangle\\[2mm]
&+\langle(\rho^*+\sigma(t))\nabla\partial_x^\alpha \sigma(t), \partial_x^\alpha w(t)\rangle-\kappa\langle\nabla\Delta\partial_x^\alpha \sigma(t), \tilde{A}(t)\partial_x^\alpha w(t)\rangle+\langle \tilde{A}(t)B(t)\nabla\partial_x^\alpha\vartheta(t), \partial_x^\alpha w(t)\rangle\\[2mm]
&=\langle\partial_x^\alpha f(t)+J_\alpha(t),\tilde{A}(t)\partial_x^\alpha w(t)\rangle,
\end{array}
\]
and
\[\begin{array}{ll}
&\displaystyle\langle\tilde{B}(t)\partial_x^\alpha \vartheta_t(t), \partial_x^\alpha \vartheta(t)\rangle-\tilde{\alpha}\langle D^*\Delta\partial_x^\alpha\vartheta(t), \tilde{B}(t)\partial_x^\alpha \vartheta(t)\rangle+\langle \tilde{A}(t)B(t)\nabla\partial_x^\alpha\vartheta(t), \partial_x^\alpha(\nabla\cdot w)(t)\rangle\\[2mm]
&=\langle\partial_x^\alpha h(t)+K_\alpha(t),\tilde{B}(t)\partial_x^\alpha \vartheta(t)\rangle,
\end{array}
\]
where \[
\begin{array}{ll}
&\displaystyle I_\alpha(t)=\sum_{\beta<\alpha}C_\alpha^\beta\partial_x^{\alpha-\beta}\left(\rho^*+\sigma(t)\right)\partial_x^\beta w(t),\\[2mm]
&\displaystyle J_\alpha(t)=\sum_{\beta<\alpha}C_\alpha^\beta\left\{\left(\partial_x^{\alpha-\beta}\frac{1}{\rho^*}\right)\partial_x^\beta (\mu\Delta w(t)+(\mu+\mu^\prime)\nabla
\left(\nabla\cdot w(t)\right))-\left(\partial_x^{\alpha-\beta}A(t)\right)\nabla\partial_x^\beta\sigma(t)\right.\\[2mm]
&\left.\hspace{30mm}-\left(\partial_x^{\alpha-\beta}B(t)\right)\nabla\partial_x^\beta\vartheta(t) \right\},\\[2mm]
&\displaystyle K_\alpha(t)=\sum_{\beta<\alpha}C_\alpha^\beta\left\{\tilde{\alpha}\left(\partial_x^{\alpha-\beta}D^*\right)\Delta\partial_x^\beta\vartheta(t)-\left(\partial_x^{\alpha-\beta}
E(t)\right)\nabla\cdot\partial_x^{\beta}w(t)\right\}.
\end{array}
\]
Canceling the  terms $\langle(\rho^*+\sigma(t))\nabla\partial_x^\alpha \sigma(t), \partial_x^\alpha w(t)\rangle$ and $\langle \tilde{A}(t)B(t)\nabla\partial_x^\alpha\vartheta(t), \partial_x^\alpha w(t)\rangle$ by adding the above three formulas and using the identities
\[
\begin{array}{ll}
&\displaystyle\langle\tilde{A}(t)\partial_x^\alpha w_t(t), \partial_x^\alpha w(t)\rangle=\frac{1}{2}\frac{d}{dt}\langle\tilde{A}(t)\partial_x^\alpha w(t), \partial_x^\alpha w(t)\rangle-\frac{1}{2}\langle\tilde{A}_t(t)\partial_x^\alpha w(t), \partial_x^\alpha w(t)\rangle\\[2mm]
&\displaystyle\langle\tilde{B}(t)\partial_x^\alpha \vartheta_t(t), \partial_x^\alpha \vartheta(t)\rangle=\frac{1}{2}\frac{d}{dt}\langle\tilde{B}(t)\partial_x^\alpha \vartheta(t), \partial_x^\alpha, \vartheta(t)\rangle-\frac{1}{2}\langle\tilde{B}_t(t)\partial_x^\alpha \vartheta(t), \partial_x^\alpha \vartheta(t)\rangle,
\end{array}
\]
we get from integration by parts that
\begin{equation}\label{3.12}
\begin{array}{ll}
&\displaystyle\frac{d}{dt}\left(\left\|\partial_x^\alpha\sigma(t)\right\|^2+\langle \tilde{A}(t)\partial_x^\alpha w(t), \partial_x^\alpha w(t)\rangle+\langle \tilde{B}(t)\partial_x^\alpha\vartheta(t), \partial_x^\alpha\vartheta(t)\rangle\right)\\[3mm]
&\displaystyle+\langle\frac{\tilde{A}(t)}{\rho^*}\nabla\partial_x^\alpha w(t), \partial_x^\alpha\nabla w(t)\rangle+\tilde{\alpha}\langle D^*\nabla\partial_x^\alpha\vartheta(t), \tilde{B}(t)\nabla\partial_x^\alpha \vartheta(t)\rangle\\[2mm]
&\leq \left|\langle\partial_x^\alpha(v^*\sigma(t)), \nabla\partial_x^\alpha\sigma(t)\rangle\right|+\tilde{\alpha}\left|\langle\nabla(D^*\tilde{B}(t))\cdot\nabla\partial_x^\alpha\vartheta(t), \partial_x^\alpha\vartheta(t)\rangle\right|\\[2mm]&
\quad+\left|\langle\nabla(\tilde{A}(t)B(t))\partial_x^\alpha\vartheta(t), \partial_x^\alpha w(t)\rangle\right|+\kappa\langle\nabla\Delta\partial_x^\alpha\sigma(t), \tilde{A}(t)\partial_x^\alpha w(t)\rangle\\[2mm]
&\quad+\left[\mu\left|\langle\nabla\left(\frac{\tilde{A}(t)}{\rho^*}\right)\partial_x^\alpha w(t), \nabla\partial_x^\alpha w(t)\rangle\right|+(\mu+\mu^\prime)\left|\langle\nabla\left(\frac{\tilde{A}(t)}{\rho^*}\right)\cdot\partial_x^\alpha w(t), \partial_x^\alpha(\nabla\cdot w)(t)\rangle\right|\right]\\[2mm]
&\quad+\left|\langle\partial_x^\alpha f(t), \tilde{A}(t)\partial_x^\alpha w(t)\rangle\right|+\left|\langle\partial_x^\alpha h(t), \tilde{B}(t)\partial_x^\alpha \vartheta(t)\rangle\right|\\[2mm]
&\quad+\left[\left|\langle I_\alpha(t), \nabla\partial_x^\alpha \sigma(t)\rangle\right|+\left|\langle J_\alpha(t), \tilde{A}(t)\partial_x^\alpha w(t)\rangle\right|+\left|\langle K_\alpha(t), \tilde{B}(t)\partial_x^\alpha \vartheta(t)\rangle\right|\right]\\[2mm]&
\quad\displaystyle+\frac{1}{2}\left|\langle\tilde{A}_t(t)\partial_x^\alpha w(t), \partial_x^\alpha w(t)\rangle\right|\displaystyle+\frac{1}{2}\left|\langle\tilde{B}_t(t)\partial_x^\alpha \vartheta(t), \partial_x^\alpha \vartheta(t)\rangle\right|\\[2mm]
&\displaystyle=I_1+I_2+\cdots +I_{10}.
\end{array}
\end{equation}
Now, we estimate $I_i, i=1,2, \ldots, 10$ term by term. First, if $\alpha=0$, employing the Hardy inequality, we have
\begin{equation}\label{3.13-1}
I_1\leq\|(1+|x|)v^*\|_{L^\infty}\left\|\frac{\sigma(t)}{|x|}\right\|\left\|\nabla\sigma(t)\right\|\leq C\epsilon\left\|\nabla\sigma(t)\right\|^2.
\end{equation}
If $1\leq|\alpha|\leq3$, using integration by parts and the Sobolev inequality, we get
\begin{equation}\label{3.13}
\begin{array}{ll}
I_1&\leq\displaystyle\sum_{\beta\leq\alpha}C_\alpha^\beta\left|\langle\partial_x^{\alpha-\beta}\nabla\cdot v^*\partial_x^\beta\sigma(t)
+\partial_x^{\alpha-\beta}v^*\cdot\partial_x^\beta\nabla\sigma(t), \partial_x^\alpha\sigma(t)\rangle\right|\\[2mm]
&\leq C\displaystyle\sum_{\beta\leq\alpha}\left\{\left\|\partial_x^{\alpha-\beta}\nabla\cdot v^*\right\|_{L^3}
\left\|\partial_x^{\beta}\sigma(t)\right\|_{L^6}\left\|\partial_x^{\alpha}\sigma(t)\right\|
+\left\|\partial_x^{\alpha-\beta} v^*\right\|_{L^6}
\left\|\partial_x^{\beta}\nabla\sigma(t)\right\|\left\|\partial_x^{\alpha}\sigma(t)\right\|_{L^3}\right\}\\[2mm]
&\leq C\epsilon\left\|\nabla\sigma(t)\right\|_{|\alpha|}^2
\end{array}
\end{equation}
where we have used the following inequalities (cf.\cite{R. Adams-1985, M. E. Talor-1996}).
\[\|u\|_{L^3}\leq\|u\|_{L^2}+\|u\|_{L^6}, \quad\|u\|_{L^6}\leq C\|\nabla u\|, \quad \forall u\in H^1(\mathbb{R}^3)\]
$I_2$ and $I_3$ can be estimated as follows
\begin{equation}\label{3.14}
\begin{array}{ll}
I_2&\leq\displaystyle C\left\|(\nabla\sigma^*, \nabla\sigma(t), \nabla\theta^*, \nabla\vartheta(t))\right\|_{L^3}\left\|\nabla\partial_x^{\alpha}\vartheta(t)\right\|\left\|\partial_x^{\alpha}\sigma(t)\right\|_{L^6}\\[2mm]
&\leq C\epsilon\left\|\nabla\partial_x^{\alpha}\vartheta(t)\right\|^2
\end{array}
\end{equation}
\begin{equation}\label{3.15}
\begin{array}{ll}
I_3&\leq\displaystyle C\left\|(1+|x|)^2(\nabla\sigma^*,  \nabla\theta^*)\right\|_{L^\infty}\left\|\frac{\partial_x^{\alpha}w(t)}{|x|}\right\|\left\|\frac{\partial_x^{\alpha}\vartheta(t)}{|x|}\right\|\\[3mm]
&\quad+C\left\|(\nabla\sigma, \nabla\vartheta)(t)\right\|\left\|\partial_x^{\alpha}w(t)\right\|_{L^3}\left\|\partial_x^{\alpha}\vartheta(t)\right\|_{L^6}\\[2mm]
&\leq C\epsilon\left\|(\nabla \sigma, \nabla w, \nabla\vartheta)(t)\right\|^2_{|\alpha|,|\alpha|,|\alpha|}
\end{array}
\end{equation}
Similar to the estimate of $I_2$, we can get
\begin{equation}\label{3.16}
I_5\leq\displaystyle C\epsilon\left\|\nabla\partial_x^{\alpha} w(t)\right\|^2.
\end{equation}
Now, we turn to estimate $I_4$. We deduce from integration by parts and the equation $(\ref{3.1})_1$
that
\begin{equation}\label{3.17}
\begin{array}{ll}
I_4&=-\displaystyle \kappa\langle\Delta\partial_x^{\alpha}\sigma(t), \nabla\cdot\left(\tilde{A}(t)\partial_x^\alpha w(t)\right)\rangle\\[2mm]
&=-\displaystyle\kappa\langle\Delta\partial_x^{\alpha}\sigma(t), (\nabla\hat{A}(t))(\rho^*+\sigma(t))\partial_x^\alpha w(t)+\hat{A}(t)\nabla\cdot\left\{(\rho^*+\sigma(t))\partial_x^\alpha w(t)\right\}\rangle\\[2mm]
&=-\displaystyle \kappa\langle\Delta\partial_x^{\alpha}\sigma(t), (\nabla\hat{A}(t))(\rho^*+\sigma(t))\partial_x^\alpha w(t)\rangle
+\kappa\langle\Delta\partial_x^{\alpha}\sigma(t), \hat{A}(t)\partial_x^\alpha \sigma_t(t)\rangle\\[2mm]
&\quad\displaystyle+\kappa\sum_{\beta<\alpha}C_\alpha^\beta\langle\Delta\partial_x^{\alpha}\sigma(t), \nabla\cdot\left\{\partial_x^{\alpha-\beta}(\rho^*+\sigma(t)) \partial_x^\beta w(t)\right\}\rangle+\kappa\langle\Delta\partial_x^{\alpha}\sigma(t), \nabla\cdot\partial_x^\alpha(\sigma(t)v^*)\rangle\\[2mm]
&
=I_{4,1}+I_{4,2}+I_{4,3}+I_{4,4}
\end{array}
\end{equation}
$I_{4,1}$ can be estimated as follows
\begin{equation}\label{3.18}
\begin{array}{ll}
I_{4, 1}&\leq C\|\Delta\partial_x^{\alpha}\sigma(t)\|\|\nabla \hat{A}(t)\|_{L^3}\|(\rho^*, \sigma(t))\|_{L^\infty}\|\partial_x^{\alpha}w(t)\|_{L^6}\\[2mm]
&\leq C\|(\nabla\rho^*, \nabla\sigma(t), \nabla\theta^*, \nabla\vartheta(t))\|_1\|(\rho^*, \sigma(t))\|_{L^\infty}\|\nabla\partial_x^\alpha w(t)\|\|\Delta\partial_x^{\alpha}\sigma(t)\|\\[2mm]
&\leq C\epsilon\left(\left\|\nabla\partial_x^{\alpha} w(t)\right\|^2+\left\|\nabla^2\partial_x^{\alpha}\sigma(t)\right\|^2\right)
\end{array}
\end{equation}
For $I_{4,2}$, using integration by parts and the equation $(\ref{3.1})_1$ again, we have
\begin{equation}\label{3.19}
\begin{array}{ll}
I_{4,2}&=\displaystyle-\frac{\kappa}{2}\frac{d}{dt}\langle\nabla\partial_x^\alpha\sigma(t), \hat{A}(t)\nabla\partial_x^\alpha\sigma(t)\rangle+\frac{\kappa}{2}\langle\nabla\partial_x^\alpha\sigma(t), \hat{A}_t(t)\nabla\partial_x^\alpha\sigma(t)\rangle
\\[2mm]
&\quad+\langle\nabla\partial_x^\alpha\sigma(t), \nabla\hat{A}(t)\nabla\cdot\left\{(\rho^*+\sigma(t))\partial_x^\alpha w(t)\right\}\rangle\\[2mm]
&\quad\displaystyle+\kappa\sum_{\beta<\alpha}C_\alpha^\beta\langle\nabla\partial_x^{\alpha}\sigma(t), \nabla\hat{A}(t)\nabla\cdot\left\{\partial_x^{\alpha-\beta}(\rho^*+\sigma(t))\partial_x^\beta w(t)\right\}\rangle\\[2mm]
&\quad+\langle\nabla\partial_x^\alpha\sigma(t), \nabla\hat{A}(t)\nabla\cdot\partial_x^\alpha (\sigma(t)v^*)\rangle\\[2mm]
&=\displaystyle-\frac{\kappa}{2}\frac{d}{dt}\langle\nabla\partial_x^\alpha\sigma(t), \hat{A}(t)\nabla\partial_x^\alpha\sigma(t)\rangle+I_{4,2}^1+I_{4,2}^2+I_{4,2}^3+I_{4,2}^4.
\end{array}
\end{equation}
To estimate $I_{4,2}^1$, we use the equation $(\ref{3.1})_1$ and $(\ref{3.1})_3$,
\begin{equation}\label{3.20}
\begin{array}{ll}
I_{4,2}^1&=\displaystyle\frac{\kappa}{2}\langle\left(\hat{A}_\rho(t)\sigma_t(t)+\hat{A}_\theta(t)\vartheta_t(t)\right)\nabla\partial_x^\alpha\sigma(t), \nabla\partial_x^\alpha\sigma(t)\rangle
\\[2mm]
&\leq C\left\|\nabla\cdot\{(\rho^*+\sigma(t))w(t)\}+\nabla\cdot(v^*\sigma(t))\right\|_{L^\infty}\left\|\nabla\partial_x^\alpha\sigma(t)\right\|^2\\[2mm]
&\quad+C\left\|h(t)-\tilde{\alpha}(D(t)-D^*)\Delta\vartheta(t)-E(t)(\nabla\cdot w)(t)\right\|_{L^\infty}\left\|\nabla\partial_x^\alpha\sigma(t)\right\|^2\\[2mm]
&\displaystyle\quad+\frac{\kappa\tilde{\alpha}}{2}\left|\langle\hat{A}_\theta(t) D(t)\Delta\vartheta(t)\nabla\partial_x^\alpha\sigma(t), \nabla\partial_x^\alpha\sigma(t)\rangle\right|\\[2mm]
&\leq C\epsilon\left\|\nabla\partial_x^\alpha\sigma(t)\right\|^2+\displaystyle\frac{\kappa\tilde{\alpha}}{2}\left|\langle\nabla\big(\hat{A}_\theta(t) D(t)\big)\cdot\nabla\vartheta(t), |\nabla\partial_x^\alpha\sigma(t)|^2\rangle\right|\\[2mm]
&\displaystyle\quad+\kappa\tilde{\alpha}\left|\langle\hat{A}_\theta(t) D(t)\nabla\vartheta(t)\cdot\nabla\partial_x^\alpha\sigma(t), \Delta\partial_x^\alpha\sigma(t)\rangle\right|\\[2mm]
&\leq C\epsilon\left\|\nabla\partial_x^\alpha\sigma(t)\right\|^2+C\left\|\nabla\vartheta(t)\right\|_{L^\infty}\left\|\nabla\partial_x^\alpha\sigma(t)\right\|
\left\|\Delta\partial_x^\alpha\sigma(t)\right\|\\[2mm]
&\displaystyle\quad+C\left\|\left(\nabla\sigma^*, \nabla\sigma(t), \nabla\vartheta^*, \nabla\vartheta(t)\right)\right\|_{L^\infty}\left\|\nabla\partial_x^\alpha\sigma(t)\right\|^2\\[2mm]
&\leq C\epsilon\left\|\nabla\partial_x^\alpha\sigma(t)\right\|_1^2.
\end{array}
\end{equation}
It follows from the H\"{o}lder inequality and the Sobolev inequality that
\begin{equation}\label{3.21}
\begin{array}{ll}
I_{4,2}^2&\leq\displaystyle C\left\|\left(\nabla\sigma^*, \nabla\sigma(t), \nabla\vartheta(t), \nabla\vartheta^*\right)\right\|_{L^6}\left\|(\nabla\sigma^*, \nabla\sigma(t))\right\|_{L^6}\left\|\partial_x^\alpha w(t)\right\|_{L^6}\left\|\nabla\partial_x^\alpha \sigma(t)\right\|\\[2mm]
&\quad+C\left\|(\rho^*, \sigma(t))\right\|_{L^\infty}\left\|\left(\nabla\sigma^*, \nabla\sigma(t), \nabla\vartheta(t), \nabla\vartheta^*\right)\right\|_{L^\infty}\left\|\partial_x^\alpha(\nabla\cdot w)(t)\right\|\left\|\nabla\partial_x^\alpha \sigma(t)\right\|\\[2mm]
&\leq C\epsilon\left(\left\|\nabla\partial_x^\alpha \sigma(t)\right\|^2+\left\|\nabla\partial_x^\alpha w(t)\right\|^2\right)
\end{array}
\end{equation}
and
\begin{equation}\label{3.22}
\begin{array}{ll}
I_{4,2}^3+I_{4,2}^4&\leq\displaystyle C\sum_{\beta<\alpha}C_\alpha^\beta\left\|\nabla\partial_x^\alpha \sigma(t)\right\|\left\|\left(\nabla\sigma^*, \nabla\sigma(t), \nabla\vartheta(t), \nabla\vartheta^*\right)\right\|_{L^6}\\[2mm]
&\quad\times\left\{\left\|\partial_x^{\alpha-\beta}\nabla(\rho^*+\sigma(t))\right\|_{L^6}\left\|\partial_x^{\beta}w(t)\right\|_{L^6}
+\left\|\partial_x^{\alpha-\beta}(\rho^*+\sigma(t))\right\|_{L^6}\left\|\nabla\partial_x^{\beta}w(t)\right\|_{L^6}\right.\\[2mm]
&\qquad\left.+\left\|\partial_x^{\alpha-\beta}\nabla\sigma(t)\right\|_{L^6}\left\|\partial_x^{\beta}v^*\right\|_{L^6}
+\left\|\partial_x^{\alpha-\beta}\sigma(t)\right\|_{L^6}\left\|\partial_x^{\beta}(\nabla\cdot v^*)\right\|_{L^6}\right\}\\[2mm]
&\leq C\epsilon\left(\left\|\nabla\partial_x^{\alpha}\sigma(t)\right\|^2_1+\left\|\nabla w(t)\right\|^2_{|\alpha|}\right).
\end{array}
\end{equation}
Combining (\ref{3.19})-(\ref{3.22}), we obtain
\begin{equation}\label{3.23}
I_{4,2}\leq\displaystyle-\frac{\kappa}{2}\frac{d}{dt}\langle\nabla\partial_x^\alpha\sigma(t), \hat{A}(t)\nabla\partial_x^\alpha\sigma(t)\rangle+C\epsilon\left(\left\|\nabla\partial_x^{\alpha}\sigma(t)\right\|^2_1+\left\|\nabla w(t)\right\|^2_{|\alpha|}\right).
\end{equation}
Similar to the estimate of $I_1$, we can also get
\begin{equation}\label{3.24}
I_{4,3}\leq C\epsilon\left(\left\|\nabla\sigma(t)\right\|^2_{|\alpha|+1}+\left\|\nabla w(t)\right\|^2_{|\alpha|}\right).
\end{equation}
\begin{equation}\label{3.25}
I_{4,4}\leq C\epsilon\left\|\nabla\sigma(t)\right\|^2_{|\alpha|+1}.
\end{equation}
Putting (\ref{3.18}), (\ref{3.23})-(\ref{3.25}) into (\ref{3.17}) gives rise to
\begin{equation}\label{3.26}
I_{4}\leq\displaystyle-\frac{\kappa}{2}\frac{d}{dt}\langle \hat{A}(t)\nabla\partial_x^\alpha\sigma(t),\nabla\partial_x^\alpha\sigma(t)\rangle+ C\epsilon\left(\left\|\nabla\sigma(t)\right\|^2_{|\alpha|+1}+\left\|\nabla w(t)\right\|^2_{|\alpha|}\right).
\end{equation}
To estimate $I_6$ and $I_7$, we use Lemma 3.1. Here, we only give the detailed estimation of $I_7$. $I_6$ can be estimated similarly. In fact, $I_7$ can be divided into the following two parts
\begin{equation}\label{3.27}
\begin{array}{ll}
I_{7}&\leq\displaystyle\left|\langle H_\alpha(t), \tilde{B}(t)\partial_x^\alpha\vartheta(t)\rangle\right|+\tilde{\alpha}\left|\langle D_1(t)\sigma(t)\Delta\partial_x^\alpha\vartheta(t), \tilde{B}(t)\partial_x^\alpha\vartheta(t)\rangle\right|\\[2mm]
&=I_{7,1}+I_{7,2}.
\end{array}
\end{equation}
For $I_{7,2}$, using integration by parts, we have
\begin{equation}\label{3.28}
\begin{array}{ll}
I_{7,2}&\leq\displaystyle\tilde{\alpha}\left|\langle\nabla(D_1(t)\sigma(t))\cdot\nabla\partial_x^\alpha\vartheta(t), \tilde{B}(t)\partial_x^\alpha\vartheta(t)\rangle\right|+\tilde{\alpha}\left|\langle D_1(t)\sigma(t)\nabla\partial_x^\alpha\vartheta(t), \nabla\left(\tilde{B}(t)\partial_x^\alpha\vartheta(t)\right)\rangle\right|\\[2mm]
&\leq\left\|(\nabla\sigma^*, \nabla\sigma(t), \nabla\vartheta^*, \nabla\vartheta(t))\right\|_{L^3}\left\|\nabla\partial_x^\alpha\vartheta(t)\right\|
\left\|\partial_x^\alpha\vartheta(t)\right\|_{L^6}+C\left\|\sigma(t)\right\|_{L^\infty}\left\|\nabla\partial_x^\alpha\vartheta(t)\right\|^2\\[2mm]
&\leq C\epsilon\left\|\nabla\partial_x^\alpha\vartheta(t)\right\|^2.
\end{array}
\end{equation}
To estimate $I_{7,1}$, we use (\ref{3.6}). If $\alpha=0$,
\begin{equation}\label{3.29}
\begin{array}{ll}
I_{7,1}&\leq\displaystyle C\bigg\{\left\|(1+|x|)^2(\nabla v^*, \nabla^2\sigma^*, \nabla\theta^*, \nabla^3\theta^*, H, G)\right\|_{L^\infty}\left\|\frac{(\sigma, w, \vartheta)(t)}{|x|}\right\|
\left\|\frac{\vartheta(t)}{|x|}\right\|\\[2mm]
&\quad+\left\|(1+|x|)v^*\right\|_{L^\infty}\left\|\frac{\vartheta(t)}{|x|}\right\|\left\|\nabla\vartheta(t)\right\|+\left\|w(t)\right\|_{L^3}\left\|\vartheta(t)\right\|_{L^6}
\|\left\|\nabla\vartheta(t)\right\|\\[2mm]
&\quad+\left\|(\nabla\sigma(t), \nabla^2\sigma(t),\nabla w(t), \nabla\sigma^*, \nabla^2\sigma^*,\nabla v^*)\right\|_{L^3}\left\|\nabla w(t)\right\|\left\| \vartheta(t)\right\|_{L^6}\\[2mm]
&\quad+\left\|(\nabla\sigma^*, \nabla\sigma(t))\right\|\left\|\nabla \sigma(t)\right\|_{L^3}\left\| \vartheta(t)\right\|_{L^6}\bigg\}\\[2mm]
&\leq C\epsilon\left\|(\nabla\sigma, \nabla w, \nabla\vartheta)(t)\right\|^2_{1,0,0}.
\end{array}
\end{equation}
and if $1\leq|\alpha|\leq3$,
\begin{equation}\label{3.30}
\begin{array}{ll}
I_{7,1}&\leq\displaystyle C\left\{\displaystyle\sum_{\nu=1}^{|\alpha|+1}\left\|(\nabla^{\nu+1}\sigma^*, \nabla^{\nu}v^*, \nabla^{\nu+1}\theta^*, \nabla^{\nu-1}H)\right\|_{L^3}\left\|\sigma(t)\right\|_{L^6}
\left\|\partial_x^\alpha\vartheta(t)\right\|\right.
\\[3mm]&\quad+\displaystyle\sum_{\nu=1}^{|\alpha|+2}\left\|\nabla^\nu\sigma(t)\right\|\left\|\partial_x^\alpha\vartheta(t)\right\|
+\displaystyle\sum_{\nu=1}^{|\alpha|+1}\left( \left\|\nabla^\nu w(t)\right\|+\left\|\nabla^\nu\vartheta(t)\right\|\right)\left\|\partial_x^\alpha\vartheta(t)\right\|\\[2mm]&\quad
+\displaystyle\sum_{\nu=1}^{|\alpha|+1}\left\|\nabla^\nu v^*\right\|_{L^3}\left\| \vartheta(t)\right\|_{L^6}\left\|\partial_x^\alpha\vartheta(t)\right\|
+\left\|\nabla^{|\alpha|+1} \theta^*\right\|_{L^3}\left\| w(t)\right\|_{L^6}\left\|\partial_x^\alpha\vartheta(t)\right\|
\\[2mm]&\quad\left.+\displaystyle\sum_{\nu=0}^{|\alpha|}\left\|(\sigma, w, \vartheta)(t)\right\|_{L^6}\left\| \nabla^\nu G\right\|_{L^3}\left\|\partial_x^\alpha\vartheta(t)\right\|+\left\|R_{H}^{|\alpha|}(t)\right\|_{L^{3/2}}\left\|\partial_x^\alpha\vartheta(t)\right\|_{L^3}
\right\}\\[2mm]
&\leq C(\epsilon+\lambda)\left\|(\nabla\sigma, \nabla w, \nabla\vartheta)(t)\right\|^2_{|\alpha|+1, |\alpha|, |\alpha|}+C\lambda^{-1}\left\|\nabla^{|\alpha|}\vartheta(t)\right\|^2.
\end{array}
\end{equation}
Combining (\ref{3.27})-(\ref{3.30}), we have
\begin{equation}\label{3.31}
I_{7}\leq C\left\{\begin{array}{ll} \epsilon\left\|(\nabla\sigma, \nabla w, \nabla\vartheta)(t)\right\|^2_{1,0,0},\quad \alpha=0,\\[2mm]
(\epsilon+\lambda)\left\|(\nabla\sigma, \nabla w, \nabla\vartheta)(t)\right\|^2_{|\alpha|+1, |\alpha|, |\alpha|}+\lambda^{-1}\left\|\nabla^{|\alpha|}\vartheta(t)\right\|^2,\quad1\leq|\alpha|\leq3.
\end{array}\right.
\end{equation}
By using the same argument as $I_7$, one can get
\begin{equation}\label{3.32}
I_{6}\leq C\left\{\begin{array}{ll} \epsilon\left\|(\nabla\sigma, \nabla w, \nabla\vartheta)(t)\right\|^2,\quad \alpha=0,\\[2mm]
(\epsilon+\lambda)\left\|(\nabla\sigma, \nabla w, \nabla\vartheta)(t)\right\|^2_{|\alpha|-1, |\alpha|, |\alpha|-1}+\lambda^{-1}\left\|\nabla^{|\alpha|}w(t)\right\|^2,\quad1\leq|\alpha|\leq3.
\end{array}\right.
\end{equation}
Similar to the estimate of $I_1$, it is easy to check that
\begin{equation}\label{3.33}
I_{8}\leq C\epsilon\left\|(\nabla\sigma, \nabla w, \nabla\vartheta)(t)\right\|^2_{|\alpha|, |\alpha|, |\alpha|}.
\end{equation}
In order to estimate $I_9$ and $I_{10}$, we use the equations $(\ref{3.1})_1$
and $(\ref{3.1})_2$ again. In fact for $I_9$,
\begin{equation}\label{3.34}
\begin{array}{ll}
2I_{9}&\leq C \displaystyle\left|\langle \tilde{A}_\rho(t)\sigma_t(t)\partial_x^\alpha w(t),\partial_x^\alpha w(t)\rangle\right|+\left|\langle \tilde{A}_\theta(t)\vartheta_t(t)\partial_x^\alpha w(t), \partial_x^\alpha w(t)\rangle\right|\\[2mm]
&=I_{9,1}+I_{9,2}
\end{array}
\end{equation}
Using $(\ref{3.1})_1$, $(\ref{3.1})_2$ and (\ref{3.6}), $I_{9,1}$ and $I_{9,2}$
can be estimated as follows
\begin{equation}\label{3.35}
\begin{array}{ll}
I_{9, 1}&= \displaystyle\left|\langle\nabla\cdot\{(\rho^*+\sigma(t))w(t)+v^*\sigma(t)\}, \tilde{A}_\rho(t)\partial_x^\alpha w(t)\partial_x^\alpha w(t)\rangle\right|\\[2mm]
&\leq C\left\|(\nabla\sigma^*,\nabla\sigma(t), \nabla w(t), \nabla v^*)\right\|\left\|(w(t),\sigma(t), \rho^*,  v^*)\right\|_{L^6}
\left\|\partial_x^\alpha w(t)\right\|^2_{L^6}\\[2mm]&
\leq C\epsilon\left\|\nabla\partial_x^\alpha w(t)\right\|^2
\end{array}
\end{equation}
\begin{equation}\label{3.36}
\begin{array}{ll}
I_{9, 2}&= \displaystyle\left|\langle\tilde{\alpha}D(t)\Delta\vartheta(t)-E(t)(\nabla\cdot w)(t)+h(t), \tilde{A}_\theta(t)\partial_x^\alpha w(t)\partial_x^\alpha w(t)\rangle\right|\\[2mm]
&\leq \tilde{\alpha}\left|\langle \nabla D(t)\nabla\vartheta(t), \tilde{A}_\theta(t)\partial_x^\alpha w(t)\partial_x^\alpha w(t)\rangle\right|+\tilde{\alpha}\left|\langle D(t)\nabla\vartheta(t), \nabla\tilde{A}_\theta(t)\partial_x^\alpha w(t)\partial_x^\alpha w(t)\rangle\right|\\[2mm]&
\quad+2\tilde{\alpha}\left|\langle D(t)\nabla\vartheta(t), \tilde{A}_\theta(t)\nabla\partial_x^\alpha w(t)\partial_x^\alpha w(t)\rangle\right|\\[2mm]
&\quad+\left|\langle-E(t)(\nabla\cdot w)(t)+H_0(t), \tilde{A}_\theta(t)\partial_x^\alpha w(t)\partial_x^\alpha w(t)\rangle\right|\\[2mm]
&\leq C\left\{\left\|\partial_x^\alpha w(t)\right\|^2_{L^6}\big(\left\|(\nabla^3\theta^*, \nabla v^*, \nabla^2\sigma(t), \nabla w(t), \nabla\vartheta(t), \nabla\theta^*, G, H)\right\|\right.\\[2mm]
&\quad\times\left\|(\sigma(t), w(t), \vartheta(t), v^*)\right\|_{L^6}+\left\|(\nabla\sigma(t), \nabla^2\sigma(t), \nabla\sigma^*, \nabla^2\sigma^*, \nabla\vartheta(t), \nabla\theta^*, \nabla w(t), \nabla v^*)\right\|\\[2mm]
&\quad\times\left\|(\nabla\sigma, \nabla w, \nabla\vartheta)(t)\right\|_{L^6}\big)+\left\|\nabla\vartheta(t)\right\|_{L^3}\left\|\nabla\partial_x^\alpha w(t)\right\|\left\|\partial_x^\alpha w(t)\right\|_{L^6}\Big\}\\[2mm]
&\leq C\epsilon\left\|\nabla\partial_x^\alpha w(t)\right\|^2
\end{array}
\end{equation}
Consequently,
\begin{equation}
\label{3.37}
I_9\leq C\epsilon\left\|\nabla\partial_x^\alpha w(t)\right\|^2.
\end{equation}
Finally, the term $I_{10}$ is estimated in way similar to that of $I_9$, and we have
\begin{equation}
\label{3.38}
I_{10}\leq C\epsilon\left\|\nabla\partial_x^\alpha \vartheta(t)\right\|^2.
\end{equation}
Combining (\ref{3.12})-(\ref{3.16}), (\ref{3.26}), (\ref{3.31})-(\ref{3.33}), (\ref{3.37}) and (\ref{3.38}), we can obtain (\ref{3.10}) and (\ref{3.11}), if we take $d_1=\min_{(\rho, \theta)\in \mathcal{S}(\bar{\rho}, \bar{\theta})}\{\frac{\mu\tilde{A}(\rho, \theta)}{\rho}\}$,
$d_2=\min_{(\rho, \theta)\in \mathcal{S}(\bar{\rho}, \bar{\theta})}(\tilde{\alpha}D(\rho, \theta)\tilde{B}(\rho, \theta)\}$ and choose $\epsilon$ and $\lambda$ small enough. This completes the proof of Lemma 3.2.

\subsection{Estimates for $\nabla\sigma(t)$ and its derivatives up to $\nabla^5\sigma(t)$}
\begin{Lemma}
Let $(\sigma, w, \vartheta)(t)\in\mathcal{C}(0, t_1; \mathcal{H}^{4,3,3})$ be a solution to the Cauchy problem (\ref{3.1})-(\ref{3.2}) for some positive constant $t_1$. Then there exists three  constants $\epsilon_0, \lambda_0>0$ and $d_3>0$ such that if $\epsilon\leq\epsilon_0$ and
$\|(\sigma, w, \vartheta)(t)\|_{4,3,3}, \|(\sigma, w, \vartheta)\|_{\mathcal{F}^{5,5,5}} \leq\epsilon$, then it holds
\begin{equation}\label{3.39}
\displaystyle\frac{d}{dt}\langle w(t), \nabla\sigma(t)\rangle
+d_3\left\|\nabla\sigma(t)\right\|^2+\kappa\left\|\nabla^2\sigma(t)\right\|^2
\leq C\|(\nabla w, \nabla\vartheta)(t)\|^2
\end{equation}
and for $1\leq k\leq3$ and any $\lambda$ with $0<\lambda<\lambda_0$,
\begin{equation}\label{3.40}
\displaystyle\frac{d}{dt}\langle \nabla^kw(t), \nabla^{k+1}\sigma(t)\rangle
+d_3\left\|\nabla^{k+1}\sigma(t)\right\|^2+\kappa\left\|\nabla^{k+2}\sigma(t)\right\|^2
\leq C\|(\nabla\sigma, \nabla w, \nabla\vartheta)(t)\|^2
\end{equation}
where the constant $C>0$ is depending only on $\bar{\rho}, \bar{\theta}, \mu, \mu^\prime, \kappa$ and $\tilde{\alpha}$.
\end{Lemma}
\noindent{\bf Proof.}  Using the Friedrichs mollifier, we may assume that $(\sigma, w, \vartheta)\in\mathcal{C}(0, t_1; \mathcal{H}^{\infty, \infty, \infty})$. For any multi-index $\alpha$ with $0\leq|\alpha|\leq3$, applying $\partial_x^\alpha$ to $(\ref{3.1})_2$, then taking the $L^2$ inner product of the resultant equations with $\partial_x^\alpha\nabla\sigma(t)$, we have
\begin{equation}\label{3.41}
\begin{array}{ll}
&\langle A(t)\nabla\partial_x^\alpha\sigma(t), \nabla\partial_x^\alpha \sigma(t)\rangle+\kappa\left\|\nabla^2\partial_x^\alpha\sigma(t)\right\|^2\\[2mm]
&=-\langle\partial_x^\alpha w_t(t), \nabla\partial_x^\alpha \sigma(t)\rangle+\left|\langle\partial_x^\alpha\left\{\frac{1}{\rho^*}\left[\mu\Delta w(t)+(\mu+\mu^\prime)\nabla
\left(\nabla\cdot\right) w(t)\right]\right\},\nabla\partial_x^\alpha\sigma(t)\rangle\right|\\[2mm]
&\quad+\displaystyle\sum_{\beta<\alpha}C_\alpha^\beta\left|\langle\partial_x^{\alpha-\beta}A(t)\partial_x^\beta\nabla\sigma(t),\nabla\partial_x^\alpha\sigma(t)\rangle\right|\\[2mm]
&\quad+\left|\langle\partial_x^{\alpha}(B(t)\nabla\vartheta(t)), \nabla\partial_x^\alpha\sigma(t)\rangle\right|+\left|\langle\partial_x^{\alpha}f(t),\nabla\partial_x^\alpha\sigma(t)\rangle\right|\\[2mm]
&=I_1+I_2+\cdots+I_5.
\end{array}
\end{equation}
For $I_1$, we deduce from integration by parts and $(\ref{3.1})_1$ that
\begin{equation}\label{3.42}
\begin{array}{ll}
I_1&=\displaystyle-\frac{d}{dt}\langle\partial_x^\alpha w(t), \nabla\partial_x^\alpha \sigma(t)\rangle-\langle\partial_x^\alpha (\nabla\cdot w)(t), \partial_x^\alpha \sigma_t(t)\rangle\\[2mm]
&=\displaystyle-\frac{d}{dt}\langle\partial_x^\alpha w(t), \nabla\partial_x^\alpha \sigma(t)\rangle+\langle\partial_x^\alpha (\nabla\cdot w)(t), \partial_x^\alpha \nabla\cdot\{(\rho^*+\sigma(t))w(t)\}\rangle\\[2mm]
&\quad+\displaystyle\langle\partial_x^\alpha (\nabla\cdot w)(t), \partial_x^\alpha \nabla\cdot(\sigma(t)v^*)\rangle\\[2mm]
&=\displaystyle-\frac{d}{dt}\langle\partial_x^\alpha w(t), \nabla\partial_x^\alpha \sigma(t)\rangle+I_{1,1}+I_{1,2}
\end{array}
\end{equation}
By using the way similar to that of (\ref{3.13}), we have
\begin{equation}\label{3.43}
I_{1,1}+I_{1,2}\leq C\epsilon\left\|\nabla^2\partial_x^\alpha \sigma(t)\right\|^2+C \epsilon\left(\left\|\nabla w(t)\right\|^2_{|\alpha|+1}+\left\|\nabla \sigma(t)\right\|^2_{|\alpha|}\right)
\end{equation}
For $I_2$, let $\alpha_0\leq\alpha$ be a multi-index with $|\alpha_0|=1$, then it follows from integration by parts and the Cauchy inequality that
\begin{equation}\label{3.44}
\begin{array}{ll}
I_2&=\displaystyle\left|\langle\partial_x^{\alpha-\alpha_0}\left\{\frac{1}{\rho^*}\left[\mu\Delta w(t)+(\mu+\mu^\prime)\nabla
\left(\nabla\cdot\right) w(t)\right]\right\},\nabla\partial_x^{\alpha+\alpha_0}\sigma(t)\rangle\right|\\[2mm]
&\leq\lambda\left\|\nabla^2\partial_x^\alpha\sigma(t)\right\|^2+C \lambda^{-1}\left\|\nabla w(t)\right\|^2_{|\alpha|}
\end{array}
\end{equation}
Using the Cauchy inequality , $I_3$ and $I_4$ can be estimated as follows
\begin{equation}\label{3.44}
I_3\leq C\epsilon\left\|\nabla\sigma(t)\right\|^2_{|\alpha|},
\end{equation}
\begin{equation}\label{3.45}
I_4\leq \lambda\left\|\nabla\partial_x^\alpha\sigma(t)\right\|^2+C \lambda^{-1}\left\|\nabla \vartheta(t)\right\|^2_{|\alpha|}.
\end{equation}
Finally, similar to the estimate of $I_7$ in the proof of Lemma 3.2, we have
\begin{equation}\label{3.47}
I_5\leq C(\epsilon+\lambda)\left\|\nabla\partial_x^\alpha\sigma(t)\right\|^2_1+C \lambda^{-1}\left\|(\nabla\sigma, \nabla w, \nabla\vartheta)(t)\right\|^2_{|\alpha|-1, |\alpha|, |\alpha|-1}.
\end{equation}
Combining (\ref{3.41})-(\ref{3.47}) and summing up $\alpha$ with $|\alpha|=k$, we can get (\ref{3.40}), if we take $d_3=\min_{(\rho, \theta)\in\mathcal{S}(\bar{\rho}, \bar{\theta})}\left\{\frac{P_\rho(\rho, \theta)}{\rho}\right\}$ and choose $\epsilon$ and $\lambda$ small enough.
For $\alpha=0$, by using the same argument as above, we can also get (\ref{3.39}). This completes the proof of Lemma 3.3.

\subsection{Proof of Proposition 3.2}
Let $(\sigma, w, \vartheta)(t)\in\mathcal{C}(0, t_1; \mathcal{H}^{4,3,3})$ be a solution to the Cauchy problem (\ref{3.1})-(\ref{3.2}) for some positive constant $t_1$. Furthermore, we assume that $\|(\sigma, w, \vartheta)(t)\|_{4,3,3}, \|(\sigma, w, \vartheta)\|_{\mathcal{F}^{5,5,5}} \leq\epsilon$, where $\epsilon>0$ is small enough such that we can use the results obtained in Lemmas 3.2-3.3. Set
\[\left[\sigma, w, \vartheta\right](t)=\left\|\sigma(t)\right\|^2+\langle\hat{A}(t)\nabla\sigma(t), \nabla\sigma(t)\rangle+\langle\tilde{A}(t)w(t), w(t)\rangle+\langle \tilde{B}(t)\vartheta(t), \vartheta(t)\rangle
\]
where $\langle\hat{A}(t)$, $\tilde{A}(t)$ and $\tilde{B}(t)$ are defined as in Lemma 3.2.

Multiplying (\ref{3.39}) with a small constant $\lambda_0$, then adding the resultant equation to (\ref{3.10}), we have
\begin{equation}
\label{3.48}
\displaystyle\frac{d}{dt}\left\{a_0\left[\sigma, w, \vartheta\right](t)+b_0\langle w(t), \nabla\sigma(t)\rangle\right\}+\left\|(\nabla\sigma, \nabla w, \nabla\vartheta)(t)\right\|^2_{1,0,0}\leq0
\end{equation}
provided that $\epsilon>0$ is small enough. Here and hereafter, $a_\nu, b_\nu>0, \nu=0,1, \cdots, 3$ denote some  constants depending only on $\bar{\rho},  \bar{\theta},\mu, \mu^\prime, \kappa$ and $\tilde{\alpha}$. Then summing up (\ref{3.40}), (\ref{3.11}) with $k=1$ and (\ref{3.48}), we get
\begin{equation}
\label{3.49}
\displaystyle\frac{d}{dt}\left\{\sum_{\nu=0}^1a_\nu\left[\nabla^\nu\sigma, \nabla^\nu w, \nabla^\nu\vartheta\right](t)+\sum_{\nu=0}^1b_\nu\langle\nabla^\nu w(t), \nabla^{\mu+1}\sigma(t)\rangle\right\}+\left\|(\nabla\sigma, \nabla w, \nabla\vartheta)(t)\right\|^2_{2,1,1}\leq0
\end{equation}
Similarly, summing up (\ref{3.40}), (\ref{3.11}) with $k=2$ and (\ref{3.49}) gives
 \begin{equation}
\label{3.50}
\displaystyle\frac{d}{dt}\left\{\sum_{\nu=0}^2a_\nu\left[\nabla^\nu\sigma, \nabla^\nu w, \nabla^\nu\vartheta\right](t)+\sum_{\nu=0}^2b_\nu\langle\nabla^\nu w(t), \nabla^{\mu+1}\sigma(t)\rangle\right\}+\left\|(\nabla\sigma, \nabla w, \nabla\vartheta)(t)\right\|^2_{3,2,2}\leq0
\end{equation}
and summing up (\ref{3.40}), (\ref{3.11}) with $k=3$ and (\ref{3.50}) gives
\begin{equation}
\label{3.51}
\displaystyle\frac{d}{dt}\left\{\sum_{\nu=0}^3a_\nu\left[\nabla^\nu\sigma, \nabla^\nu w, \nabla^\nu\vartheta\right](t)+\sum_{\nu=0}^3b_\nu\langle\nabla^\nu w(t), \nabla^{\mu+1}\sigma(t)\rangle\right\}+\left\|(\nabla\sigma, \nabla w, \nabla\vartheta)(t)\right\|^2_{4,3,3}\leq0
\end{equation}
Integrating (\ref{3.51}) with respect to $t$ over $[0, t]$, we have
\begin{equation}
\label{3.52}
\displaystyle N\left[\sigma, w, \vartheta\right](t)+\int_0^t\left\|(\nabla\sigma, \nabla w, \nabla\vartheta)(s)\right\|^2_{4,3,3}\,ds\leq N\left[\nabla^\nu\sigma, \nabla^\nu w, \nabla^\nu\vartheta\right](0)
\end{equation}
for any $t\in[0,t_1]$, where \[ N\left[\sigma, w, \vartheta\right](t)\equiv\sum_{\nu=0}^3a_\nu\left[\nabla^\nu\sigma, \nabla^\nu w, \nabla^\nu\vartheta\right](t)+\sum_{\nu=0}^3b_\nu\langle\nabla^\nu w(t), \nabla^{\mu+1}\sigma(t)\rangle,\quad t>0\]

Denote $B_0=\min_{(\rho, \theta)\in\mathcal{S}(\bar{\rho}, \bar{\theta})}\left\{\hat{A}(\rho, \theta), \tilde{A}(\rho, \theta), \tilde{B}(\rho, \theta), 1\right\}$ and $B_1=\max_{(\rho, \theta)\in\mathcal{S}(\bar{\rho}, \bar{\theta})}\{\hat{A}(\rho, \theta), \tilde{A}(\rho, \theta),\\ \tilde{B}(\rho, \theta), 1\}$. Since we may assume without loss of generality that $a_\nu\leq a_{\nu-1}$ and $b_\nu\leq a_\nu\min\{B_0, 1\}/4$ for $\nu=1,2,3$, it follows from simple calculation that
\begin{equation}\label{3.53}
\displaystyle \frac{\alpha_3}{4}B_0\left\|(\sigma, w, \vartheta)(t)\right\|^2_{4,3,3}\leq N\left[\sigma, w, \vartheta\right](t)\leq2\alpha_0\left\|(\sigma, w, \vartheta)(t)\right\|^2_{4,3,3}
\end{equation}
for each $t\in[0,t_1]$. Combining $(\ref{3.52})$ and $(\ref{3.53})$, we get (\ref{3.3}). This completes the proof of Proposition 3.2.

Hence, by Propositions 3.1 and 3.2, we finally arrive at the conclusion of Theorem 1.2.

\bigbreak
\begin{center}
{\bf Acknowledgment}
\end{center}
This work was supported by a grant
from the National Natural Science Foundation of China under
contract 10925103 and ``the Fundamental Research Funds for the Central Universities".

\end{document}